\documentclass[12pt]{article}
\usepackage{amsthm, amsmath, amssymb, enumerate}
\usepackage{fullpage}
\usepackage{amscd}
\usepackage[colorlinks=true]{hyperref}
\usepackage{MnSymbol}

\begin{document}
\title{Modular Heights of Quaternionic Shimura Curves}
\author{Xinyi Yuan}
\maketitle

\theoremstyle{plain}

\newtheorem{thm}{Theorem}[section]
\newtheorem{theorem}[thm]{Theorem}
\newtheorem{cor}[thm]{Corollary}
\newtheorem{corollary}[thm]{Corollary}
\newtheorem{lem}[thm]{Lemma}
\newtheorem{lemma}[thm]{Lemma}
\newtheorem{pro}[thm]{Proposition}
\newtheorem{proposition}[thm]{Proposition}
\newtheorem{prop}[thm]{Proposition}
\newtheorem{definition}[thm]{Definition}
\newtheorem{assumption}[thm]{Assumption}

\newtheorem*{thmm}{Theorem}
\newtheorem*{conj}{Conjecture}
\newtheorem*{notation}{Notation}
\newtheorem*{corr}{Corollary}

\theoremstyle{remark} 
\newtheorem{remark}[thm]{Remark}
\newtheorem{example}[thm]{Example}
\newtheorem{remarks}[thm]{Remarks}
\newtheorem{problem}[thm]{Problem}
\newtheorem{exercise}[thm]{Exercise}
\newtheorem{situation}[thm]{Situation}

\numberwithin{equation}{subsection}

\newcommand{\ZZ}{\mathbb{Z}}
\newcommand{\CC}{\mathbb{C}}
\newcommand{\QQ}{\mathbb{Q}}
\newcommand{\RR}{\mathbb{R}}
\renewcommand{\AA}{\mathbb{A}}

\newcommand{\HH}{\mathcal{H}}     

\newcommand{\inverse}{^{-1}}          

\newcommand{\gl}{\mathrm{GL}_2}             
\newcommand{\sll}{\mathrm{SL}_2}            
\newcommand{\adele}{\mathbb{A}}             
\newcommand{\finiteadele}{\mathbb{A}_f}     
\newcommand{\af}{\mathbb{A}_f}              
\newcommand{\across}{\mathbb{A}^{\times}}      
\newcommand{\afcross}{\mathbb{A}_f^{\times}}   

\newcommand{\gla}{\mathrm{GL}_2(\mathbb{A})}       
\newcommand{\glaf}{\mathrm{GL}_2(\mathbb{A}_f)}    
\newcommand{\glf}{\mathrm{GL}_2(F)}                
\newcommand{\glfv}{\mathrm{GL}_2(F_v)}             
\newcommand{\glofv}{\mathrm{GL}_2(O_{F_v})}        

\newcommand{\ofv}{O_{F_v}}                         
\newcommand{\oev}{O_{E_v}}                         
\newcommand{\evcross}{E_v^{\times}}                
\newcommand{\fvcross}{F_v^{\times}}                

\newcommand{\adelef}{\mathbb{A}_F}             
\newcommand{\adelee}{\mathbb{A}_E}             
\newcommand{\aecross}{\mathbb{A}_E^{\times}}   

\newcommand{\OCS}{\overline{\mathcal{S}}}
\renewcommand{\Pr}{\mathcal{P}r}

\newcommand{\jv}{\mathfrak{j}_v}
\newcommand{\fg}{\mathfrak{g}}
\newcommand{\fj}{\mathfrak{j}}

\newcommand{\vv}{\mathbb{V}}                      
\newcommand{\bb}{\mathbb{B}}                      
\newcommand{\bbf}{\mathbb{B}_f}                   
\newcommand{\bfcross}{\mathbb{B}_f^\times}        
\newcommand{\ba}{\mathbb{B}_{\mathbb{A}}}          
\newcommand{\baf}{\mathbb{B}_{\mathbb{A}_f}}       
\newcommand{\bv}{{\mathbb{B}_v}}                   
\newcommand{\bacross}{\mathbb{B}_{\mathbb{A}}^{\times}}      
\newcommand{\bafcross}{\mathbb{B}_{\mathbb{A}_f}^{\times}}   
\newcommand{\bvcross}{\mathbb{B}_v^{\times}}               

\newcommand{\ad}{\mathrm{ad}}            
\newcommand{\NT}{\mathrm{NT}}         
\newcommand{\nonsplit}{\mathrm{nonsplit}}         
\newcommand{\Pet}{\mathrm{Pet}}         
\newcommand{\Fal}{\mathrm{Fal}}         
\newcommand{\Norm}{\mathrm{Norm}}         

\newcommand{\lb}{\mathcal{L}}   
\newcommand{\DD}{\mathcal{D}}   

\newcommand{\quasilim}{\widetilde\lim_{s\rightarrow 0}}   
\newcommand{\pr}{\mathcal{P}r}   

\newcommand{\CMU}{\mathrm{CM}_U}             
\newcommand{\eend}{\mathrm{End}}             
\newcommand{\eendd}{\mathrm{End}^{\circ}}    

\newcommand{\sumu}{\sum_{u \in \mu_U^2 \bs F\cross}}

\newcommand{\supp}{\mathrm{supp}}
\newcommand{\cross}{^{\times}}
\newcommand{\der}{\frac{d}{ds}|_{s=0}}   

\newcommand{\pair}[1]{\langle {#1} \rangle}
\newcommand{\wpair}[1]{\left\{{#1}\right\}}
\newcommand\wh{\widehat}
\newcommand\Spf{\mathrm{Spf}}
\newcommand{\lra}{{\longrightarrow}}

\newcommand{\Ei}{\mathrm{Ei}} 

\newcommand{\sumyu}{\sum_{(y,u)}}

\newcommand{\matrixx}[4]
{\left( \begin{array}{cc}
  #1 &  #2  \\
  #3 &  #4  \\
 \end{array}\right)}        

\newcommand{\barint}{\mbox{$ave \int$}}  
\def\barint_#1{\mathchoice
            {\mathop{\vrule width 6pt
height 3 pt depth -2.5pt
                    \kern -8.8pt
\intop}\nolimits_{#1}}%
            {\mathop{\vrule width 5pt height
3 pt depth -2.6pt
                    \kern -6.5pt
\intop}\nolimits_{#1}}%
            {\mathop{\vrule width 5pt height
3 pt depth -2.6pt
                    \kern -6pt
\intop}\nolimits_{#1}}%
            {\mathop{\vrule width 5pt height
3 pt depth -2.6pt
          \kern -6pt \intop}\nolimits_{#1}}}

\newcommand{\BA}{{\mathbb {A}}}
\newcommand{\BB}{{\mathbb {B}}}
\newcommand{\BC}{{\mathbb {C}}}
\newcommand{\BD}{{\mathbb {D}}}
\newcommand{\BE}{{\mathbb {E}}}
\newcommand{\BF}{{\mathbb {F}}}
\newcommand{\BG}{{\mathbb {G}}}
\newcommand{\BH}{{\mathbb {H}}}
\newcommand{\BI}{{\mathbb {I}}}
\newcommand{\BJ}{{\mathbb {J}}}
\newcommand{\BK}{{\mathbb {K}}}
\newcommand{\BL}{{\mathbb {L}}}
\newcommand{\BM}{{\mathbb {M}}}
\newcommand{\BN}{{\mathbb {N}}}
\newcommand{\BO}{{\mathbb {O}}}
\newcommand{\BP}{{\mathbb {P}}}
\newcommand{\BQ}{{\mathbb {Q}}}
\newcommand{\BR}{{\mathbb {R}}}
\newcommand{\BS}{{\mathbb {S}}}
\newcommand{\BT}{{\mathbb {T}}}
\newcommand{\BU}{{\mathbb {U}}}
\newcommand{\BV}{{\mathbb {V}}}
\newcommand{\BW}{{\mathbb {W}}}
\newcommand{\BX}{{\mathbb {X}}}
\newcommand{\BY}{{\mathbb {Y}}}
\newcommand{\BZ}{{\mathbb {Z}}}

\newcommand{\CA}{{\mathcal {A}}}
\newcommand{\CB}{{\mathcal {B}}}
\renewcommand{\CD}{{\mathcal{D}}}
\newcommand{\CE}{{\mathcal {E}}}
\newcommand{\CF}{{\mathcal {F}}}
\newcommand{\CG}{{\mathcal {G}}}
\newcommand{\CH}{{\mathcal {H}}}
\newcommand{\CI}{{\mathcal {I}}}
\newcommand{\CJ}{{\mathcal {J}}}
\newcommand{\CK}{{\mathcal {K}}}
\newcommand{\CL}{{\mathcal {L}}}
\newcommand{\CM}{{\mathcal {M}}}
\newcommand{\CN}{{\mathcal {N}}}
\newcommand{\CO}{{\mathcal {O}}}
\newcommand{\CP}{{\mathcal {P}}}
\newcommand{\CQ}{{\mathcal {Q}}}
\newcommand{\CR }{{\mathcal {R}}}
\newcommand{\CS}{{\mathcal {S}}}
\newcommand{\CT}{{\mathcal {T}}}
\newcommand{\CU}{{\mathcal {U}}}
\newcommand{\CV}{{\mathcal {V}}}
\newcommand{\CW}{{\mathcal {W}}}
\newcommand{\CX}{{\mathcal {X}}}
\newcommand{\CY}{{\mathcal {Y}}}
\newcommand{\CZ}{{\mathcal {Z}}}

\newcommand{\RA}{{\mathrm {A}}}
\newcommand{\RB}{{\mathrm {B}}}
\newcommand{\RC}{{\mathrm {C}}}
\newcommand{\RD}{{\mathrm {D}}}
\newcommand{\RE}{{\mathrm {E}}}
\newcommand{\RF}{{\mathrm {F}}}
\newcommand{\RG}{{\mathrm {G}}}
\newcommand{\RH}{{\mathrm {H}}}
\newcommand{\RI}{{\mathrm {I}}}
\newcommand{\RJ}{{\mathrm {J}}}
\newcommand{\RK}{{\mathrm {K}}}
\newcommand{\RL}{{\mathrm {L}}}
\newcommand{\RM}{{\mathrm {M}}}
\newcommand{\RN}{{\mathrm {N}}}
\newcommand{\RO}{{\mathrm {O}}}
\newcommand{\RP}{{\mathrm {P}}}
\newcommand{\RQ}{{\mathrm {Q}}}
\newcommand{\RS}{{\mathrm {S}}}
\newcommand{\RT}{{\mathrm {T}}}
\newcommand{\RU}{{\mathrm {U}}}
\newcommand{\RV}{{\mathrm {V}}}
\newcommand{\RW}{{\mathrm {W}}}
\newcommand{\RX}{{\mathrm {X}}}
\newcommand{\RY}{{\mathrm {Y}}}
\newcommand{\RZ}{{\mathrm {Z}}}

\newcommand{\ga}{{\frak a}}
\newcommand{\gb}{{\frak b}}
\newcommand{\gc}{{\frak c}}
\newcommand{\gd}{{\frak d}}
\newcommand{\gf}{{\frak f}}
\newcommand{\gh}{{\frak h}}
\newcommand{\gi}{{\frak i}}
\newcommand{\gj}{{\frak j}}
\newcommand{\gk}{{\frak k}}
\newcommand{\gm}{{\frak m}}
\newcommand{\gn}{{\frak n}}
\newcommand{\go}{{\frak o}}
\newcommand{\gp}{{\frak p}}
\newcommand{\gq}{{\frak q}}
\newcommand{\gr}{{\frak r}}
\newcommand{\gs}{{\frak s}}
\newcommand{\gt}{{\frak t}}
\newcommand{\gu}{{\frak u}}
\newcommand{\gv}{{\frak v}}
\newcommand{\gw}{{\frak w}}
\newcommand{\gx}{{\frak x}}
\newcommand{\gy}{{\frak y}}
\newcommand{\gz}{{\frak z}}

\newcommand{\ab}{{\mathrm{ab}}}
\newcommand{\Ad}{{\mathrm{Ad}}}
\newcommand{\an}{{\mathrm{an}}}
\newcommand{\Aut}{{\mathrm{Aut}}}

\newcommand{\Br}{{\mathrm{Br}}}
\newcommand{\bs}{\backslash}
\newcommand{\bbs}{\|\cdot\|}

\newcommand{\Ch}{{\mathrm{Ch}}}
\newcommand{\cod}{{\mathrm{cod}}}
\newcommand{\cont}{{\mathrm{cont}}}
\newcommand{\cl}{{\mathrm{cl}}}
\newcommand{\criso}{{\mathrm{criso}}}

\newcommand{\dR}{{\mathrm{dR}}}
\newcommand{\disc}{{\mathrm{disc}}}
\newcommand{\Div}{{\mathrm{Div}}}
\renewcommand{\div}{{\mathrm{div}}}

\newcommand{\Eis}{{\mathrm{Eis}}}
\newcommand{\End}{{\mathrm{End}}}

\newcommand{\Frob}{{\mathrm{Frob}}}

\newcommand{\Gal}{{\mathrm{Gal}}}
\newcommand{\GL}{{\mathrm{GL}}}
\newcommand{\GO}{{\mathrm{GO}}}
\newcommand{\GSO}{{\mathrm{GSO}}}
\newcommand{\GSp}{{\mathrm{GSp}}}
\newcommand{\GSpin}{{\mathrm{GSpin}}}
\newcommand{\GU}{{\mathrm{GU}}}
\newcommand{\BGU}{{\mathbb{GU}}}

\newcommand{\Hom}{{\mathrm{Hom}}}
\newcommand{\Hol}{{\mathrm{Hol}}}
\newcommand{\HC}{{\mathrm{HC}}}

\renewcommand{\Im}{{\mathrm{Im}}}
\newcommand{\Ind}{{\mathrm{Ind}}}
\newcommand{\inv}{{\mathrm{inv}}}
\newcommand{\Isom}{{\mathrm{Isom}}}

\newcommand{\Jac}{{\mathrm{Jac}}}
\newcommand{\JL}{{\mathrm{JL}}}

\newcommand{\Ker}{{\mathrm{Ker}}}
\newcommand{\KS}{{\mathrm{KS}}}

\newcommand{\Lie}{{\mathrm{Lie}}}

\newcommand{\new}{{\mathrm{new}}}
\newcommand{\NS}{{\mathrm{NS}}}

\newcommand{\ord}{{\mathrm{ord}}}
\newcommand{\ol}{\overline}

\newcommand{\rank}{{\mathrm{rank}}}

\newcommand{\PGL}{{\mathrm{PGL}}}
\newcommand{\PSL}{{\mathrm{PSL}}}
\newcommand{\Pic}{\mathrm{Pic}}
\newcommand{\Prep}{\mathrm{Prep}}
\newcommand{\Proj}{\mathrm{Proj}}

\renewcommand{\Re}{{\mathrm{Re}}}
\newcommand{\red}{{\mathrm{red}}}
\newcommand{\sm}{{\mathrm{sm}}}
\newcommand{\sing}{{\mathrm{sing}}}
\newcommand{\reg}{{\mathrm{reg}}}
\newcommand{\Rep}{{\mathrm{Rep}}}
\newcommand{\Res}{{\mathrm{Res}}}

\newcommand{\Sel}{{\mathrm{Sel}}}
\font\cyr=wncyr10  \newcommand{\Sha}{\hbox{\cyr X}}
\newcommand{\SL}{{\mathrm{SL}}}
\newcommand{\SO}{{\mathrm{SO}}}
\newcommand{\Sp}{\mathrm{Sp}}
\newcommand{\Spec}{{\mathrm{Spec}}}
\newcommand{\Sym}{{\mathrm{Sym}}}
\newcommand{\sgn}{{\mathrm{sgn}}}
\newcommand{\Supp}{{\mathrm{Supp}}}

\newcommand{\tor}{{\mathrm{tor}}}
\newcommand{\tr}{{\mathrm{tr}}}

\newcommand{\ur}{{\mathrm{ur}}}

\newcommand{\vol}{{\mathrm{vol}}}

\newcommand{\wt}{\widetilde}
\newcommand{\pp}{\frac{\partial\bar\partial}{\pi i}}
\newcommand{\intn}[1]{\left( {#1} \right)}
\newcommand{\norm}[1]{\|{#1}\|}
\newcommand{\sfrac}[2]{\left( \frac {#1}{#2}\right)}
\newcommand{\ds}{\displaystyle}
\newcommand{\ov}{\overline}
\newcommand{\incl}{\hookrightarrow}
\newcommand{\imp}{\Longrightarrow}
\newcommand{\lto}{\longmapsto}
\newcommand{\iso}{\overset \sim \lra}

\tableofcontents


\section{Introduction}

The goal of this paper is to prove a formula expressing the modular height of a quaternionic Shimura curve over a totally real number field in terms of the logarithmic derivative of the Dedekind zeta function of the totally real number field. 
Our proof is based on the work Yuan--Zhang--Zhang \cite{YZZ} on the Gross--Zagier formula, and the work Yuan--Zhang \cite{YZ} on the averaged Colmez conjecture. All these works are in turn inspired by the Pioneering work Gross--Zagier \cite{GZ} and some philosophies of Kudla's program.

In the following, let us state the exact formula, compare it with other similar formulas, and  explain our idea of proof.

\subsection{Modular height of the Shimura curve} \label{sec intro shimura curve}

Let $F$ be a totally real number field. 
Let $\Sigma$ be a finite set of places of $F$ containing all the archimedean places and having an odd cardinality $|\Sigma|$. 
Denote by $\Sigma_f$ the subset of non-archimedean places in $\Sigma$.
Let $\BB$ be the totally definite incoherent quaternion algebra over the adele ring
$\BA=\BA_F$ with ramification set $\Sigma$. 

Let $U=\prod_{v\nmid\infty}U_v$ be an open compact subgroup of $\BB_f^\times$ such that $U_v$ is maximal at every $v\in \Sigma_f$. 
Our main theorem concerns the case that $U$ is maximal, but we allow $U$ to be more general before the statement of the main theorem.
Let $X_U$ be the associated \emph{Shimura curve} over $F$, which is a projective and smooth curve over $F$ descended from the analytic quotient
$$
X_{U,\sigma}(\BC)=(B(\sigma)^\times \bs \CH^\pm\times \BB_f^\times/U) \cup \{\rm cusps\},
$$
where $\sigma:F\to\BC$ is any archimedean place of $F$ and $B(\sigma)$ is the quaternion algebra over $F$ with ramification set $\Sigma\setminus \{\sigma\}$.
Note that $X_U$ is defined as the corresponding coarse moduli scheme, which is a projective and smooth curve over $F$. 
See \cite[\S1.2.1]{YZZ} for more details. 

Let $L_U$ be the \emph{Hodge bundle} of $X_U$ corresponding to modular forms of weight 2. It is a $\QQ$-line bundle over $X_U$, i.e. an element of 
$\Pic(X_U)\otimes_\ZZ\QQ$, defined by
$$
L_{U}= \omega_{X_{U}/F} \otimes\CO_{X_U}\Big(\sum_{Q\in X_U(\overline F)} (1-e_Q^{-1}) Q\Big).
$$
Here $\omega_{X_{U}/F}$ is the canonical bundle of $X_U$ over $F$, and
for each $Q\in X_U(\overline F)$, the ramification index $e_Q$ is described as
follows.
Fix an embedding $\bar\sigma:\overline F\to \CC$ extending $\sigma:F\to \CC$, so that $Q$ is also viewed as a point of $X_{U,\sigma}(\BC)$.
If $Q$ is a cusp in $X_{U,\sigma}(\BC)$ under the above complex uniformization, then $e_Q=\infty$ and $1-e_Q^{-1}=1$.
If $Q$ is not a cusp in that sense, then the connected component of $Q$ in $X_{U,\sigma}(\BC)$ can be written as a quotient $\Gamma\bs\CH^*$ for a discrete group $\Gamma$, and  $e_Q$ is the ramification index of any preimage of $Q$ under the map 
$\CH\to \Gamma\bs\CH$. 
One can check that $e_Q$ does not depend on the choices of $(\sigma, \bar\sigma)$ or the preimage, and that 
$e_Q$ is Galois invariant, so $L_U$ is indeed defined over $F$. 
See \cite[\S3.1.3]{YZZ} for more details.

Let $\CX_U$ be the \emph{canonical integral model} of $X_U$ over $O_F$, as reviewed in \cite[\S4.2]{YZ}. Note that we always assume that $U$ is maximal at every $v\in \Sigma_f$. 
If $|\Sigma|=1$ (or equivalently $F=\QQ$ and $\Sigma=\{\infty\}$), then $X_U$ is a modular curve, $X_U\simeq \BP^1_\QQ$ via the $j$-function, and 
$\CX_U\simeq\BP^1_\ZZ$ under this identification. 
We refer to Deligne--Rapoport \cite{DR} for a thorough theory of this situation.
If $|\Sigma|>1$, $\CX_U$ is a projective, flat, normal, semistable and $\QQ$-factorial arithmetic surface over $O_F$, defined as quotients of the canonical integral models of 
 $\CX_{U'}$ for sufficiently small open compact subgroups $U'$ of $U$.
We refer to Carayol \cite{Ca} and Boutot--Zink \cite{BZ} for integral models for sufficiently small level groups, and refer to \S\ref{sec shimura curve} for the quotient process.
Only the part of the integral model $\CX_U$ above places $v$ such that $U_v$ is maximal is essential in this paper.

Let $\CL_U$ be the \emph{canonical integral model} of $L_U$ over $\CX_U$ as reviewed in \cite[\S4.2]{YZ}. 
Then $\CL_U$ is a $\QQ$-line bundle on $\CX_U$ constructed as follows. If $U$ is sufficiently small in the sense that $U$ is contained in $U(N)=(1+N O_{\bb_f})^\times$ for some integer $N\geq 3$ and some maximal order $O_{\bb_f}$ of $\bb_f$, then over the open subscheme 
$\CX_{U, S(U)}=\CX_{U}\times_{\Spec\, O_F} S(U)$ with
$$S(U)=\Spec\, O_F\setminus \{v: U_v \text{ is not maximal}\},$$ 
we have 
$$
\CL_U|_{\CX_{U, S(U)}}= \omega_{\CX_{U}/O_F} \otimes \CO_{\CX_{U, S(U)}}\Big(\sum_{Q\in X_U} (1-e_Q^{-1}) \CQ\Big).
$$
Here $\omega_{\CX_{U}/O_F}$ is the relative dualizing sheaf, the summation is through closed points $Q$ of $X_U$, $\CQ$ is the Zariski closure of $Q$ in $\CX_{U,S(U)}$, and $e_Q$ is the ramification index of any point of $X_U(\ol F)$ corresponding to $Q$.
If $U$ is a maximal open compact subgroup of $\bfcross$, for any sufficiently small normal open compact subgroup $U'$ of $\bfcross$ contained in $U$, we have
$$
\CX_{U, S(U')}= \CX_{U',S(U')}/(U/U'), \quad
\CL_{U}|_{\CX_{U, S(U')}}= N_{\pi} (\CL_{U'}|_{\CX_{U', S(U')}})^{\otimes (1/\deg(\pi))},
$$
where
$N_{\pi}:\Pic(\CX_{U', S(U')})\to \Pic(\CX_{U, S(U')})$ is the norm map with respect to the natural map $\pi: \CX_{U', S(U')}\to \CX_{U, S(U')}$. 
Varying $U'$, we glue $\{\CL_{U}|_{\CX_{U, S(U')}}\}_{U'}$ together to form the $\QQ$-line bundle $\CL_U$ over $\CX_U$ for maximal $U$.
For general $U$, we take an embedding $U\subset U_0$ into a maximal $U_0$, and than define $\CL_U$ to be the pull-back of $\CL_{U_0}$ via the natural map $\CX_U \to \CX_{U_0}$. 

At any archimedean place $\sigma:F\to \BC$, 
 the \emph{Petersson metric} of $\CL_U$ is given by 
$$\|f(\tau)d\tau\|_{\mathrm{Pet}}=2\, \Im(\tau) |f(\tau)|,$$
where $\tau$ is the standard coordinate function on $\CH\subset \CC$, and $f(\tau)$ is any meromorphic modular form of weight 2 over $X_{U,\sigma}(\CC)$. 
Thus we have \emph{the arithmetic Hodge bundle}
$$\ol\CL_U=(\CL_U, \{\|\cdot\|_\sigma\}_\sigma).$$

The \emph{modular height} of $X_U$ with respect to the arithmetic Hodge bundle $\ol\CL_U$ is defined to be
$$
h_{\ol\CL_U}(X_U)
= \frac{\widehat\deg(\hat c_1(\ol\CL_U)^2)}{2\deg(L_{U})}.
$$
Here $\deg(L_{U})$ is the degree over the generic fiber $X_U$, and the
numerator is the arithmetic self-intersection number on the arithmetic surface 
$\CX_U$ in the setting of Arakelov geometry. 
Note that if $|\Sigma|>1$, then $\ol\CL_U$ is a hermitian $\QQ$-line bundle over $\CX_U$, and the self-intersection number essentially follows from the theory of Gillet--Soul\'e \cite{GS};
if $|\Sigma|=1$, then the metric has a logarithmic singularity along the cusp, and the intersection number is defined in the framework of Bost \cite{Bo} or K\"uhn \cite{Kuh}. 

By the projection formula, $h_{\ol\CL_U}(X_U)$ is independent of $U$. 
However, as $\CL_U$ is less canonical at places $v$ such that $U_v$ is not maximal,  we will assume that $U$ is maximal at every $v$ in the main theorem and afterwards.

For any non-archimedean place $v$ of $F$, denote by $N_v$ the norm of $v$. 
Recall the Dedekind zeta function
$$
\zeta_{F}(s)=\prod_{v\nmid\infty} (1-N_v^{-s})^{-1}.
$$
The functional equation switches the values and derivatives of $\zeta_{F}(s)$ between $-1$ and $2$.
The goal of this paper is to prove the following formula. 

\begin{thm}[modular height]\label{main}
Let $U\subset \BB_f^\times$ be a \emph{maximal} open compact subgroup. 
Then
\begin{eqnarray*}
h_{\ol \CL_U}(X_U)
=-\frac{\zeta_{F}'(-1)}{\zeta_{F}(-1)} -\frac12[F:\QQ] 
+\sum_{v\in\Sigma_f}  \frac{3N_v-1}{4(N_v-1)}\log N_v.
\end{eqnarray*}
\end{thm}

If $F=\QQ$ and $\Sigma=\{\infty\}$, the formula was proved by Bost (un-published) and K\"uhn (cf. \cite[Theorem 6.1]{Kuh}); if $F=\QQ$ and $|\Sigma|>1$, the formula was proved by Kudla--Rapoport--Yang (cf. \cite[Theorem 1.0.5]{KRY2}).

Denote by $h_F$ the class number of $O_F$. 
A classical formula of Vign\'eras \cite{Vi} gives  
\begin{eqnarray*}
\deg(L_{U})
= 4\cdot h_F \cdot (-2)^{-[F:\QQ]} \cdot \zeta_{F}(-1) \cdot
\prod_{v\in \Sigma_f}(N_v-1).
\end{eqnarray*}
This is also an easy consequence of the formula in the remark right after \cite[Proposition 4.2]{YZZ}.
Theorem \ref{main} is an arithmetic version of this formula. It computes the arithmetic degree instead of the geometric degree, and the result is given by the logarithmic derivative at $-1$ instead of the value at $-1$. 

The relation between these two formulas is similar to the relation between the Gross--Zagier formula and the Waldspurger formula (as fully explored in \cite{YZZ}), and is also similar to the relation between the averaged Colmez conjecture and the class number formula (as treated in \cite{YZZ}). 

In Kudla's program, it is crucial to extend the (modular) generating series of CM cycles over a Shimura variety to a (modular) generating series of arithmetic cycles over a reasonable integral model. An idea of S. Zhang \cite[\S3.5]{Zh2} to treat this problem is to apply his notion of admissible arithmetic extensions. This approach relies on concrete results on arithmetic intersection numbers, so our main formula fits  this setting naturally.
Inspired by S. Zhang's idea, Qiu \cite{Qi} solved the problem for generating series of divisors over unitary Shimura varieties under some assumptions,  and his argument is based on many computational results of this paper.

\subsection{The case $F=\QQ$ and other similar formulas}

If $F=\BQ$ and $\Sigma=\{\infty\}$, or equivalently if $X_U$ is the usual modular curve, then the formula of Bost and K\"uhn \cite[Theorem 6.1]{Kuh} agrees with our formula by \cite[Theorem 5.3, Remark 5.4]{Yu2}.

If $F=\BQ$ and $|\Sigma|>1$, the formula in \cite[Theorem 1.0.5]{KRY2} of Kudla--Rapoport--Yang is equivalent to 
\begin{eqnarray*}
 h_{\widehat\omega_0}(X_U) = -\frac{\zeta_{\BQ}'(-1)}{\zeta_{\BQ}(-1)}
-\frac12 + \sum_{p\in\Sigma_f}  \frac{p+1}{4(p-1)}\log p.
\end{eqnarray*}
This formula is compatible with our formula. 
In fact, the right-hand side of the formula differs from that of ours by $\ds\frac12 \log d_\BB$, and $\ds h_{\overline\CL_U}(X_U)= h_{\widehat\omega_0}(X_U)+\frac12 \log d_\BB$  by
the explicit results on the Kodaira--Spencer map in
 \cite[Theorem 2.2, Remark 2.3]{Yu2}. 

There are many formulas of similar flavor in the literature. 
Besides the above mentioned works of Bost, K\"uhn and Kudla--Rapoport--Yang, 
Bruinier--Burgos--K\"uhn \cite{BBK} proved a modular height formula for Hilbert modular surfaces, 
H\"ormann \cite{Ho} proved a modular height formula up to $\log \QQ_{>0}$ 
for Shimura varieties of orthogonal types over $\QQ$, 
and Bruinier--Howard \cite{BH} recently proved a modular height formula 
for Shimura varieties of unitary types over $\QQ$.
The formulas of \cite{Ho, BH} are based on the formulas of Bost, K\"uhn and 
Kudla--Rapoport--Yang.

In a slightly different direction,
Freixas--Sankaran \cite{FS} proved some other formulas for intersections of more general Chern classes over Hilbert modular surfaces. 
Finally, we refer to Maillot--R\"ossler \cite{MR1,MR2} for far-reaching conjectures generalizing these formulas.

Our formula is primitive in that it involves Dedekind zeta functions of general totally real fields, while the previous known formulas involve Dedekind zeta functions of $\QQ$ and quadratic fields.
Moreover, based on our formula and a strategy in the flavor of \cite{Ho, BH}, it is promising to prove similar formulas for high-dimensional Shimura varieties associated to orthogonal groups of signature $(n,2), (n+2,0), \cdots, (n+2,0)$ or unitary groups of signature $(n,1), (n+1,0), \cdots, (n+1,0)$ over totally real fields.

\subsection{Modular height of a CM point}

Our proof of Theorem \ref{main} is inspired by the works \cite{YZZ,YZ}. 
In the proof, we need to pick an auxiliary CM point, and the height of this point is also relevant to our treatment. 
Let us first review a formula in \cite{YZ} which is related to our main theorem. 

Let $E$ be a totally imaginary quadratic extension over $F$. 
Assume that there is an embedding $\BA_E\hookrightarrow \bb$ of $\BA$-algebras such that the image of $\wh O_E^\times$ lies in the maximal compact subgroup $U$.

Let $P_U\in X_{U,\sigma}(\CC)$ be the CM point represented by $[\tau_0,1]$ under the complex uniformization at $\sigma:F\to \CC$, where $\tau_0$ is the unique fixed point of $E^\times$ in $\CH$.
Fix embeddings $E^\ab\to \overline F\to \CC$ compatible with $\sigma:F\to \CC$.
By the definition of the canonical model, we actually have $P_U\in X_U(E^\ab)$ via the embedding.
The modular height of $P_U$ is defined by
$$h_{\ol\CL_U}(P_U):=\frac{1}{\deg(P_U)}\wh\deg(\ol\CL_U|_{\bar P_U}),$$
where $\bar P_U$ denotes the Zariski closure of the image of $P_U$ in $\CX_U$, and 
$\deg(P_U)$ is the degree of the field of definition of of $P_U$ over $F$. 
By \cite[Theorem 1.7]{YZ}, we have the following formula.

\begin{thm}\label{height CM}
Assume that  
there is no non-archimedean place of $F$ ramified in both $E$ and $\BB$. Then
$$ h_{\ol\CL_U}(P_U)
=-\frac{L_f'(0,\eta)}{L_f(0,\eta)} +\frac 12 \log  \frac{d_\BB}{d_{E/F}}.
$$
Here $\eta:F^\times\bs\AA_F^\times\to \{\pm1\}$ is the quadratic character associated to the quadratic extension $E/F$, 
$d_\BB=\prod_{v\in \Sigma_f}N_v$ is the absolute discriminant of $\BB$, and $d_{E/F}$ is the norm of the relative discriminant of ${E/F}$.
\end{thm}

Theorem \ref{height CM} is one of the two steps in the proof of the averaged Colmez conjecture of \cite{YZ}. The averaged Colmez conjecture was proved independently by Andreatta--Goren--Howard--Madapusi-Pera \cite{AGHM}, and plays a crucial role in the final solution of the Andr\'e--Oort conjecture of Tsimerman \cite{Ts}. 

In the case $F=\QQ$ and $\Sigma=\{\infty\}$, Theorem \ref{height CM} is equivalent to the classical Chowla--Selberg formula proved in \cite{CS}. We refer to \cite[\S3.3]{Yu1} for many equivalent forms of the Chowla--Selberg formula.

\subsection{Kronecker's limit formula}

Both the Bost--K\"uhn formula and the Chowla--Selberg formula are easy consequences of the more classical Kronecker limit formula. 

In fact, by \cite[Prop. 5.2]{Kuh}, the Kronecker limit formula asserts that
$$
-\log |\Delta(\tau)^2 \mathrm{Im}(\tau)^{12}|
=4\pi \lim_{s\to 1} (E(\tau, s)-\varphi(s)),
$$
where 
$$
\Delta(\tau)= q\prod_{n=1}^\infty (1-q^n)^{24}, \quad q=e^{2\pi i\tau}
$$
is the modular discriminant function,  
$$E(\tau, s)=\frac12\sum_{c,d\in \ZZ,\ \gcd(c,d)=1} 
\frac{\mathrm{Im}(\tau)^s}{|c\tau+d|^{2s}}
$$
is the classical non-holomorphic Eisenstein series, and 
$$
\frac{\pi}{3}\varphi(s)=\frac{1}{s-1}+ 2-2\log(4\pi)-24\zeta_{\QQ}'(-1)+O(s-1).
$$

In particular, $\Delta(\tau)$ induces a global section of $L_U^{\otimes 6}$ over the modular curve $X_U=X_0(1)$. Then we can use this section to compute 
$h_{\ol\CL_U}(X_U)$ and $h_{\ol\CL_U}(P_U)$. 
 
Integrating $-\log |\Delta(\tau)^2 \mathrm{Im}(\tau)^{12}|$
over $X_U(\BC)$ with respect to the Poincare measure $y^{-2}dxdy$, the Kronecker limit formula implies the Bost--K\"uhn formula. This is essentially the proof of K\"uhn \cite{Kuh}. 

Averaging $-\log |\Delta(\tau)^2 \mathrm{Im}(\tau)^{12}|$ over the Galois orbit of the CM point $P_U$, the Kronecker limit formula implies the Chowla--Selberg formula. This is essentially the proof in Weil \cite{We}. 

In summary, in the case $|\Sigma|=1$, both Theorem \ref{main} and Theorem \ref{height CM}
are consequences of the Kronecker limit formula. 

On the other hand, there is no analogous formulation of the Kronecker limit formula over totally real fields, since there is no explicit modular form over a quaternionic Shimura curve to replace the classical modular discriminant function $\Delta$. 
Hence, the above proof of the theorem does not work in the general case. 

Our proofs of Theorem \ref{main} and Theorem \ref{height CM} are extensions of the treatment of \cite{YZZ}. The original goal of \cite{YZZ} is to prove the Gross--Zagier formula over Shimura curves, but the method was enhanced in \cite{YZ} to prove Theorem \ref{height CM}, and now we can further enhance the method to prove
Theorem \ref{main}. Note that our proof of Theorem \ref{main} in the case $F=\QQ$ is  different from those of \cite{Kuh, KRY2}. 

It is interesting that in both the classical proofs and our current proofs, 
Theorem \ref{main} and Theorem \ref{height CM} are always put in the same framework.

\subsection{Idea of proof}

Now we sketch our proof of Theorem \ref{main}.
It is an extension of the proof of the Gross--Zagier formula in \cite{YZZ} and the proof of the averaged Colmez conjecture in \cite{YZ}. 
To have a setup compatible with those in \cite{YZZ,YZ}, we first choose a CM extension $E$ over $F$ as in Theorem \ref{height CM}, though $E$ is irrelevant to the final statement of Theorem \ref{main}. 

\subsubsection*{The degeneracy assumptions}

Recall that the Gross--Zagier formula is an identity between the derivative of the Rankin--Selberg $L$-function of 
a Hilbert modular form  and the height of a CM point on a modular abelian variety. This formula is proved 
by a comparison  of a derivative  series  $\Pr I'(0, g, \phi)$ with a geometric series $2Z(g, (1,1), \phi)$ parametrized by certain modified 
Schwartz function $\phi\in \ol\CS (\BB\times \BA^\times )$. More precisely, we have proved that  
the difference 
$$\mathcal D (g, \phi)=\Pr I'(0, g, \phi)-2Z(g, (1,1), \phi), \qquad g\in \GL_2(\BA_F)$$
is perpendicular to the relevant cusp form.

The matching for the ``main terms'' of $\mathcal D (g, \phi)$ eventually implies the Gross--Zagier formula in \cite{YZZ}. In this process, many assumptions on the choice of $\phi$ 
in \cite[\S5.2.1]{YZZ} are made to ``annihilate'' the ``degenerate terms'', which simplifies the calculations dramatically and forces the computational results to satisfy the conditions of an approximation argument. The ``strictest'' degeneracy assumptions involved are 
\cite[Assumption 5.3, Assumption 5.4]{YZZ}.
The assumptions are not harmful for the Gross--Zagier formula, as proved in
\cite[Theorem 5.7]{YZZ}.

Nonetheless, if we allow the Schwartz function to be more general, the matching process will actually give us more formulas. In fact, after removing \cite[Assumption 5.3]{YZZ}, we obtain a matching of some ``degenerate terms'', which eventually implies Theorem \ref{height CM}. This is the work of \cite[Part II]{YZ}. 

In the current paper, we remove both \cite[Assumption 5.3, Assumption 5.4]{YZZ} when considering the matching of the series $\Pr I'(0, g, \phi)$ and $2Z(g, (1,1), \phi)$. Then we finally obtain an extra identity, which eventually implies Theorem \ref{main}. 
Our precise choice of the Schwartz functions is given in \S\ref{choices}.

From \cite{YZZ} to \cite{YZ}, and from \cite{YZ} to the current paper, each step removes a degeneracy assumption, which causes two significant problems. The first problem is that more terms appear in the comparison, which incur far more involved local computations. 
This is eventually overcome by patience and carefulness. 
The second problem is how to obtain exact identity from the ``partial matching'' of the two series; i.e., the matching of ``all but finitely many'' terms of the two series. 
In \cite{YZZ}, this problem is solved by the method of approximation (cf. \cite[\S1.5.10]{YZZ}). 
In \cite{YZ}, this problem is solved by the theory of pseudo-theta series (cf. \cite[\S6]{YZ}), which is an extension of the method of approximation.
In the current paper, the theory of pseudo-theta series is not sufficient for the comparison. Our solution is to introduce a new notion of \emph{pseudo-Eisenstein series}, and generalize 
\cite[Lemma 6.1]{YZ}, the key matching principle of pseudo-theta series, to include both pseudo-theta series and pseudo-Eisenstein series.

In the following, we review the derivative series $\Pr I'(0, g, \phi)$ and the height series $Z(g, (1,1), \phi)$ and introduce some new ingredients of our proof.

\subsubsection*{Derivative series}
By the reduced norm $q$, the incoherent quaternion algebra $\BB$ is viewed as a quadratic space over $\BA=\BA_F$. 
Then we have a modified space $\ol\CS (\BB\times \BA^\times)$ of Schwartz functions, which has a Weil representation $r$ by $\GL_2(\BA)\times\BB^\times\times \BB^\times$.
Strictly speaking, the representation $r$ is induced by the canonical homomorphism $\BB^\times\times \BB^\times\to \GO(\BB,q)$ defined by sending $(b_1,b_2)\in \BB^\times\times \BB^\times$ to the automorphism of $\BB$ given by $x\mapsto b_1xb_2^{-1}$. 
We refer to \S\ref{sec theta eisenstein} (or the original \cite[\S 2.1, \S2.2, \S 4.1]{YZZ}) for more details. 

For each  $\phi\in \ol\CS (\BB\times \BA^\times )$ invariant under an open compact subgroup $U\times U$ of $\BB_f^\times\times \BB_f^\times$, 
we have a mixed theta--Eisenstein series
$$I(s, g, \phi)=\sum _{u\in \mu _U^2\bs F^\times} \sum _{\gamma \in P^1(F)\bs \SL_2(F)}
\delta (\gamma g)^s \sum _{x_1\in E} r(\gamma g)\phi (x_1, u),$$
where $\mu _U=F^\times \cap U$, and $P^1$ is the upper triangular  subgroup of $\SL_2$.

The derivative $I'(0, g, \phi)$ of $I(s, g, \phi)$ at $s=0$ is an automorphic form in $g\in \gla$. 
Let $\Pr I'(0, g, \phi)$ be the \emph{holomorphic projection} of the derivative $I'(0, g, \phi)$. 
This holomorphic projection is just the orthogonal projection from the space of automorphic forms to the space of cuspidal and holomorphic automorphic forms of parallel weight two with respect to the Petersson inner product. 

In \S\ref{sec derivative series}, we decompose $\Pr I'(0, g, \phi)$ into a sum of ``local terms'', and compute all the relevant local components. 
Most of the terms are computed in \cite{YZZ, YZ}. 
However, as in \S\ref{sec 3.1}, 
a new extra term $\Pr' \CJ'(0, g, \phi)$ appears in the expression of $\Pr I'(0, g, \phi)$.
This term comes from the overly fast growth of $I'(0, g, \phi)$ in the computation of the holomorphic projection. It was zero under 
\cite[Assumption 5.4]{YZZ}, but its non-vanishing is crucial to the treatment here. 
The extra term $\Pr' \CJ'(0, g, \phi)$ is 
computed in Proposition \ref{analytic series extra}, and its local component computed in 
Lemma \ref{local explicit}(1) gives $\zeta_v'(2)/\zeta_v(2)$ at almost all places $v$.
The sum over all places gives the global logarithmic derivative $\zeta_F'(2)/\zeta_F(2)$, which is the main term on the right-hand side of Theorem \ref{main}.

\subsubsection*{Height series}
For any $\phi\in \ol\CS (\BB\times \BA^\times)$ invariant under $U\times U$, we have a generating series of Hecke operators on the Shimura curve $X_U$:
$$Z(g, \phi)_U=Z_0(g, \phi)+w_U\sum _{a\in F^\times}\sum _{x\in U\bs \BB_f^\times /U}
r(g)\phi (x, aq(x)^{-1})Z(x)_U,$$
where $w_U=|\{\pm1\}\cap U|$ and every $Z(x)_U$ is a divisor of $X_U\times X_U$ associated to the Hecke operator corresponding to the double coset $UxU$. 
The constant term $Z_0(g, \phi)$ does not play any essential role in this paper, and we denote by $Z_*(g, \phi)$ the sum of the other terms.
By \cite[Theorem 3.17]{YZZ}, this series is absolutely convergent and defines an 
automorphic form in $g\in \GL_2(\BA)$ with coefficients in $\Pic (X_U\times X_U)_\BC$.

Recall that $P_U\in X_U(E^\ab)$ is the CM point represented by $[\tau_0,1]$ under the complex uniformization, where $\tau_0$ is the unique fixed point of $E^\times$ in the upper half plane $\CH$. More generally, we have a CM point $t=[\tau_0,t]$ for any 
$t\in E^\times(\af)$. 
Let $t^\circ=t-\xi_t$
be the divisor in $\Pic(X_{U,\overline F})\otimes_\ZZ\QQ$ of degree zero on every connected component. Here the normalized Hodge class $\displaystyle\xi_{t}=\frac{1}{\deg(L_{U,t})}L_{U,t}$, where $L_{U,t}$ is the restriction of the Hodge bundle $L_U$ to the connected component of $X_{U,\overline F}$ containing $t$.
Then we can form a height series 
$$Z(g, (t_1,t_2),\phi)=\pair{Z(g, \phi)_U t_1^\circ, \ t_2^\circ}_\NT,$$
where the right-hand side is the Neron--Tate height pairing.

In \S\ref{sec height series}, we decompose the height series $Z(g, (t_1,t_2),\phi)$ into a sum of ``local terms'', and compute all the relevant local components. 
For simplicity, we might write $Z(g, (t_1,t_2))$ for $Z(g, (t_1,t_2),\phi)$ by suppressing the dependence on $\phi$ in this paper.
The starting point is the decomposition
$$Z(g, (t_1, t_2))
=\pair{Z_*(g,\phi) t_1,  t_2} -\pair{Z_*(g,\phi) \xi_{t_1}, t_2}
+\pair{Z_*(g,\phi)\xi_{t_1}, \xi_{t_2}} -\pair{Z_*(g,\phi) t_1, \xi_{t_2}}.$$
The first term is computed in \cite{YZZ, YZ}. 
The remaining three terms are further computed in \S\ref{sec 4.3}. 
These three terms are zero under 
\cite[Assumption 5.4]{YZZ}, but their non-vanishing is crucial to the treatment here. 
In particular, by Proposition \ref{geometric series extra1}, the term $\pair{Z_*(g,\phi)\xi_{t_1}, \xi_{t_2}}$ is equal to an Eisenstein series times $\pair{\xi_{t_2}, \xi_{t_2}}$, and thus it is an easy multiple of $h_{\ol \CL_U}(X_U)$. 
This gives the main term on the left-hand side of Theorem \ref{main}.

\subsubsection*{Pseudo-Eisenstein series}

The notion of pseudo-Eisenstein series is parallel to that of pseudo-theta series of \cite{YZ}.
To illustrate the idea, we sketch the idea of both notions for $\SL_2$, while those for $\GL_2$, which are the ones we really need, can be introduced similarly.  

Let $(V,q)$ be a quadratic space over a totally real number field $F$, assumed to be even-dimensional for simplicity. Let $\phi\in \CS (V(\BA))$ be a Schwartz function. Then we have an action of $g\in \SL_2(\BA)$ on $\phi$ via the Weil representation. 

Start with the theta series 
$$\theta(g,\phi)=\sum _{x\in V} r(g)\phi(x), \qquad g\in \SL_2(\BA).$$
Let $S$ be a finite set of non-archimedean places of $F$. In 
$r(g)\phi(x)=r(g_S)\phi_S(x) r(g^S)\phi^S(x)$,
if we replace $r(g_S)\phi_S(x)$ by a locally constant function $\phi_S'(g, x)$ of 
 $(g, x)\in \GL_2(F_S)\times V(F_S)$, then we obtain a 
 \textit{pseudo-theta series}
 $$A^{(S)}_{\phi'}(g)=
\sum _{x\in V}\phi_S'(g, x)r(g)\phi^S (x), \qquad g\in \SL_2(\BA).$$
Note that $A^{(S)}_{\phi'}$ is not automorphic in general. 
More general types of pseudo-theta series are introduced in \cite[\S6]{YZ}
and reviewed in \S\ref{sec key lemma}. 

We say that the pseudo-theta series $A^{(S)}_{\phi'}(g)$ is \textit{non-singular}  if $\phi_S'(1, x)$ (for $g=1$) is actually a Schwarz function of $x\in V(F_S)$. 
In this case, we form a true theta series 
 $$\theta_{A^{(S)}}(g)=
\sum _{x\in V}r(g)\phi_S'(1, x)r(g)\phi^S (x), \qquad g\in \SL_2(\BA).$$
It is automorphic and approximates the original series in the sense that $A^{(S)}_{\phi'}(g)=\theta_{A^{(S)}}(g)$ as long as $g_S=1$. 

Now we start with the Siegel--Eisenstein series 
\begin{eqnarray*}
E(s,g, \phi) =\sum_{\gamma \in P^1(F)\bs \SL_2(F)}
\delta (\gamma g)^s r(\gamma g)\phi(0), \quad g\in\SL_2(\BA).
\end{eqnarray*}
The non-constant part of $E(s, g,\phi)$ has a Fourier expansion 
$$
E_*(s, g, \phi)=\sum_{a\in F^\times} W_a(s,g,\phi),
$$
where the Whittaker function is defined by
\begin{eqnarray*}
W_a(s,g,\phi) = \int_{\adele} \delta(wn(b)g)^s \ r(wn(b)g)\phi(0)
\psi(-ab) db, \quad a\in F.
\end{eqnarray*}
We define the local Whittaker functions similarly. 
For our purpose, we only care about the behavior at $s=0$.
Let $S$ be a finite set of non-archimedean places of $F$. 
In 
$W_a(0,g,\phi)=W_{a,S}(0,g,\phi_S)W_a^S(0,g,\phi^S)$,
if we replace $W_{a,S}(0,g,\phi_S)$ by a locally constant function
$B_{a,S}(g)$ of $(a,g)\in  F_S\cross \times \SL_2(F_S)$, then we obtain a 
\textit{pseudo-Eisenstein series}
$$B^{(S)}_{\phi}(g)=\sum_{a\in F^\times} B_{a,S}(g) W_a^S(0,g,\phi^S), \ \quad g\in \SL_2(\adele).$$

Pseudo-Eisenstein series arise naturally in derivatives of Eisenstein series.
In fact, the derivative of $E_*(s, g, \phi)$ at $s=0$ is
$$
E_*'(0, g, \phi)= \sum_{a\in F^\times} \sum_v W_{a,v}'(0,g,\phi)W_a^v(0,g,\phi^v).
$$
For every non-archimedean $v$, the ``$v$-part''
$$\sum_{a\in F^\times} 
W_{a,v}'(0,g,\phi)W_a^v(0,g,\phi^v)$$
is a pseudo-Eisenstein series. 

We say that the pseudo-Eisenstein series $B^{(S)}_{\phi}(g)$ is \textit{non-singular}  if for every $v\in S$, there exist $\phi_v^+\in \CS(V_v^+)$ and $\phi_v^-\in \CS(V_v^-)$ such that 
$$B_{a,v}(1)=W_{a,v}(0,1,\phi_v^+)+ W_{a,v}(0,1,\phi_v^-),
 \quad \forall a\in F_v\cross.$$
Here $\{V_v^+,V_v^-\}$ is the set of (one or two) quadratic spaces over $F_v$ with the same dimension and the same discriminant as $V_v$.
In this case, we form a linear combination of true Eisenstein series 
$$
E_{B}(s,g)=\sum_{\epsilon: S\to \{\pm\}}
E(s,g, \phi_S^{\epsilon}\otimes \phi^S), \qquad g\in \SL_2(\BA).
$$
It approximates the original series in the sense that 
$B_{\phi}^{(S)}(g)$ is equal to the non-constant part of $E_{B}(0,g)$
as long as $g_S=1$. 

The above pair $(\phi_v^+, \phi_v^-)$ (if it exists) is generally not uniquely determined by $B_{a,v}(1)$, but  
the Eisenstein series 
$E_{B}(s,g)$, as a function of $g\in \SL_2(F_v)$ and $s\in \CC$, is uniquely determined by
$B_{\phi}^{(S)}(g)$.
We refer to Lemma \ref{whittaker image1} for more details (in a slightly different setting).

The key result in this pseudo theory is Lemma \ref{pseudo}, as an extension of \cite[Lemma 6.1]{YZ}. It asserts that if an automorphic form is equal to a finite linear combination 
of non-singular pseudo-theta series and non-singular pseudo-Eisenstein series, then it is actually equal to the finite linear combination 
of the corresponding theta series and Eisenstein series.

\subsubsection*{The comparison}

Go back to the difference 
$$\mathcal D (g, \phi)=\Pr I'(0, g, \phi)-2Z(g, (1,1), \phi), \qquad g\in \GL_2(\BA_F).$$
By the computational result of \S\ref{sec derivative series} and \S\ref{sec height series}, we eventually see that $\mathcal D (g, \phi)$ is a finite linear combination of non-singular pseudo-theta series and non-singular pseudo-Eisenstein series. 

By Lemma \ref{pseudo}, $\mathcal D (g, \phi)$ is actually equal to the finite linear combination 
of the corresponding theta series and Eisenstein series. 
Note that $\mathcal D (g, \phi)$ is cuspidal, so the linear combination of the corresponding constant terms is zero. This gives a nontrivial relation involving the major terms of Theorem \ref{main}. 
It suffices to take $g$ to be a specific matrix to make the relation precise. 
Take $g=(g_v)_v\in\gla$ with $g_v=1$ for $v\notin \Sigma_f$ and $g_v=w$ for 
$v\in \Sigma_f$. 
After explicit computation, the nontrivial relation becomes 
$$
d_0  \sumu r(g)\phi(0,u)=0.
$$
Here $d_0$ is the difference of two sides of Theorem \ref{main}. 
This proves the theorem. 

Note that if we take $g=1$, 
then the nontrivial relation becomes $0=0$, since $\phi_v(0,u)=0$ for any $v\in \Sigma_f$ by our choice $\phi_v=1_{O_\bv^\times\times\ofv\cross}$ in \S\ref{choices}. 
As $r(w)\phi_v(0,u)\neq 0$,
 we choose $g_v$ to be $w$ for $v\in \Sigma_f$ instead. This serves the purpose, but incurs more computations about evaluating $g_v=w$ and about averaging of many local terms.

\subsection{Notations and conventions} \label{sec notation}
Most of the notations of this paper are compatible with those in \cite{YZZ, YZ}. The basic notations are as in \cite[\S1.6]{YZZ}. 

In particular, $F$ is a fixed totally real number field, and $E$ is a fixed totally imaginary quadratic extension of $F$. 
Denote by $\BA=\BA_F=\prod_v' F_v$ the adele ring of $F$. 
As in  \cite[\S1.6]{YZZ}, we normalize the character $\psi=\oplus_v\psi_v:F\bs \BA\to \BC^\times$. 
We introduce the Weil representation based on $\psi$, and choose a precise Haar measure on each relevant algebraic group locally everywhere.

For a non-archimedean place $v$ of $F$, make the following notations. 
\begin{enumerate}[(1)]
\item $p_v$ denotes the maximal ideal of $O_{F_v}$;
\item $N_v$ denotes the order of the residue field $O_{F_v}/p_v$;  
\item $d_v\in F_v$ denotes the local different of $F$ over $\BQ$;
\item $D_v\in F_v$ denotes the discriminant of the quadratic extension $E_v$ in $F_v$.
\end{enumerate}
Note that $d_v$ and $D_v$ are only well-defined up to multiplication by $O_{F_v}^\times$, but we will only use their valuations at $v$ and the ideals of $O_{F_v}$ generated by them.

For convenience, we recall the matrix notation:
\begin{align*}
 m(a)&=\matrixx{a}{}{}{a^{-1}},\quad  d(a)=\matrixx{1}{}{}{a}, \quad
d^*(a)=\matrixx{a}{}{}{1}\\
 n(b)&=\matrixx{1}{b}{}{1}, 
 \quad w=\matrixx{}{1}{-1}{}.
\end{align*}
We denote by $P\subset \gl$ and $P^1\subset \sll$ the subgroups of upper
triangular matrices, and by $N$ the standard unipotent subgroup of them. 

For any place $v$ of $F$, the character
$$\delta_v: P(F_v)\longrightarrow \RR\cross, \quad
\matrixx{a}{b}{}{d} \longmapsto \left| \frac ad \right|_v^{\frac 12}$$
extends to a function $\delta_v: \gl(F_v)\rightarrow \RR\cross$ by the Iwasawa
decomposition.
For the global field $F$, the product $\delta=\prod_v \delta_v$ gives a function
on $\gl(\adele)$.

If $v$ is a real place, we define a function 
$$\rho_v:
\gl(F_v)\longrightarrow \BC, \quad
g_v\longmapsto e^{i\theta},$$
where
$$g_v=  \matrixx{a}{b}{}{d}  \matrixx{\cos \theta}{\sin \theta}{-\sin \theta}{\cos \theta} $$
is in the form of the Iwasawa decomposition, where we require $a>0$ so that the decomposition is unique.
For the global field $F$, the product $\rho_\infty=\prod_{v\mid\infty} \rho_v$ gives a function
on $\gl(\adele)$, which ignores the non-archimedean components.

The following are all the conventions of this paper that are different from those of \cite{YZZ, YZ}, while only (3) is a major difference which brings extra computations.
\begin{enumerate}[(1)]
\item The Petersson metric on $\CL_U$ is defined by $\|d\tau\|_{\rm Pet}=2\, \Im(\tau)$ in \cite{YZ} and the current paper, while it is defined by 
$\|d\tau\|_{\rm Pet}=4\pi\, \Im(\tau)$ in \cite{YZZ}.
This discrepancy does not affect our applying results of \cite{YZZ}, since only the curvature form of the Petersson metric is crucial in \cite{YZZ}.

\item 
In the current paper, $\zeta_F(s)$ denotes the usual Dedekind zeta function (without Gamma factors),  $\tilde \zeta_F(s)$ denotes the completed Zeta function (with Gamma factors),
and $L(s,\eta)$ denotes the completed L-function (with Gamma factors) of the quadratic character $\eta$. 
In \cite{YZZ,YZ}, both $\zeta_F(s)$ and $L(s,\eta)$ denote the completed L-functions (with Gamma factors).

\item Our choice of $(U,\phi)$ in \S\ref{choices} is different from those in \cite{YZZ,YZ} due to the dropping of the degeneracy assumptions. 
Moreover, \cite{YZ} and the current paper eventually assume that $U$ is maximal compact, while \cite{YZZ} does not.  
We will mention this difference and its effect from time to time.

\item This paper and \cite{YZZ} do not assume $|\Sigma|>1$, while \cite{YZ} assumes $|\Sigma|>1$. Most results of \cite{YZ} actually hold in the case $|\Sigma|=1$.

\end{enumerate}

\subsubsection*{Acknowledgment}
The author is indebted to Shou-Wu Zhang and Wei Zhang, as the current paper is inspired by the long-term joint works of the author with them on the Gross--Zagier formula and the averaged Colmez conjecture. 

The author would like to thank Congling Qiu for pointing out a few mistakes in an early version of this paper, and thank Tuoping Du, Benedict H. Gross, Ulf K\"uhn, Yifeng Liu, Vincent Maillot, and Michael Rapoport for helpful communications.
Finally, the author is grateful to the anonymous referee for so many valuable comments or suggestions to revise this paper.

The author is supported by a grant from the National Science Foundation
of China (grant NO. 12250004) and the Xplorer Prize from the New
Cornerstone Science Foundation.

\section{Pseudo-Eisenstein series} \label{sec pseudo}

In \cite[\S 6]{YZ}, the notion of pseudo-theta series is introduced, and its crucial property in \cite[Lemma 6.1]{YZ} is the key to get a clean identity from the matching of the major terms. The goal of this section is to introduce a notion of pseudo-Eisenstein series and extend \cite[Lemma 6.1]{YZ} to a result including both pseudo-theta series and pseudo-Eisenstein series.
 
Throughout this section, let $F$ be a totally real number field, and $\BA$  the adele ring of $F$. We will use the terminologies of \S\ref{sec notation} and \cite[\S 6.1]{YZ} freely.

\subsection{Theta series and Eisenstein series} \label{sec theta eisenstein}

The goal here is to recall Weil representations, theta series and Eisenstein series following
\cite{YZZ} and \cite[\S 6.1]{YZ}.

\subsubsection*{Weil representations}

Let us first review  Weil representations associated to  quadratic space and quaternion algebras. We will only introduce the local case, and apply similar conventions for the global base. Our conventions are compatible with those in \cite[\S 2.1, \S2.2, \S 4.1]{YZZ}.

Let $k$ be a local field.
Let $(V,q)$ be a quadratic space over $k$ of even dimension. 
If $k=\RR$, we further assume that $(V,q)$ is positive definite.
We do not consider the case $k=\CC$ in this paper.

If $k$ is non-archimedean, denote by 
$\OCS(V\times k^{\times})=\CS(V\times k^{\times})$ the usual space of locally constant, compactly supported, and complex-valued functions on
$V\times k^{\times}$. 

If $k=\RR$, then $\OCS(V\times k^{\times})$ is the space of
functions on $V\times k^{\times} $ of the form
$$\phi(x,u)=\left(P_1(uq(x))+\sgn (u)P_2(uq(x))\right)e^{-2\pi |u|q(x)},$$
where $P_1, P_2$ are any polynomials of complex coefficients, and $\sgn (u)=u/|u|$ denotes the sign of
$u$.
The standard Schwartz function in this case is the standard Gaussian function
$$
\phi(x,u)= \begin{cases}
e^{-2\pi uq(x)}, & u>0, \\
0, & u<0.
\end{cases}
$$

In \cite[\S 2.1.3]{YZZ}, 
the Weil representation on the usual space $\CS(V)$ is extended to a representation of 
$\GL_2(k)\times \GO(V)$ on  
$\OCS(V\times k^{\times})$. 
Note that the actions of $\GL_2(k)$ and $\GO(V)$
commute with each other. 
This extension is originally from Waldspurger \cite{Wa}. 
As a convention, the action of $(g,h)\in \GL_2(k)\times \GO(V)$ on  
$\phi\in \OCS(V\times k^{\times})$ is denoted by $r(g,h)\phi$.
We also write $r(g)\phi=r(g,1)\phi$ and $r(h)\phi=r(1,h)\phi$.

The Weil representation behaves well under the direct sum of orthogonal spaces.
Assume that there is an orthogonal decomposition $V=V_1\oplus V_2$. 
Assume that $\phi=\phi_1\otimes\phi_2$ for 
$\phi\in  \overline \CS (V\times k^\times)$ and $\phi_i\in  \overline \CS (V_i\times k^\times)$ 
(with $i=1,2$) in the sense that 
$$
\phi(x_1+x_2,u)=\phi_1(x_1,u)\phi_2(x_2,u), \quad x_i\in V_i,\ u\in k^\times.
$$
Then the Weil representation splits as
$$
r(g)\phi(x_1+x_2,u)=r(g)\phi_1(x_1,u)\,r(g)\phi_2(x_2,u), \quad g\in\gla.
$$
Here $r(g)\phi, r(g)\phi_1, r(g)\phi_2$ are Weil representations corresponding to the quadratic spaces $V, V_1, V_2$ respectively.

Now we consider quadratic spaces coming from quaternion algebras. 
Assume that $V=B$ for a quaternion algebra $B$ over $k$, and assume that $q:B\to k$ is the reduced norm.  
If $k=\RR$, we further assume that $B$ is the Hamiltonian algebra so that $V$ is positive definite.
Let
$B^\times \times B^\times $ act on $V$ by
$$x\longmapsto h_1xh_2^{-1}, \qquad x\in V,\quad h_1,h_2\in B^\times.$$
This induces a natural map 
$$
\tau:B^\times \times B^\times \lra \GO(V).
$$
For $g\in \gl(k)$ and $h_1, h_2\in B^\times$,  denote 
$$r(g,(h_1, h_2))\phi=r(g,\tau(h_1, h_2))\phi, \quad
r(h_1, h_2)\phi=r(1,\tau(h_1, h_2))\phi.$$ 
We refer to \cite[\S2.2]{YZZ} for the kernel and the cokernel of $\tau$, but we do not need them in this paper. 

Finally, let $k'$ be either $k\oplus k$ or a separable quadratic field extension of $k$. 
Assume that there is an embedding $k'\to B$ of $k$-algebras. 
There is an element $j\in B^\times$ such that $jtj^{-1}=\bar t$ for any $t\in k'$. Here $\bar t$ denotes the action on $t$ by the unique non-trivial automorphism of $k'$ over $k$. 
This gives an orthogonal decomposition 
$$V=V_1\oplus V_2, \quad V_1=k', \ V_2=k'j,$$
where $V_1, V_2$ are endowed with the induced quadratic forms. Note that $k'j$, as a subset of $B$, does not depend on the choice of $j$.

Assume that $\phi=\phi_1\otimes\phi_2$ for 
$\phi\in  \overline \CS (V\times k^\times)$ and 
$\phi_i\in  \overline \CS (V_i\times k^\times)$ in the above sense.
Then the Weil representation splits as
$$
r(g, (t_1, t_2))\phi(x_1+x_2,u)=r(g, (t_1, t_2))\phi_1(x_1,u)\cdot r(g, (t_1, t_2))\phi_2(x_2,u), \quad g\in\gla, t_1, t_2\in k^\times.
$$
Here $r(g, (t_1, t_2))\phi$ is defined for the quadratic space $V$ by viewing $(t_1,t_2)\in B^\times\times B^\times$, and 
$$r(g, (t_1, t_2))\phi_i(x_i,u)=r(g, \tau_i(t_1, t_2))\phi_i(x_i,u)$$
 via the map $\tau_i: k^\times\times k^\times \to \GO(V_i)$ induced by the action 
 $(t_1, t_2)\circ x_i=t_1 x_i t_2^{-1}$ of $k^\times\times k^\times$ on $V_i$. 

For convenience, we also take the convention
$$
r(t_1, t_2)\phi(x,u)=r(1, (t_1, t_2))\phi(x,u), \quad
r(t_1, t_2)\phi_i(x_i,u)=r(1, (t_1, t_2))\phi_i(x_i,u).
$$

\subsubsection*{Theta series}
Let $(V,q)$ be a positive definite quadratic space over a totally real number field $F$.
Assume that $\dim V$ is even in the following, which is always satisfied in our application. 

Denote $V_v=V\otimes_FF_v$ for any place $v$ of $F$. 
Denote $V(\BA)=V\otimes_F \BA\simeq \prod_v' V_v$ for the restricted product. 
Let 
$$
\OCS(V(\adele)\times \adele^{\times})
=\otimes_v' \OCS(V_v\times F_v^{\times})
$$
be the space of Schwartz functions.
Here to defined the restricted tensor product, we need a distinguished vector in $\OCS(V_v\times F_v^{\times})$ for all but finitely many $v$. 
For that, fix a full $O_F$-lattice $\Lambda$ of $V$. 
For every non-archimedean place $v$, $\Lambda_v=\Lambda\otimes_{O_F}O_{F_v}$ is a full $O_{F_v}$-lattice of $V_v$.
Take the distinguished vector in $\OCS(V_v\times F_v^{\times})$ to be the characteristic function of $\Lambda_v\times O_{F_v}\cross$. 

By the product, we have the Weil representation of $\GL_2(\BA)\times \GO(V(\BA))$ on  
$\OCS(V(\adele)\times \adele^{\times})$. 
Note that the actions of $\GL_2(\BA)$ and $\GO(V(\BA))$
commute with each other. 

Take any $\phi\in  \overline \CS (V(\adele)\times \BA^\times)$.
There is the partial theta series
$$
\theta(g,u,\phi)= \sum_{x\in V} r(g)\phi(x,u), \quad g\in\gla, \ u\in \across.
$$
If $u\in F\cross$, it is invariant under the left multiplication of $\SL_2(F)$ on $g$.

To get an automorphic form on $\gla$, we define 
\begin{equation*}
\theta(g, \phi)_K
=\sum_{u\in \mu_K^2\backslash F\cross} \theta(g,u, \phi)
= \sum_{u\in \mu_K^2\backslash F\cross} \sum_{x\in V} r(g)\phi(x,u),  \quad g\in\gla.
\end{equation*}
Here $\mu_K=F\cross \cap K$, and $K$ is any open compact subgroup of 
$\GO(V(\adele_f))$ such that $\phi_f$ is invariant under the action of $K$ by the Weil representation.
The summation is well-defined and absolutely convergent.
The result $\theta(g, \phi)_K$ is an automorphic form in $g\in \gl(\adele)$, and
$\theta(g, r(h)\phi)_K$ is an automorphic form in $(g,h)\in \gl(\adele)\times
\GO(V(\adele))$. See \cite[\S 4.1.3]{YZZ} for more details.

Furthermore, if the infinite component $\phi_\infty$ is standard, i.e., for any archimedean place $v$,
$$
\phi_v(x,u)= \begin{cases}
e^{-2\pi uq(x)}, & u>0, \\
0, & u<0, 
\end{cases}
$$
then $\theta(g, \phi)_K$ is holomorphic of parallel weight $\frac 12 \dim V$.

The Weil representation and theta series behaves well under direct sum of orthogonal spaces.
Assume that there is an orthogonal decomposition $V=V_1\oplus V_2$ and a decomposition $\phi=\phi_1\otimes\phi_2$ for 
$\phi\in  \overline \CS (V(\adele)\times \BA^\times)$ and $\phi_i\in  \overline \CS (V_i(\adele)\times \BA^\times)$ in the above sense; i.e.  
$$
\phi(x_1+x_2,u)=\phi_1(x_1,u)\phi_2(x_2,u).
$$
Then the splitting of the Weil representation gives a natural splitting
$$
\theta(g,u,\phi)=\theta(g,u,\phi_1)\theta(g,u,\phi_2), \quad g\in\gla, \ u\in \across.
$$

\subsubsection*{Eisenstein series}

In the above setting of $\phi\in \overline \CS (V(\adele)\times \BA^\times)$ for a quadratic space $(V,q)$ over $F$, we can define an Eisenstein series $E(s,g, \phi)$. Then $E(s,g, \phi)$ and $\theta(g, \phi)$ are related by 
the Siegel--Weil formula. On the other hand, we also have Eisenstein series associated to incoherent quadratic collections in the sense of Kudla \cite{Kud}. 
For convenience, we introduce the notion of adelic quadratic spaces to include both cases by dropping the last condition of \cite[Definition 2.1]{Kud}. 

A collection $\{(\BV_v, q_v)\}_v$ of quadratic spaces 
$(\BV_v, q_v)$ over $F_v$ indexed by the set of places $v$ of $F$ is called \emph{adelic} if it satisfies the following conditions:
\begin{enumerate}[(1)]
\item There is a quadratic space $(V,q_0)$ over $F$ and a finite set $S$ of places of $F$ such that there is an isomorphism $(V_v,q_0) \to (\BV_v, q)$ for all  $v\notin S$;
\item For any place $v$ of $F$, the quadratic spaces $(V_v,q_0)$ and
$(\BV_v, q)$ have the same dimension and the same discriminant. 
\end{enumerate}
In that case, we obtain a quadratic space 
$$
(\vv, q):={\prod_v}' (\vv_v, q_v)
$$
over $\BA$. Here the restricted product makes sense by condition (1).  
We call $(\vv, q)$ an \emph{adelic quadratic space over $\BA$}. 
The dimension $\dim\vv\in \ZZ$ and the quadratic character $\chi_{(\vv,q)}:F^\times\bs \across\to \BC^\times$ are defined to be those of $(V,q)$. Its Hasse invariant is defined to be
$$
\epsilon(\vv,q):=\prod_v \epsilon(\vv_v,q_v).
$$
We say that the adelic quadratic space $(\vv, q)$ is \emph{coherent} 
(resp. \emph{incoherent}) if $\epsilon(\vv,q)=1$ (resp. $\epsilon(\vv,q)=-1$).
Note that $(\vv, q)$ is coherent if 
it is isomorphic to $(V(\BA), q_0)$ for some quadratic space $(V,q_0)$ over $F$. 

Let $(\vv, q)$ be an adelic quadratic space over $\BA$ which is positive definite at all archimedean places. 
For simplicity, we still assume that $\dim\vv$ is even. 
The space 
$$\OCS(\vv\times \adele^{\times})=\otimes_v'\OCS(\vv_v\times F_v^{\times})$$ 
is defined.
To explain the restricted tensor, for any $v\notin S$, we have an isomorphism $V_v\to \vv_v$ as above, so the distinguished vectors of $\{\OCS(V_v\times F_v^{\times})\}_{v\notin S}$ transfer to those of
$\{\OCS(\vv_v\times F_v^{\times})\}_{v\notin S}$. 

The Weil representation of $\GL_2(\BA)\times \GO(\vv)$ on  
$\OCS(\vv\times \adele^{\times})$ is defined by local products as in the coherent case.

Let $\phi\in \OCS(\BV\times \across)$ be a Schwartz function. Recall the associated partial Siegel Eisenstein series
\begin{eqnarray*}
E(s,g, u, \phi) &=& \sum_{\gamma \in P(F)\bs \GL_2(F)}
\delta (\gamma g)^s r(\gamma g)\phi(0, u) \\
&=& \sum_{\gamma \in P^1(F)\bs \SL_2(F)}
\delta (\gamma g)^s r(\gamma g)\phi(0, u), \quad g\in\gla, \ u\in \BA^\times.
\end{eqnarray*}
Here $P^1$ (resp. $P$) denotes the algebraic subgroup of upper triangle matrices in $\SL_2$ (resp. $\GL_2$), and $\delta$ is the standard modulus function as in \cite[\S1.6.6]{YZZ}. 
If $u\in F\cross$, it is invariant under the left multiplication of $\SL_2(F)$ on $g$, and it has a meromorphic continuation to $s\in\BC$
and a functional equation with center $s=1-\frac{\dim V}{2}$.

To get an automorphic form on $\gla$, we define 
\begin{equation*}
E(s,g,\phi) _K
=\sum_{u\in \mu_K^2\backslash F\cross} E(s,g,u, \phi),  \quad g\in\gla.
\end{equation*}
Here as before, $\mu_K=F\cross \cap K$, and $K$ is any open compact subgroup of 
$\GO(\adele_f)$ such that $\phi_f$ is invariant under the action of $K$ by the Weil representation.
It is easy to see that $E(s,g,\phi)_K$ is invariant under the left multiplication of 
$\GL_2(F)$ on $g$.
See \cite[\S 4.1.4]{YZZ}.

The Eisenstein series $E(s,g,u,\phi)$ has the standard Fourier expansion
$$E(s,g,u,\phi)=\delta(g)^s  r(g)\phi(0,u)+\sum_{a \in F} W_a(s,
g,u,\phi).$$
Here the Whittaker function is given by
\begin{eqnarray*}
W_a(s,g,u,\phi) = \int_{\adele} \delta(wn(b)g)^s \ r(wn(b)g)\phi(0,u)
\psi(-ab) db, \quad a\in F, \ u\in F\cross.
\end{eqnarray*}
Here $w=\matrixx{}{1}{-1}{}$ and $n(b)=\matrixx{1}{b}{}{1}$. 
For each place $v$ of $F$ and any $\phi_v\in \OCS(\vv_v\times \fvcross)$, we also introduce the local Whittaker function 
\begin{eqnarray*}
W_{a,v}(s,g,u,\phi_{v}) = \int_{F_v} \delta(wn(b)g)^s \
r(wn(b)g)\phi_{v}(0,u)
\psi_v(-ab) db, \quad a\in F_v, \ u\in F_v\cross.
\end{eqnarray*}

We have a splitting 
$$
E(s,g,u,\phi)=E_0(s,g,u,\phi)+E_*(s,g,u,\phi),
$$
where the constant term 
$$
E_0(s,g,u,\phi)=\delta(g)^s  r(g)\phi(0,u)+ W_0(s, g,u), 
$$
and the non-constant part 
$$
E_*(s,g,u,\phi)=\sum_{a \in F^\times} W_a(s, g,u,\phi).
$$
Similarly, we have 
$$
E(s,g,\phi)_K=E_0(s,g, \phi)_K+E_*(s,g,\phi)_K,
$$
where the constant term
$$
E_0(s,g,\phi) _K
=\sum_{u\in \mu_K^2\backslash F\cross} E_0(s,g,u, \phi)
$$
and the non-constant part
$$
E_*(s,g,\phi) _K
=\sum_{u\in \mu_K^2\backslash F\cross} E_*(s,g,u, \phi).
$$

If $(\BV, q)$ is coherent, then we can express $E(0, g, \phi)$ and $E(0, g, u,\phi)$
in terms of the theta series by the Siegel--Weil formula (in most convergent
cases).

If $(\BV, q)$ is incoherent, then there is no theta series available. However,
we can still express $W_{a,v}(0, g, u, \phi_v)$ in terms of certain average of the Schwartz function $\phi_v$. See the local Siegel--Weil formula in 
\cite[Theorem 2.2]{YZZ}. See also the examples of incoherent Eisenstein series (for $\SL_2$ or $\wt\SL_2$) in \cite[\S2.5]{YZZ}.

\subsection{Pseudo-Eisenstein series}
\label{sec pseudo}

Here we introduce the notion of pseudo-Eisenstein series, which is parallel to the notion of pseudo-theta series in \cite[\S 6.2]{YZ}. Note that the term 
``pseudo-Eisenstein series'' is also used in the literature as an unrelated terminology. 

\subsubsection*{Definition}
Let $(V, q)$ be an even-dimensional quadratic space over a totally real number field $F$, positive definite at all archimedean places. 
Let $S$ be a fixed finite set of non-archimedean place of $F$, and
$$\phi^S=\otimes_{w\notin S} \phi_w \in \overline\CS(V(\adele^S)\times \adele^{S, \times})$$ 
be a Schwartz
function with standard archimedean components.
A \textit{pseudo-Eisenstein series} is a series of the form
$$B^{(S)}_{\phi}(g)=\sum_{u\in \mu^2\backslash F\cross} \sum_{a\in F^\times} B_{a,S}(g,u) W_a^S(0,g,u,\phi^S), \ \quad g\in \gl(\adele).$$
We explain the notations as follows:
\begin{itemize}
\item $W_a^S(0,g,u,\phi^S)=\prod_{w\notin S} W_{a,w}(0,g,u,\phi_w)$ is the product of the local Whittaker functions defined before.
\item $B_{a,S}(g,u)=\prod_{v\in S} B_{a,v}(g,u)$ is the product of the local terms.
\item For any $v\in S$, the function
 $$B_{\bullet,v}(\bullet, \bullet): F_v\cross \times \gl(F_v) \times F_v\cross \rightarrow \BC$$
 is locally constant. 
It is smooth in $g$ in the sense that there is an open compact subgroup $K_v$ of $\gl(F_v)$ such that
 $$B_{a,v}(g\kappa,u)=B_{a,v}(g,u),
 \quad \forall (a, g,u)\in F_v\cross\times \gl(F_v)\times  F_v\cross,\  \kappa\in K_v.$$
 It is compactly supported in $u$ in the sense that there is a compact subset $D_g$ of $F_v^\times$ depending on $g$ (but independent of $a$) such that $B_{a,v}(g,u)=0$ for any $(a,g,u)$ with $u\notin D_g$.  
\item $\mu$ is a subgroup of $O_F\cross$ of finite index which acts trivially on the variable $u$ of $B_{a,v}(g,u)$ and $W_{a,w}(0,g,u,\phi_w)$ for every non-archimedean $w\notin S$ and $v\in S$.
\item For any $g\in \gla$, the double sum is absolutely convergent.
\end{itemize}
Note that $B^{(S)}_{\phi}(g)$ does not have a ``constant term'' in the sense that 
the summation is over $a\in F^\times$. 

\begin{example}
Let $(\vv, q)$ be an adelic quadratic space over $\BA$ which is positive definite at infinity, and $\phi\in \OCS(\BV\times \across)$ be a Schwartz function which is standard at infinity.
Consider the non-constant part 
$$
E_*(s, g, \phi)=\sum_{u\in \mu_K^2\backslash F\cross} \sum_{a\in F^\times} W_a(s,g,u,\phi).
$$
Its derivative at $s=0$ is
$$
E_*'(0, g, \phi)=\sum_{u\in \mu_K^2\backslash F\cross} \sum_{a\in F^\times} \sum_v W_{a,v}'(0,g,u,\phi)W_a^v(0,g,u,\phi^v).
$$
For every non-archimedean $v$, the ``$v$-part''
$$E_*'(0, g, \phi)(v)=\sum_{u\in \mu_K^2\backslash F\cross} \sum_{a\in F^\times} 
W_{a,v}'(0,g,u,\phi)W_a^v(0,g,u,\phi^v)$$
is a pseudo-Eisenstein series if it is absolutely convergent. 
For archimedean $v$, the ``$v$-part'' is not a pseudo-Eisenstein series by our definition, but a holomorphic projection will convert it to a multiple of 
$E_*(0, g, \phi)$. 
\end{example}

\subsubsection*{Non-singular pseudo-Eisenstein series}

Let $B^{(S)}_{\phi}(g)$ be the pseudo-Eisenstein series associated to $(V,q)$ as above. 
For every $v\in S$, there are one or two quadratic spaces over $F_v$ up to isomorphism with the same dimension and the same discriminant as $(V_v, q)$. 
Order them by $(V_v^+, q^+)$ and $(V_v^-, q^-)$ so that their Hasse invariants 
$\epsilon(V_v^+, q^+)=1$ and $\epsilon(V_v^-, q^-)=-1$.
If there is only one such space, which happens when $V$ is isomorphic to the 2-dimensional hyperbolic space over $F_v$, ignore the notation $V_v^-$.  

The pseudo-Eisenstein series $B^{(S)}_{\phi}(g)$ is called \textit{non-singular} if for every $v\in S$, there exist $\phi_v^+\in \OCS(V_v^+\times F_v^\times)$ and $\phi_v^-\in \OCS(V_v^-\times F_v^\times)$ such that 
$$B_{a,v}(1,u)=W_{a,v}(0,1,u,\phi_v^+)+ W_{a,v}(0,1,u,\phi_v^-),
 \quad \forall (a,u)\in F_v\cross\times  F_v\cross.$$
We take the convention that $W_{a,v}(0,1,u,\phi_v^-)=0$ if $V_v^-$ does not exist. 
Note that the equality is only for $g_v=1$. 
Once this is true, replacing $B_{a,v}(g,u)$ by $W_{a,v}(0,g,u,\phi_v^+)+ W_{a,v}(0,g,u,\phi_v^-)$ in $B^{(S)}_{\phi}(g)$, we see that 
$$
B^{(S)}_{\phi}(g)=\sum_{\epsilon: S\to \{\pm\}}
E_*(0,g, \phi_S^{\epsilon}\otimes \phi^S), \quad \ \forall g\in 1_S \gl(\BA^S). 
$$
Here $\phi_S^{\epsilon}=\otimes_{v\in S}\phi_v^{\epsilon(v)}$ is the Schwartz function associated to the adelic quadratic space $V_S^\epsilon\times V(\BA^v)$ with $V_S^{\epsilon}=\otimes_{v\in S}V_v^{\epsilon(v)}$, and 
$E_*(0,g, \phi_S^{\epsilon}\otimes \phi^S)$ denotes the non-constant part of the Eisenstein series $E_*(0,g, \phi_S^{\epsilon}\otimes \phi^S)$. 
If $V_v^-$ does not exist, take the convention that $\epsilon(v)=+$ for every $\epsilon$. 

This is the counterpart of the approximation formula for pseudo-theta series in \cite[\S6.2]{YZ}. Hence, it is convenience to denote 
$$
E_{B}(g)=E_{B^{(S)}_\phi}(g)=\sum_{\epsilon: S\to \{\pm\}}
E(0,g, \phi_S^{\epsilon}\otimes \phi^S).
$$
It is called \emph{the Eisenstein series associated to $B^{(S)}_{\phi}(g)$}. 

One can also formulate the terminology of pseudo-Eisenstein series for $\SL_2$ based on Schwartz functions in $\CS(V(\BA))$ instead of $\overline\CS(V(\BA)\times \across)$. 
It can be done based on principal series of $\SL_2$, which is really what an Eisenstein series needs.  
However, we stick with the current formulation because it fits our application.

\subsubsection*{Uniqueness of the Eisenstein series}

Resume the above notation for $B_\phi^{(S)}, (V_v^+, V_v^-), $ and $(\phi_v^+, \phi_v^-)$.
Assume that we are in the non-singular situation, so $(\phi_v^+, \phi_v^-)$ exists for $v\in S$.

The pair $(\phi_v^+, \phi_v^-)$ is generally not uniquely determined by $B_{a,v}(1,u)$. However, the following lemma asserts that the function 
$r(g_v)\phi_v^+(0,u)+ r(g_v)\phi_v^-(0,u)$ of $(g,u)\in \gl(F_v)\times F_v^\times$ is actually uniquely determined by $B_{a,v}(1,u)$. 
Thus the principal series 
$\sum_{\epsilon: S\to \{\pm\}}
\delta(g)^s r(g)(\phi_S^{\epsilon}\otimes \phi^S)(0,u)$ used to define the Eisenstein series 
$E_{B^{(S)}}(s,g)$ is also uniquely determined by $(\phi_v^+, \phi_v^-)$. 
As a consequence, the Eisenstein series 
$E_{B^{(S)}}(s,g)$, as a function of $g\in \GL_2(F_v)$ and $s\in \CC$, is uniquely determined by $B_{\phi}^{(S)}(g)$.

\begin{lem} \label{whittaker image1}
For any pair $(\phi_v^+, \phi_v^-)$ as above,  the function
$$f(g,u)=r(g)\phi_v^+(0,u)+r(g)\phi_v^-(0,u), \quad (g,u)\in \gl(F_v)\times F_v^\times$$
 is uniquely determined by the function
$$B_{a,v}(1,u)=W_{a,v}(0,1,u,\phi_v^+)+ W_{a,v}(0,1,u,\phi_v^-), \quad (a,u)\in F_v\cross\times  F_v\cross.$$
In particular, the relation gives
$$r(w)\phi_v^+(0,u)+r(w)\phi_v^-(0,u) = \int_{F_v}  B_{a,v}(1,u) da, \quad \
u\in F_v^\times. $$
\end{lem}

\begin{proof}
Recall the matrix notations $w,m(a),n(b),d(c)$ introduced in \S\ref{sec notation}. 

For any $c\in F_v^\times$, we have 
$f(d(c)g,u)=|c|^{-(\dim V_v^+)/4}f(g,c^{-1}u)$ by properties of the Weil representation.
Thus $f(g,u)=r(g)\phi_v^+(0,u)+r(g)\phi_v^-(0,u)$ is determined by its restriction to 
$\SL_2(F_v)\times F_v^\times$. For fixed $u$, it is a principal series of 
$g\in \SL_2(F_v)$. Then this is a classical result closely related to Kirillov models and has nothing to do with Weil representations. In fact, we need to recover $f(g,u)$ from 
\begin{eqnarray*}
B_{a,v}(1,u) = \int_{F_v} f(wn(b), u)
\psi(-ab) db, \quad a\in F, \ u\in F\cross.
\end{eqnarray*}
Observe that $B_{a,v}(1,u)$ as a function of $a\in F_v$ is the Fourier transform of $f(wn(b), u)$ as a function of $b\in F_v$. 
Thus we can recover $f(wn(b), u)$ by the Fourier inversion formula. 
Then $f(m(a)n(b')wn(b), u)$ can be recovered for any $a\in F_v^\times$ and $b'\in F_v$. 
But the set $m(a)n(b')wn(b)$ is dense in $\SL_2(F_v)$, as can be seen from the Bruhat decomposition. This determines all values of $f(g,u)$. 
In particular, the Fourier inversion formula gives 
$$f(w,u) = \int_{F_v}  B_{a,v}(1,u) da. $$
This finishes the proof.
\end{proof}

\subsection{The key lemma} \label{sec key lemma}

The goal of this section is to prove a key result (cf. Lemma \ref{pseudo}) for pseudo-theta series and pseudo-Eisenstein series. 
As a preparation, we first introduce approximation of 
pseudo-Eisenstein series, and also recall basic notions and approximation of pseudo-theta series.

\subsubsection*{Approximation by the Eisenstein series}

Let $B^{(S)}_{\phi}(g)$ be a non-singular pseudo-Eisenstein series as in the last subsection.
Resume the above notations for its corresponding Eisenstein series
$$
E_{B}(g)=\sum_{\epsilon: S\to \{\pm\}}
E(0,g, \phi_S^{\epsilon}\otimes \phi^S). 
$$
We already have an approximation
$$
B^{(S)}_{\phi}(g)=E_{B,*}(g), \quad \ \forall g\in 1_S \gl(\BA^S). 
$$
Here 
$$
E_{B,*}(g)=\sum_{\epsilon: S\to \{\pm\}}
E_*(0,g, \phi_S^{\epsilon}\otimes \phi^S)
$$
is the non-constant part of $E_{B}(g)$. 

In the following, we will approximate the constant term
$$
E_{B,0}(g)=\sum_{\epsilon: S\to \{\pm\}}
E_0(0,g, \phi_S^{\epsilon}\otimes \phi^S)
$$
by a function built from automorphic forms. 
Note that by definition, for every archimedean place $v$, 
the Schwartz function $\phi_v$ is standard in $\OCS(V_v\times F_v^\times)$ for the positive definite quadratic space $V_v$. 

By linearity, it suffices to treat a general Eisenstein series $E(0,g,\phi)_K$ associated to an adelic quadratic space $\vv$ over $\BA$ of even dimension $d$. 
Here we assume that at every archimedean place $v$, the quadratic space $\vv_v$ is positive definite, and
the Schwartz function $\phi_v$ is standard in $\OCS(\vv_v\times F_v^\times)$ . 
Recall the constant term
$$
E_0(0,g,\phi)_K= \sum_{u\in \mu_K^2\backslash F\cross}  r(g)\phi(0,u)+ 
\sum_{u\in \mu_K^2\backslash F\cross} W_{0}(0,g,u,\phi).
$$

We first treat the first term on the right-hand side. 
For any archimedean place $v$, as $\phi_v$ is standard, 
$r(g_v)\phi_v(0,u)$ gives the standard Whittaker function as in \cite[\S4.1.1]{YZZ}.
In particular,
$$
r(g_v)\phi_v(0,u)= 
\rho_v(g_v)^{\frac{d}{2}}\delta_v(g_v)^{\frac{d}{2}} 
\phi_v(0,\det(g_v)^{-1}u),\quad g_v\in \gl(F_v).
$$
Here $\delta_v$ and $\rho_v$ are introduced in \S\ref{sec notation}.

We claim that there is a finite set $S'$ of non-archimedean places of $F$ such that for every non-archimedean $v\notin S'$, 
$$
r(g_v)\phi_v(0,u)= 
\delta_v(g_v)^{\frac{d}{2}} 
\phi_v(0,\det(g_v)^{-1}u),\quad g_v\in \gl(F_v).
$$
In fact, there is a quadratic space $(V,q_0)$ over $F$ such that there is an isomorphism $(V_v,q_0) \to (\BV_v, q)$ for all but finitely many places $v$. 
Fix a full $O_F$-lattice $\Lambda$ in $V$. Then we can take $S'$ such that for any non-archimedean $v\notin S'$, the following holds:
\begin{enumerate}[(1)]
\item $v$ is unramified over $\QQ$; 
\item
the isomorphism $(V_v,q_0) \to (\BV_v, q)$ is available for $v$;
\item
the lattice $\Lambda_v=\Lambda\otimes_{O_F}O_{F_v}$ of $V_v$ is self-dual with respect to $\psi_v$;
\item
the Schwartz function $\phi_v$ is the standard characteristic function of $\Lambda_v\times O_{F_v}^\times$ via the isomorphism $(V_v,q_0) \to (\BV_v, q)$. 
\end{enumerate}
Then the claim follows form explicit computation of Weil representation.

Combining all the places together, we have
$$
\sum_{u\in \mu_K^2\backslash F\cross}  r(g)\phi(0,u)= 
\rho_\infty(g)^{\frac{d}{2}}\delta(g)^{\frac{d}{2}} \sum_{u\in \mu_K^2\backslash F\cross}  
\phi(0,\det(g)^{-1}u), \quad g\in 1_{S'}\gl(\adele^{S'}).
$$
Note that the summation on the right-hand side is automorphic in $g\in\gla$; i.e.
it is invariant under the left action of $\glf$ on $g\in\gla$. 
This finishes approximating the first term of on the right-hand side of $E_0(0,g,\phi)_K$.

For the second term on the right-hand side, by the above result for all places $v\notin S'$, $W_{0,v}(0,g,u,\phi_v)$ is a multiple of 
$\delta(g_v)^{2-d} r(g)\phi_v(0,u)$, as a basic result of intertwining operators of principal series. 
Consequently, we have a similar approximation 
$$
\sum_{u\in \mu_K^2\backslash F\cross} W_{0}(0,g,u,\phi)= 
\rho_\infty(g)^{\frac{d}{2}}\delta(g)^{2-\frac{d}{2}} 
\sum_{u\in \mu_K^2\backslash F\cross}  
W_{0}(0,1,\det(g)^{-1}u,\phi), \quad g\in 1_{S'}\gl(\adele^{S'}).
$$
The summation on the right-hand side is again automorphic in $g\in\gla$.

Return to the 
non-singular pseudo-Eisenstein series $B^{(S)}_{\phi}(g)$. 
Our conclusion gives an approximation
$$
B^{(S)}_{\phi}(g)=E_{B}(g)+
\rho_\infty(g)^{\frac{d}{2}}\delta(g)^{\frac{d}{2}} f_1(g)
+
\rho_\infty(g)^{\frac{d}{2}}\delta(g)^{2-\frac{d}{2}} f_2(g)
, \quad g\in 1_{S'}\gl(\adele^{S'}).
$$
Here $E_{B}(g)$ is the full Eisenstein series, and $f_1(g), f_2(g)$ are automorphic in $g\in\gla$. 
Both sides are invariant under the action of some (non-empty) open compact subgroup $K'_S$ of $\gl(\adele_{S'})$, 
so the approximation actually holds for all $g\in K'_{S'}\gl(\adele^{S'})$.

\subsubsection*{Review of pseudo-theta series}

For convenience, we briefly review the notion of pseudo-theta series and its approximation in
\cite[\S 6.2]{YZ}.

Let $V$ be a positive definite
quadratic space over $F$, and $V_0\subset V_1 \subset V$ be two
subspaces over $F$ with induced quadratic forms. All spaces are assumed to be even-dimensional. We allow $V_0$ to be the empty set $\emptyset$, which is not a subspace in the usual sense. Let $S$ be a finite set of non-archimedean places of $F$, and
$\phi^S \in \overline\CS(V(\adele^S)\times \adele^{S \times})$ be a Schwartz
function with standard infinite components.

A \textit{pseudo-theta series} is a series of the form
$$A^{(S)}_{\phi'}(g)=\sum_{u\in \mu^2\backslash F\cross} \sum_{x\in V_1-V_0} \phi_S'(g,x,u) r_{_V}(g)\phi^S(x,u), \ \quad g\in \gl(\adele).$$
We explain the notations as follows:
\begin{itemize}
\item The Weil representation $r_{_V}$ is not attached to  the space $V_1$
but to the space $V$;
\item $\phi_S'(g,x,u)=\prod_{v\in S} \phi'_v(g_v,x_v,u_v)$ as local product;
\item For each $v\in S$, the function
 $$\phi'_v: \gl(F_v)\times (V_1-V_0)(F_v)  \times F_v\cross \rightarrow \BC$$
 is locally constant. And it is smooth in the sense that there is an open compact subgroup $K_v$ of $\gl(F_v)$ such that
 $$\phi'_v(g\kappa,x,u)=\phi'_v(g,x,u),
 \quad \forall (g,x,u)\in \gl(F_v)\times (V_1-V_0)(F_v)\times F_v\cross,\  \kappa\in K_v.$$
\item $\mu$ is a subgroup of $O_F\cross$ with finite index such that $\phi^S(x,u)$ and $\phi'_S(g,x,u)$ are invariant under the action $\alpha: (x,u)\mapsto (\alpha x, \alpha^{-2} u)$ for any $\alpha\in \mu$. This condition makes the summation well-defined.
\item For any $v\in S$ and $g\in \gl(F_v)$, the support of $\phi'_v(g,\cdot,\cdot)$ in $(V_1-V_0)(F_v)\times F_v\cross$ is bounded. This condition makes the sum convergent.

\end{itemize}

The pseudo-theta series $A^{(S)}$ sitting on the triple $V_0\subset
V_1 \subset V$ is called \textit{non-degenerate} (resp. \textit{degenerate}) if $V_1=V$ (resp. $V_1\neq V$ ). It is called
\textit{non-singular} if for each $v\in S$, the local component
$\phi'_v(1,x,u)$ can be extended to a Schwartz function on
$V_1(F_v)\times F_v\cross$.

Assume that $A^{(S)}_{\phi'}$ is non-singular. Then there are two
usual theta series associated to $A^{(S)}$. View
$\phi'_v(1,\cdot,\cdot)$ as a Schwartz function on $V_1(F_v)\times
F_v\cross$ for each $v\in S$, and $\phi_w$ as a Schwartz function on
$V_1(F_w)\times F_w\cross$ for each $w\notin S$. Then the theta
series
$$\theta_{A,1}(g)=\sum_{u\in \mu^2\backslash F\cross} \sum_{x\in V_1} r_{_{V_1}}(g)\phi_S'(1,x,u) r_{_{V_1}}(g)\phi^S(x,u)$$
is called \textit{the outer theta series associated to}
$A^{(S)}_{\phi'}$. 
Note that the Weil representation $r_{V_1}$ is based on the quadratic space $V_1$. 
Replacing the space $V_1$ by $V_0$, we get the
theta series
$$\theta_{A,0}(g)=\sum_{u\in \mu^2\backslash F\cross} \sum_{x\in V_0} r_{_{V_0}}(g)\phi_S'(1,x,u) r_{_{V_0}}(g)\phi^S(x,u).$$
We call it \textit{the inner theta series associated to}
$A^{(S)}_{\phi'}$. We set $\theta_{A,0}=0$ if $V_0$ is empty.

We introduce these theta series because the difference between
$\theta_{A,1}$ and $\theta_{A,0}$ somehow approximates $A^{(S)}$. 
More precisely, there exist a finite set $S'$ of non-archimedean places of $F$ 
and a (non-empty) open compact subgroup $K'_S$ of $\gl(\adele_{S'})$
satisfying
$$
A^{(S)}_{\phi'}(g)=\rho_\infty(g)^{\frac{d-d_1}{2}}\delta(g)^{\frac{d-d_1}{2}}
\theta_{A,1}(g) -
\rho_\infty(g)^{\frac{d-d_0}{2}}\delta(g)^{\frac{d-d_0}{2}}
\theta_{A,0}(g), \quad 
g\in K'_{S'}\gl(\adele^{S'}). 
$$
Here $d=\dim V$ and $d_i=\dim V_i$. 
The rough idea of the approximation is as follows. 
Approximate
$$A^{(S)}_{\phi'}(g)=\sum_{u\in \mu^2\backslash F\cross} \sum_{x\in V_1-V_0} \phi_S'(g,x,u) r_{_V}(g)\phi^S(x,u)$$
 by
$$
\sum_{u\in \mu^2\backslash F\cross} \sum_{x\in V_1-V_0} r_V(g_S)\phi_S'(1,x,u) r_{_V}(g^S)\phi^S(x,u),
$$
split the later as 
$$
\sum_{u\in \mu^2\backslash F\cross}\sum_{x\in V_1-V_0}=\sum_{u\in \mu^2\backslash F\cross}\sum_{x\in V_1}-\sum_{u\in \mu^2\backslash F\cross}\sum_{x\in V_0},
$$
and then compare the two double sums of the right-hand side with $\theta_{A,1}$ and $\theta_{A,0}$  respectively. 
We refer to \cite[\S6.2]{YZ} for a proof of this fact.

\subsubsection*{Key lemma}

The following is a generalization of \cite[Lemma 6.1(1)]{YZ} to a sum of pseudo-theta series and pseudo-Eisenstein series. 
There is also a generalization of \cite[Lemma 6.1(2)]{YZ}, but we omit it due to the complexity of the statement. 

\begin{lem}\label{pseudo}
Let $\{A_\ell^{(S_\ell)}\}_\ell$ be a finite set of non-singular
pseudo-theta series sitting on even-dimensional quadratic spaces $V_{\ell,0}\subset
V_{\ell,1}\subset V_{\ell}$. 
Let $\{B_j^{(S_j')}\}_j$ be a finite set of non-singular
pseudo-Eisenstein series sitting on even-dimensional quadratic spaces $V_{j}'$. 
Assume that the sum 
$$f(g)=\sum_\ell A_\ell^{(S_\ell)}(g)+ \sum_j B_j^{(S_j')}(g)$$
is automorphic for $g\in\gl(\adele)$. Then
$$f(g)=\sum_{\ell\in L_{0,1}} \theta_{A_\ell,1}(g)+ \sum_j E_{B_j}(g)$$
Here $L_{0,1}$ is the set of $\ell$ such that $V_{\ell,1}=V_{\ell}$ or equivalently $A_\ell^{(S_\ell)}$ is non-degenerate.
\end{lem}

\begin{proof}
The proof is similar to that of \cite[Lemma 6.1(1)]{YZ}, by taking extra care of the pseudo-Eisenstein series. 
In fact, in the equation 
$$f-\sum_\ell
A_\ell^{(S_\ell)}-\sum_j B_j^{(S_j')}=0,$$
 replace all $A_\ell^{(S_\ell)}$ and $B_j^{(S_j')}$ by their approximations described above.
After recollecting these theta
series according to the powers of $\rho_\infty(g)$ and $\delta(g)$, we
end up with an equation of the form
\begin{eqnarray*}
\sum_{(k,k')\in I} \rho_\infty(g)^{k}\delta(g)^{k'} f_{k,k'}(g)=0, \quad \
\forall g\in K_{S}\gl(\adele^{S}).  
\end{eqnarray*}
Here $I$ is a finite subset of $(k,k')\in \ZZ^2$ with $k\geq 0$ and $2|(k-k')$, 
 $S$ is some finite set of non-archimedean places of $F$, 
 $K_S$ is an open compact subgroup of $\gl(\adele_{S})$, 
 and $f_{k,k'}$ is some
automorphic form on $\gl(\adele)$ coming from combinations of $f$, the corresponding theta series, the corresponding Eisenstein series, and approximation of constant terms of the Eisenstein series. 
In particular, 
$$f_{0,0}=f-\sum_{\ell\in L_{0,1}} \theta_{A_\ell,1}- \sum_j E_{B_j}$$
is the term we care about. 

We are going to prove that $f_{k,k'}=0$ for every $(k,k')\in I$. 
The proof is similar to that of \cite[Lemma 6.1]{YZ}, but it has a slight complication due to the possibility that $k\neq k'$. 
In fact, it suffices to show $f_{k, k'}(g_0)=0$ for all $g_0\in \gl(\adele_f^{S})$,
since $\gl(F)\gl(\adele_f^{S})$ is dense in $\gl(\adele)$. Fix such a
$g_0\in \gl(\adele_f^{S})$. For any $g \in \gl(F)\cap
K_S\gl(\adele^S)$, we have
$$\sum_{(k, k')\in I}\rho_\infty(gg_0)^{k}\delta(gg_0)^{k'} f_{k, k'}(gg_0)=0,$$
and thus
$$\sum_{(k, k')\in I}\rho_\infty(g)^{k}\delta(gg_0)^{k'} f_{k, k'}(g_0)=0$$
by the automorphy.
These are viewed as linear equations of $(f_{k,k'}(g_0))_{(k, k')\in I}$. To show that the solutions are zero, we only need to find sufficiently
many $g$ to get plenty of independent equations. We first find some
special $g$ to simplify the equation.

The intersection $K_S\gl(\adele^S) \cap g_0 \gl(\widehat
O_F)g_0^{-1}$ is still an open compact subgroup of $\gl(\adele)$.
For any $g\in \gl(F)\cap (K_S\gl(\adele^S) \cap g_0 \gl(\widehat
O_F)g_0^{-1})$, we have
$$gg_0=g_0 \cdot g_0^{-1} g g_0\in g_0\gl(\widehat O_F).$$
Then $\delta_f(gg_0)=\delta_f(g_0)$, and our linear equation
simplifies as
$$\sum_{(k, k')\in I}\rho_\infty(g)^{k}\delta_\infty(g)^{k'} \delta_f(g_0)^{k'}f_{k, k'}(g_0)=0.$$
It is reduced to find sufficiently many $g$ such that the equation implies that every term $\delta_f(g_0)^{k'}f_{k, k'}(g_0)=0$.

To be more explicit, set $g=\matrixx{1}{a}{b}{1+ab}$ for 
$a,b\in \ZZ$. Note that the proof of \cite[Lemma 6.1]{YZ} takes a simpler matrix corresponding to $a=0$ of the current one.
There is a positive integer $N$ such that if $a,b$ are divisible by $N$, then
$g$ lies in the intersection $\gl(F)\cap (K_S\gl(\adele^S) \cap g_0 \gl(\widehat
O_F)g_0^{-1})$ as required. 
The Iwasawa decomposition in $\gl(\RR)$ gives
$$\matrixx{1}{a}{b}{1+ab}=  \matrixx{\frac{1}{\sqrt{c}}}{\frac{a+b+a^2b}{\sqrt{c}}}{}{_{\sqrt{c}}}  \matrixx{\frac{1+ab}{\sqrt{c}}}{-\frac{b}{\sqrt{c}}}{\frac{b}{\sqrt{c}}}{ \frac{1+ab}{\sqrt{c}}} ,$$ 
where $c=(1+ab)^2+b^2$.
Then we have
$$\rho_v(g)= \frac{1+ab-ib}{\sqrt{c}},\quad
\delta_v(g)=\frac{1}{\sqrt{c}}$$
at every archimedean place $v$ of $F$. 
Thus the equation becomes
$$\sum_{(k, k')\in I} \left(\frac{1+ab-ib}{\sqrt{c}}\right)^{nk} \left(\frac{1}{\sqrt{c}}\right)^{nk'} \delta_f(g_0)^{k'}f_{k, k'}(g_0)=0\quad
\forall a,b\in N\ZZ.$$
Here $n=[F:\QQ]$. 

Taking advantage of the property $2|(k-k')$,  set $k''=-(k+k')/2$. 
The relation becomes
$$\sum_{(k, k')} (1+ab-ib)^{nk} c^{nk''} \delta_f(g_0)^{k'}f_{k, k'}(g_0)=0.$$
Note $c=(1+ab-ib)(1+ab+ib)$. So the problem is reduced to the following one. 
\emph{
Let $J$ be a finite subset of $\ZZ^2$. 
Assume that $h:J\to \CC$ is a map satisfying the equation
$$\sum_{(k_1, k_2)\in J} (1+ab+ib)^{k_1}(1+ab-ib)^{k_2} h(k_1,k_2)=0, \quad
\forall a,b\in N\ZZ.$$
Then $h(k_1,k_2)=0$ for every $(k_1, k_2)\in J$.}

To prove the result, by multiplying both sides of the equation by $(1+ab+ib)^m(1+ab-ib)^m$ for a sufficiently large integer $m$, we can assume that $k_1\geq 0, k_2\geq0$ for every $(k_1, k_2)\in J$. 
Now the left-hand side of the equation is a polynomial of $a,b$, which vanishes for all $a,b\in N\ZZ$. 
This implies that the left-hand side is 0 as a polynomial of $a,b$, so its value is 0 for any $a,b\in\CC$. 
Set $(x,y)=(1+ab, ib)$. 
We see that 
$$\sum_{(k_1, k_2)\in J} (x+y)^{k_1}(x-y)^{k_2} h(k_1,k_2)=0$$
for any $(x,y)\in \CC^2$ with $y\neq 0$.
Then it is an identity of polynomials of $x,y$.
A further relation $(z,w)=(x+y,x-y)$ implies that 
$$\sum_{(k_1, k_2)\in J} z^{k_1}w^{k_2} h(k_1,k_2)=0$$
as a polynomial of $(z,w)$.
It follows that 
 $h(k_1,k_2)=0$ for every $(k_1, k_2)\in J$.
This finishes the proof.
\end{proof}

\subsection{Example by local quaternion algebras}

In the case of quaternion algebras, we are going to figure out some important class of functions $B_{a,v}(1,u)$ which make the pseudo-Eisenstein series non-singular. 

Let $v$ be a non-archimedean place of $F$. 
Let $(M_2(F_v),q)$ (resp. $(D_v,q)$) be the matrix algebra (resp. the unique quaternion division algebra) over $F_v$ with the reduced norm.
Consider the map 
$$
\CW_v: \OCS(M_2(F_v) \times F_v^\times) \oplus \OCS(D_v\times F_v^\times)
\lra 
C^\infty(F_v^\times \times F_v^\times)
$$
given by 
$$
(\phi^+, \phi^-) \longmapsto W_{a,v}(0,1,u,\phi^+)+ W_{a,v}(0,1,u,\phi^-).
$$
Here $C^\infty(F_v^\times \times F_v^\times)$ denotes the space of locally constant functions with complex values, and the last expression is viewed as a function of 
$(a,u)\in F_v^\times \times F_v^\times$.

\begin{lem} \label{whittaker image2}
 Let $\Psi\in C^\infty(F_v^\times \times F_v^\times)$ be a linear combination of the function $1_{O_{F_v} \times O_{F_v^\times}}$ and a locally constant and compactly supported function on $F_v^\times \times F_v^\times$. 
Then $\Psi$ has a preimage $(\phi^+, \phi^-)$ satisfying
$$\phi^+(0,u)+\phi^-(0,u)=0, \quad \forall u\in F_v\cross.$$
\end{lem}
\begin{proof}
The problem is immediately reduced to two cases: 
\begin{itemize}
\item[(a)] $\Psi$ is a locally constant and compactly supported function on $F_v^\times \times F_v^\times$;
\item[(b)] $\Psi=1_{O_{F_v} \times O_{F_v^\times}}$.
\end{itemize}
As preparation, recall that 
the local Siegel-Weil formula in \cite[Proposition 2.9(2)]{YZZ} gives
$$
W_{a,v}(0,1,u, \phi)=
 \epsilon(B_v)\ |a|_v \int_{B_v^1} \phi(hx_a,u)dh, \quad a,u\in F_v^\times.
$$
Here $B_v^1=\{x\in B_v:q(x)=1\}$, the pair
$(B_v, \phi)$ can be either  $(M_2(F_v), \phi^+)$ or $(D_v, \phi^-)$, and $x_a\in B_v$ is any element satisfying 
$uq(x_a)=a.$  

Now we treat case (a). 
Note that $D_v^1$ is compact. 
We will actually find a preimage of the form $(0,\phi^-)$, where $\phi^-$ is invariant under the action of $D_v^1$.
In fact, the local Siegel--Weil formula gives
$$
\phi^-(x,u)=-\frac{1}{\vol(D_v^1)|uq(x)|_v}W_{uq(x),v}(0,1,u, \phi^-)
=-\frac{1}{\vol(D_v^1) |uq(x)|_v}\Psi(uq(x),u).
$$ 
It is a Schwartz function since $\Psi(a,u)$ is assumed to be compactly supported in $a$. It is also clear that $\phi^-(0,u)=0$ for any $u\in F_v\cross$.

For case (b), the local Siegel--Weil formula gives 
$$
W_{a,v}(0,1,u, 1_{O_{D_v}\times \ofv\cross})
=-|d_v|^{\frac{3}{2}} N_v^{-1}(1+N_v^{-1})\cdot |a|_v\cdot 
1_{O_{F_v}\times O_{F_v^\times}}(a,u).
$$ 
Here $O_{D_v}$ denotes the maximal order of $D_v$, 
 $d_v\in F_v$ is the local different of $F$ over $\BQ$,
 and
$\vol(D_v^1)=|d_v|^{\frac{3}{2}}N_v^{-1}(1+N_v^{-1})$
as normalized in \cite[\S1.6.2]{YZZ}.

On the other hand, 
$$
W_{a,v}(0,1,u, 1_{M_2(O_{F_v})\times \ofv\cross})
=|a|_v\cdot \vol(\SL_2(O_{F_v})) \cdot |\SL_2(O_{F_v})\bs M_2(O_{F_v})(a)|\cdot
1_{O_{F_v}\times O_{F_v^\times}}(a,u).
$$ 
Here $M_2(O_{F_v})(a)$ denotes matrices in $M_2(O_{F_v})$ of 
determinant $a$. 
Set $r=v(a)\geq 0$, and denote
$$
M_2(\ofv)_r=\{x \in M_2(\ofv)_r: v(\det(x))=r\}.
$$
Then we have 
$$
\SL_2(O_{F_v})\bs M_2(O_{F_v})(a)=\GL_2(O_{F_v})\bs M_2(O_{F_v})_r. 
$$
The last coset corresponds exactly to the classical Hecke correspondence $T(p_v^r)$, and its order is just 
$1+N_v+\cdots+ N_v^r$. 
Combine 
$\vol(\SL_2(O_k))=|d_v|^{\frac{3}{2}}(1-N_v^{-2})$
as normalized in \cite[\S1.6.2]{YZZ}.
We end up with 
$$
W_{a,v}(0,1,u, 1_{M_2(O_{F_v})\times \ofv\cross})
=|d_v|^{\frac{3}{2}} N_v^{-1}(1+N_v^{-1})\cdot (N_v-|a|_v)\cdot
1_{O_{F_v}\times O_{F_v^\times}}(a,u).
$$
The linear combination of these two expressions gives a preimage 
$$
\phi^+=
|d_v|^{-\frac{3}{2}} (1+N_v^{-1})^{-1}\cdot 1_{M_2(O_{F_v})\times \ofv\cross}, 
\quad\
\phi^-=-|d_v|^{-\frac{3}{2}} (1+N_v^{-1})^{-1}\cdot 1_{O_{D_v}\times \ofv\cross}.
$$
It is clear that $\phi^+(0,u)+\phi^-(0,u)=0$ for any $u\in F_v\cross$ in this case.
\end{proof}

\begin{remark}
In part (2), the result $\phi^+(0,u)+\phi^-(0,u)=0$ in case (b) is not as random as what our computational proof suggests.
In fact, we claim that for any image 
$$\Psi(a,u)=W_{a,v}(0,1,u,\phi^+)+ W_{a,v}(0,1,u,\phi^-),$$ 
if $\Psi$ can be extended to a locally constant and compactly supported function on $F_v \times F_v^\times$ (instead of the more restrictive $F_v^\times \times F_v^\times$), 
 then $\phi^+(0,u)+\phi^-(0,u)=0$.
 For a proof, for $b\in F_v^\times$, set 
$$
g=n(b)m(-b)wn(b)=\matrixx{1}{}{b^{-1}}{1}.
$$
The right hand side goes to 1 as the valuation $v(b)\to -\infty$. 
We have 
$$
r(g)\phi^+(0,u)+r(g)\phi^-(0,u)
=|b|_v^2\cdot \big(
r(wn(b))\phi^+(0,u)+r(wn(b))\phi^-(0,u) \big).
$$
Note that $r(wn(b))\phi^+(0,u)+r(wn(b))\phi^-(0,u)$ is the Fourier transform of $\Psi(a,u)$, so it is also a locally constant and compactly supported function in $b\in F_v$. 
In particular, it is zero if $v(b)$ is sufficiently negative. 
This proves $\phi^+(0,u)+\phi^-(0,u)=0$.
\end{remark}

\section{Derivative series} \label{sec derivative series}

The goal of this section is to study the holomorphic projection of the derivative of some mixed Eisenstein--theta series.
This section is based on \cite[\S7]{YZ} and \cite[\S 6]{YZZ}, but the situation is more complicated since we do not have \cite[Assumption 7.1]{YZ} or equivalently \cite[Assumption 5.4]{YZZ}.

\subsection{Derivative series} \label{sec 3.1}

Let $F$ be a totally real field, and $E$ be a totally imaginary quadratic extension of $F$. Denote by $\BA=\BA_F$ the ring of adeles of $F$.
Let $\BB$ be a totally definite incoherent quaternion algebra over 
$\BA$ with an embedding $E_\BA\to \BB$ of $\BA$-algebras.

Fix a Schwartz function $\phi\in \ol\CS (\BB\times \BA^\times)$ invariant under $U\times U$ for some open compact subgroup $U$ of $\bfcross$. 
Start with the mixed theta-Eisenstein series
\begin{equation*} 
I(s, g, \phi)_U
 = \sumu \sum_{\gamma \in P^1(F)\bs \SL_2(F)}
\delta (\gamma g)^s \sum_{x_1\in E} r(\gamma g)\phi (x_1, u), \quad
g\in\gla.
\end{equation*}
It was first introduced in \cite[\S 5.1.1]{YZZ}. 

There exists an element $\fj\in \BB^\times$ such that $\BB=E_\BA+E_\BA\fj$ is an orthogonal decomposition. This follows from the local version described in \S\ref{sec theta eisenstein}. 
Assume further that $\phi=\phi_1\otimes \phi_2$ for $\phi_1\in \ol\CS (E_\BA\times \BA^\times)$
and $\phi_2\in\ol\CS (E_\BA\fj\times \BA^\times)$
in the style of \S\ref{sec theta eisenstein}; i.e.
$$
\phi(x_1+x_2,u)=\phi_1(x_1,u)\phi_1(x_2,u), \quad x_1\in E_\BA,\ x_2\in E_\BA\fj,\ u\in \across.
$$
Then the splitting of the Weil representation gives a splitting
\begin{equation*}
I(s, g, \phi)_U=\sum_{u\in\mu_U^2\bs F\cross} \theta(g, u, \phi_1)\ E(s,g, u,
\phi_2),
\end{equation*}
where for any $g\in\gla$, the theta series and the Eisenstein series are given
by
\begin{eqnarray*}
\theta(g, u, \phi_1)&=& \sum_{x_1\in E} r(g)\phi_1(x_1, u), \\
E(s,g, u, \phi_2)&=& \sum_{\gamma \in P^1(F)\bs \SL_2(F)}
\delta (\gamma g)^s r(\gamma g)\phi_2(0, u).
\end{eqnarray*}

The derivative series $\Pr I'(0,g,\phi)$ is the holomorphic projection of the derivative $I'(0,g,\phi)$ of $I(s,g,\phi)$. 
We will start with some general results about the holomorphic projection.

\subsubsection*{Holomorphic projection}

Recall that the holomorphic projection is the orthogonal projection  
$$\pr: \CA(\gl(\adele), \omega)\longrightarrow \CA_0^{(2)}(\gl(\adele), \omega)$$
with respect to the Petersson inner product. 
Here $\omega:F^\times\bs\across\to \BC^\times$ is a Hecke character with trivial archimedean components,  
$\CA(\gl(\adele), \omega)$ is the space of
automorphic forms of central character $\omega$, and 
$\CA_0^{(2)}(\gl(\adele), \omega)$ is the subspace of holomorphic
cusp forms of parallel weight two.
It induces a projection
$$\pr: \bigoplus_\omega\CA(\gl(\adele), \omega)\longrightarrow 
\bigoplus_\omega\CA_0^{(2)}(\gl(\adele), \omega).$$

As in \cite[\S7.1]{YZ}, by decomposing $I'(0,g,\phi)$ into a finite direct sum of automorphic forms with (distinct) central characters, 
we see that $ I'(0,g,\phi)$ lies in $\oplus_\omega\CA(\gl(\adele), \omega)$.
Thus the holomorphic projection 
$\pr I'(0,g,\phi)$ is a well-defined holomorphic cusp form of parallel weight two in $g\in\gla$. 
We are still going to apply the formula in \cite[Proposition 6.12]{YZZ} to compute $\pr I'(0,g,\phi)$.

To recall \cite[Proposition 6.12]{YZZ}, we start with the operator $\pr'$ defined right after the proposition. 
For convenience, we first introduce the corresponding operator $\pr'_\psi $ for Whittaker functions. 
For any (Whittaker) function $\alpha:\gl(\RR)\to \BC$ with $\alpha(n(b)g)=\psi(b)\alpha(g)$ for any $b\in \RR$ and  $g\in\gl(\RR)$, define 
$$(\pr'_\psi \alpha)(g):=4\pi W^{(2)}(g)\cdot \quasilim \int_{Z(\RR)N(\RR)\bs \GL_2(\RR)} \delta(h)^s 
\alpha(h)\overline{W^{(2)}(h)}dh, 
$$
if the right-hand side is convergent. 
Here $W^{(2)}(g)$ is the standard Whittaker function of weight two as in \cite[\S4.1.1]{YZZ}, and $\quasilim$ is the constant term  in the Laurent expansion at $s=0$. 
The definition extends to global Whittaker functions $\alpha:\gl(\adele)\to \BC$ by  
$$(\pr'_\psi \alpha)(g)=(4\pi)^{[F:\QQ]}W_\infty^{(2)}(g_\infty)\cdot \quasilim \int_{Z(F_\infty)N(F_\infty)\bs \GL_2(F_\infty)} \delta(h)^s 
\alpha(g_fh)\overline{W^{(2)}(h)}dh
$$
if it is convergent. 

For any function $f:\gl(\adele)\to \BC$, we first take the Whittaker function 
$$
f_\psi(g)=\int_{N(F)\bs N(\adele)} f(n(b)g)\psi(-b)db,
$$
and set 
$$
(\pr'f)(g)=\sum_{a\in F^\times} (\pr'_\psi f_\psi)(d^*(a)g),
$$
if both are convergent in suitable sense. Here  $d^*(a)=\matrixx{a}{}{}{1}$ is as in \S\ref{sec notation}. 

Finally, \cite[Proposition 6.12]{YZZ} asserts that if $f$ is an automorphic form satisfying certain growth condition, then
$$
\pr f=\pr'f.
$$
In other words, the above formula really computes the holomorphic projection of $f$. 

Go back to $\pr I'(0,g,\phi)$.
In our previous works, \cite[Assumption 7.1]{YZ} or 
\cite[Assumption 5.4]{YZZ} makes $I'(0,g,\phi)$ satisfy the growth condition of 
\cite[Proposition 6.12]{YZZ}, but we do not make the assumption here, and we will see that the growth condition is not satisfied after dropping the assumption. 
Then the final result has an extra term contributed by the growth of 
$I'(0,g,\phi)$, as remarked in \cite[\S6.4.3]{YZZ}.

To track the growth of $I'(0,g,\phi)$, we are going to apply \cite[Lemma 6.13]{YZZ}.
Then we recall the absolute constant term
$$I_{00}(s, g, \phi)=\sumu I_{00}(s, g, u, \phi),$$
where
$$I_{00}(s, g, u, \phi)= \theta_0(g,u,\phi_1) E_0(s,g,u,\phi_2)$$
is the product of the constant terms. 
It follows that
$$I_{00}(0, g, u, \phi)= \delta(g)^s r(g)\phi(0,u)+r(g)\phi_1(0,u) W_0(s,
g,u,\phi_2)$$
is the sum of two principal series (by restricting to $\SL_2(\BA)$). 
For any $g'=\matrixx{a}{b}{}{d}$ with $a,d\in\AA^\times, b\in \AA$, we have
$$
\delta(g'g)^s r(g'g)\phi(0,u)=\left|\frac{a}{d}\right|_\AA^{\frac s2+2}\delta(g)^s r(g)\phi(0,(ad)^{-1}u)
$$
and 
$$
r(g'g)\phi_1(0,u) W_0(s,g'g,u,\phi_2)
=\left|\frac{a}{d}\right|_\AA^{-\frac s2+2}
r(g)\phi_1(0,(ad)^{-1}u) W_0(s,g,(ad)^{-1}u,\phi_2).
$$
Here $|c|_\AA=\prod_v |c_v|_v$ for $c\in \across$.

Let $\CJ(s, g, u, \phi)$ be the
Eisenstein series formed by $I_{00}(s, g, u,\phi)$; i.e.
$$\CJ(s, g, u, \phi)=\sum _{\gamma \in P^1(F)\bs \SL_2(F)}I_{00}(s, \gamma g, u,\phi).$$
Denote
$$\CJ(s, g, \phi)_U=\sumu \CJ(s, g, u, \phi),$$
which will also be abbreviated as $\CJ(s, g, \phi)$.

By \cite[Lemma 6.13]{YZZ}, the difference 
$$I'(0,g,\phi)_U-I_{00}'(0,g,\phi)_U$$
satisfies the growth condition of \cite[Proposition 6.12]{YZZ}. 
Note that the original statement is about the twisting of the difference by some character $\chi$, but a similar proof by decomposing the difference in terms of central characters works for the current situation. 
A similar argument proves that 
$$\CJ'(0,g,\phi)_U-I_{00}'(0,g,\phi)_U$$
also satisfies the growth condition. 
As a consequence, 
$$I'(0,g,\phi)_U-\CJ'(0,g,\phi)_U$$
satisfies the growth condition.

Therefore, we have
\begin{multline*}
\pr(I'(0, g,\phi))=\pr(I'(0, g,\phi)-\CJ'(0, g, \phi)) \\
=\pr'(I'(0, g,\phi)-\CJ'(0, g, \phi)) 
=\pr'(I'(0, g,\phi))-\pr'(\CJ'(0, g, \phi)),
\end{multline*}
where the operator $\pr'$ is defined by the algorithm as recalled above. The term $\pr'(I'(0, g,\phi))$ is computed exactly as in \cite[Theorem 7.2]{YZ}. 
For its importance, we summarize the result in the following, and we will recall the notations after the statement. 

\begin{thm} \label{analytic series}
Assume that $\phi$ is standard at infinity. Then 
$$
\pr I'(0, g,\phi)_U
=\pr' I'(0, g,\phi)_U-\pr'\CJ'(0, g, \phi)_U,
$$
where $\pr' I'(0, g,\phi)_U$ has the same expression as that of $\pr I'(0, g,\phi)_U$ in \cite[Theorem 7.2]{YZ}. 
Namely, 
\begin{align*}
\pr' I'(0,g,\phi)_U
=& -\sum_{v|\infty}\overline {I'}(0,g,\phi)(v)-\sum_{v\nmid\infty \
\nonsplit}I'(0,g,\phi)(v)\\
& -c_1 \sumu \sum_{y\in E^\times} r(g)\phi(y,u) 
 -\sum_{v\nmid\infty} \sumu\sum_{y\in E^\times}
c_{\phi_v}(g,y,u)\, r(g)\phi^v(y,u)\\
&+ \sumu \sum_{y\in E^\times} (2\log \delta_f(g_f)+ \log|uq(y)|_f)\ r(g)\phi(y,u),
\end{align*}
and the right-hand side is explained in the following.
\begin{itemize}
\item[(1)]
For any archimedean $v$,
\begin{align*}
\overline {I'}(0,g,\phi)(v) &= 2 \barint_{C_U}
 \overline \CK^{(v)}_{\phi}(g,(t,t))dt, \\
\overline \CK^{(v)}_{\phi}(g,(t_1,t_2)) &= w_U \sum_{a\in F\cross}
\quasilim
\sum_{y\in \mu_U\backslash (B(v)_+\cross -E\cross)}  r(g,(t_1,t_2)) \phi(y)_a\
k_{v, s}(y), \\
k_{v, s}(y) &= \frac{\Gamma(s+1)}{2(4\pi)^{s}}
\int_1^{\infty} \frac{1}{t(1-\lambda(y)t)^{s+1}}  dt,
\end{align*}
where $\lambda(y)=q(y_2)/q(y)$ is viewed as an element of $F_v$.
\item[(2)]
For any non-archimedean $v$ which is nonsplit in $E$,
\begin{align*}
I'(0,g,\phi)(v) &= 2 \barint_{C_U} \CK^{(v)}_{\phi}(g,(t,t))dt, \\
 \CK^{(v)}_{\phi}(g,(t_1,t_2))
&= \sum_{u\in \mu_U^2\backslash F\cross} \sum_{y\in B(v)-E}
k_{r(t_1,t_2)\phi_v}(g,y,u)  r(g,(t_1,t_2)) \phi^v(y,u), \\
k_{\phi_v}(g,y,u)&= \frac{L(1,\eta_v)}{\mathrm{vol}(E_v^1)} r(g)
\phi_{1,v}(y_1,u) {W_{uq(y_2),v}^\circ}'(0,g,u,\phi_{2,v}), \quad y_2\neq 0.
\end{align*}
Here the last identity holds under the relation $\phi_v= \phi_{1,v}\otimes \phi_{2,v}$, and the definition extends by linearity to general $\phi_v$.
\item[(3)]
The constant 
$$ c_1
=2\frac{L'_f(0,\eta)}{L_f(0,\eta)} +\log|d_E/d_F|.$$
\item[(4)] 
Under the relation $\phi_v= \phi_{1,v}\otimes \phi_{2,v}$, for $y\in E_v$ and $u\in \fvcross$,
$$c_{\phi_v}(g,y,u)=r(g)\phi_{1,v}(y,u) W_{0,v}^{\circ}\ '(0,g,u, \phi_{2,v}) 
+ \log \delta(g_v)r(g)\phi_v(y,u).$$
The definition extends by linearity to general $\phi_v$.
\end{itemize}

\end{thm}

All the notations of the theorem are compatible with those of \cite{YZZ,YZ}. 
For convenience of readers, we recall them in the following. 
\begin{enumerate}[(a)]
\item In (1) and (2), the set
$C_U=E\cross\bs E^\times_{\af} / (U\cap E^\times_{\af})$ 
is a finite abelian group, and the averaged integration $\ds\barint_{C_U}$ is just the usual average over this finite group. 
\item
For any place $v$ of $F$ nonsplit in $E$, $B(v)$ denotes the nearby quaternion algebra over $F$ obtained by switching the Hasse invariant of $\BB$ at $v$. 
Fix an embedding $E\to B(v)$ of $F$-algebras. 
For any place $w\neq v$, fix an isomorphism 
$i_w:B(v)_w\simeq \BB_w$ compatible with the embeddings $E\to B(v)$ and $E_\AA\to \BB$.
This gives identifications $B(v)_{\AA^v}\simeq \BB^v$ and 
$\OCS(B(v)_{\AA^v}\times (\AA^v)^\times)\simeq \OCS(\bb^v\times (\AA^v)^\times). $
Thus we can view $\phi^v$ as an element of $\OCS(\bb^v\times (\AA^v)^\times)$.
These identifications also appear in \S\ref{choices}. 

\item 
The Weil representation $r(g,(t_1,t_2))$ takes the convention of the local case in \S\ref{sec theta eisenstein}. 
Moreover, in (1), we further take the convention 
$$r(g,(t_1,t_2)) \phi(y)_a=r(g^v,(t_1^v,t_2^v)) \phi^v(y, aq(y)^{-1})\, W^{(2)}_{a_v}(g_v).$$
Here $W^{(2)}_{a_v}(g_v)$ is the standard holomorphic Whittaker function of weight two as in \cite[\S4.1.1]{YZZ}. 

\item In (1),  $B(v)_+^\times$ denotes the subgroup of elements of $B(v)$ with totally positive norms. 

\item In (1) and (2), with the fixed embedding $E\to B(v)$, and we have an orthogonal decomposition $B(v)=E+E j(v)$. Under the decomposition, we can write $y=y_1+y_2$ for $y\in B(v)$. This explains $y_1, y_2$ (depending on $y$), and also explains $\lambda(y)=q(y_2)/q(y)$.

\item
In (2) and (4), the term ${W_{a,v}^\circ}(s,g,u,\phi_{2,v})$ is normalized as in \cite[\S7.1]{YZ}. Namely, 
for $a\in F_v\cross$, define
\begin{eqnarray*}
W_{a,v}^{\circ}(s, g,u,\phi_{2,v})
= \gamma_{u,v}^{-1}  W_{a,v}(s, g,u,\phi_{2,v}).
\end{eqnarray*}
Here $\gamma_{u,v}$ is the Weil index of $(E_v\jv, uq)$, where $\BB_v=E_v+E_v\jv$ is an orthogonal decomposition. 
For $a=0$, define
\begin{eqnarray*}
W_{0,v}^{\circ}(s, g,u,\phi_{2,v})
= \gamma_{u,v}^{-1}\frac{L(s+1,\eta_v)}{L(s,\eta_v)}
|D_v|^{-\frac{1}{2}}|d_v|^{-\frac{1}{2}} W_{0,v}(s, g,u,\phi_{2,v}).
\end{eqnarray*}
Here as in \S\ref{sec notation}, $d_v\in F_v$ is the local different of $F$ over $\BQ$,
and $D_v\in F_v$ is the discriminant of the quadratic extension $E_v$ of $F_v$.
\end{enumerate}

\subsubsection*{Contribution of the Eisenstein series}

Now we compute the term $\pr'\CJ'(0, g, \phi)$.

For any $\phi \in\overline\CS(\BB\times \BA^\times)$ of the form
$\phi=\phi_1\otimes \phi_2$,
$$I_{00}'(0, g, u, \phi)= \log\delta(g)\ r(g)\phi(0,u)+r(g)\phi_1(0,u) W_0'(0,
g,u,\phi_2).$$
Recall that the computation in \cite[Proposition 6.7]{YZZ} or that in \cite[\S7.1, p. 586]{YZ} gives
\begin{align*}
W_0'(0,g,u,\phi_2) = -c_0 r(g)\phi_2(0,u) - \sum_v 
r(g^v)\phi_{2}^v(0,u) W_{0,v}^{\circ}\ '(0,g_v,u,\phi_{2,v})
\end{align*}
with the constant
$$c_0=\der\left(\log \frac{L(s,\eta)}{L(s+1,\eta)}\right)=2\frac{L'(0,\eta)}{L(0,\eta)} +\log|d_E/d_F|.$$
Here $L(s,\eta)$ is the completed L-function with gamma factors, and we have used the functional equation 
$$
L(1-s,\eta)=|d_E/d_F|^{s-\frac12} L(s,\eta).
$$
It follows that
\begin{eqnarray*}
&& I_{00}'(0, g, u, \phi) \\
&=& \log\delta(g)\cdot r(g)\phi(0,u)- c_0 r(g)\phi(0,u)
-\sum_v r(g^v)\phi^{v}(0,u)\cdot r(g_v)\phi_{1,v}(0,u) W_{0,v}^{\circ}\
'(0,g_v,u,\phi_{2,v})\\
&=& 2\log\delta(g)\cdot r(g)\phi(0,u)- c_0 r(g)\phi(0,u)
-\sum_v r(g^v)\phi^{v}(0,u)c_{\phi_v}(g,0,u),
\end{eqnarray*}
where as in Theorem \ref{analytic series}(4), for $y\in E_v$ (including the case $y=0$) and $g\in \gl(F_v)$,
\begin{eqnarray*}
c_{\phi_v}(g,y,u)&=& r(g)\phi_{1,v}(y,u) W_{0,v}^{\circ}\, '(0,g,u,\phi_{2,v}) 
+ \log \delta(g_v)r(g)\phi_v(y,u) \\
&=& r(g)\phi_{1,v}(y,u) \big( W_{0,v}^{\circ}\, '(0,g,u,\phi_{2,v}) 
+ \log \delta(g_v)r(g)\phi_{2,v}(0,u) \big).
\end{eqnarray*}
Here the sum on $v$ is actually a finite sum, since $c_{\phi_v}(g,y,u)=0$ for all archimedean $v$ and for all but finitely many non-archimedean $v$ (uniformly in $(g,y,u)$) by 
\cite[Proposition 6.1(1)]{YZZ} and 
\cite[Lemma 7.6(1)]{YZ}. 

One checks that $c_{\phi_v}(g,0,u)$ is a principal series in the sense that
$$
c_{\phi_v}(m(a)n(b)g,0,u)=|a|_v^{2}\ c_{\phi_v}(g,0,u), \quad a\in F_v^\times, \ b\in F_v.
$$
This is a consequence of the basic fact
$$
W_{0,v}^{\circ}(s,m(a)n(b)g,u,\phi_{2,v})=|a|_v^{1-s}\eta_v(a) W_{0,v}^{\circ}(s,g,u,\phi_{2,v})
$$
and the result
$$
W_{0,v}^{\circ}(0,g,u,\phi_{2,v})=r(g)\phi_{2,v}(0,u)
$$
of \cite[Proposition 6.1]{YZZ}.

Then we introduce Eisenstein series
\begin{align*}
E(s, g, u, \phi)=&\sum _{\gamma \in P^1(F)\bs \SL_2(F)} \delta(\gamma g)^s r(\gamma g)\phi(0,u),\\
C(s, g, u, \phi)(v)=& \sum _{\gamma \in P^1(F)\bs \SL_2(F)} 
\delta(\gamma g)^s  c_{\phi_v}(\gamma g,0,u)\ r(\gamma g^v)\phi^v(0,u),
\end{align*}
and 
\begin{align*}
E(s, g, \phi)_U=&\sumu E(s, g, u,\phi),\\
C(s, g, \phi)_U(v)=&\sumu C(s, g, u,\phi)(v).
\end{align*}
Denote also 
\begin{align*}
C(s, g, u, \phi)=& \sum_{v\nmid\infty} C(s, g, u, \phi)(v),\\
C(s, g, \phi)_U=& \sum_{v\nmid\infty} C(s, g,\phi)_U(v).
\end{align*}
Note that both summations have only finitely many nonzero terms
since $c_{\phi_v}(g,0,u)=0$ for all but finitely many $v$ as mentioned above. 
We will usually supress the sub-index $U$ in $E(s, g, \phi)_U$ and $C(s, g, \phi)_U$.
 
Return to the derivative 
$$\CJ'(0,g,\phi)=\sumu \CJ'(0,g,u,\phi)$$
with 
$$
\CJ'(0,g,u,\phi)=\sum_{\gamma \in P^1(F)\bs \SL_2(F)} I_{00}'(0, g, u, \phi). 
$$
Then the Eisenstein series constructed based on
\begin{eqnarray*}
 I_{00}'(0, g, u, \phi) 
= 2\log\delta(g)\cdot r(g)\phi(0,u)- c_0 r(g)\phi(0,u)
-\sum_v r(g^v)\phi^{v}(0,u)c_{\phi_v}(g,0,u)
\end{eqnarray*}
is given by
\begin{align*}
\CJ'(0,g,\phi) = 2 E'(0,g,\phi)-c_0 E(0,g,\phi)- C(0, g,\phi).
\end{align*}
Applying the formula for $\pr'$, we have the following expression.

\begin{pro} \label{analytic series extra}
Assume that $\phi$ is standard at infinity. Then 
$$
\pr'\CJ'(0, g, \phi)
=- (c_0+(1+\log 4)[F:\QQ] ) E_*(0,g,\phi)- C_*(0, g,\phi)
+2\sum_{v\nmid\infty} E'(0,g,\phi)(v).$$
Here $E_*$ and $C_*$ are the non-constant parts of the Eisenstein series $E$ and $C$.
For any $v\nmid \infty$,
$$ E'(0,g,\phi)(v)
=\sumu E'(0,g,u,\phi)(v),$$
where
$$ E'(0,g,u,\phi)(v)
=\sum_{a\in F^\times} W_a^v(0,g,a^{-1}u,\phi^v)
\left(W_{a,v}'(0,g,a^{-1}u,\phi_v)-\frac12\log|a|_v\cdot W_{a,v}(0,g,a^{-1}u,\phi_v)\right).$$
\end{pro}

\begin{proof}
Recall  the non-constant part
$$
E_*(s, g, \phi)=\sumu E_*(s, g, u,\phi).
$$
Recall the Fourier expansions
$$
E_*(s, g, u,\phi)=\sum_{a\in F^\times} W_a(s, g, u,\phi)
$$
and 
$$
E_*(s, g,\phi)=\sum_{a\in F^\times} W_a(s, g, \phi).
$$
We can recover $E_*(s, g, \phi)$ from the Whittaker functions at $a=1$ by 
$$
E_*(s, g,\phi)= \sum_{a\in F^\times} W_1(s, d^*(a)g,\phi)
=\sumu \sum_{a\in F^\times} W_1(s, d^*(a)g, u,\phi).
$$

By linearity, 
$$
\pr'\CJ'(0, g, \phi)
=2 \pr'E'(0,g,\phi)-c_0 \pr' E(0,g,\phi)- \pr' C(0, g,\phi)
$$
Note that the Whittaker function $W_1(0,g,u,\phi)$ of $E(0,g,u,\phi)$ is a product of local Whittaker functions, and the local component $W_{1,v}(0,g,u,\phi)$ at every archimedean $v$ is a multiple of the standard Whittaker function $W^{(2)}$ of weight two. Then the Whittaker function 
$W_1(0,g,\phi)$ of $E(0,g,\phi)$ is still a 
product of the standard Whittaker functions at archimedean places with a finite part. 
Then as in the proof of \cite[Proposition 6.12]{YZZ}, an explicit calculation shows that the holomorphic
projection $\pr'$ does not change $W_1(0,g, \phi)$. 
As a consequence,
$$ \pr' E(0,g,\phi)= E_*(0,g,\phi). $$
Similarly, 
$$ \pr' C(0,g,\phi)= C_*(0,g,\phi). $$

For $\pr' E'(0,g,\phi)$, start with the Whittaker function
$$
W_1'(0,g,u,\phi)=\sumu W_1'(0,g,u,\phi)
=\sumu \sum_{v} W_{1,v}'(0,g,u,\phi_v) W_1^v(0,g,u,\phi^v).
$$
Then it amounts to apply $\pr'_\psi $ to $W_{1,v}'(0,g,u,\phi_v) W_1^v(0,g,u,\phi)$ for each place $v$ of $F$.

If $v|\infty$, then 
$$\pr'_\psi W_{1,v}'(0,g,u,\phi_v)=c_3 W_{1,v}(0,g,u,\phi_v)$$
for some constant $c_3$. 
It follows that 
$$
\pr'_\psi \left(W_{1,v}'(0,g,u,\phi_v) W_1^v(0,g,u,\phi^v) \right)
=c_3 W_{1}(0,g,u,\phi).
$$
Recovering its contribution to the whole series, we get 
$$
c_3 E_*(0,g,\phi)=c_3 \sumu E_*(0,g,u,\phi).
$$
Furthermore, the constant $c_3=-\frac 12(1+\log 4)$
is computed in Lemma \ref{constant by holomorphic projection} below.

If $v\nmid \infty$, then 
$\pr'_\psi $ does not change 
$W_{1,v}'(0,g,u,\phi_v) W_1^v(0,g,u,\phi^v)$ since it is already holomorphic. 
However, when getting back to the whole series, its contribution is 
$$
\sum_{a\in F^\times} 
W_{1,v}'(0,d^*(a)g,u,\phi_v) W_1^v(0,d^*(a)g,u,\phi^v).
$$
Apply the basic result
$$
W_{1,v'}(s,d^*(a)g,u,\phi_{v'})
=|a|_{v'}^{-\frac s2}W_{a,v'}(s,g,a^{-1}u,\phi_{v'}),
$$
which can be verified using $wn(b)d^*(a)=d(a)wn(a^{-1}b)$. 
We have
$$
W_1^v(0,d^*(a)g,u,\phi^v)
=W_a^v(0,g,a^{-1}u,\phi^v),
$$
and 
$$
W_{1,v}'(0,d^*(a)g,u,\phi_v) 
=W_{a,v}'(0,g,a^{-1}u,\phi_v)-\frac12\log|a|_v\cdot W_{a,v}(0,g,a^{-1}u,\phi_v).
$$
Then the result follows. 
\end{proof}

\subsection{Choice of the Schwartz function} \label{choices}

To make further explicit local computations, we need to put conditions on $(F,E,\BB, U)$ and 
specify the Schwartz function $\phi$. 
We will see that our choices are slightly different from that of \cite[\S7.2]{YZ}.
Throughout this paper, we also assume the basic conditions, and assume part of the restrictive conditions from time to time.

\subsubsection*{Basic conditions} 

Start with the setup of Theorem \ref{main} and Theorem \ref{height CM}.

Let $F$ be a totally real field, and $E$ be a totally imaginary quadratic extension of $F$. Let $\BB$ be a totally definite incoherent quaternion algebra over $\BA=\BA_F$ with an embedding $E_\BA\to \BB$ of $\BA$-algebras.
Then the ramification set $\Sigma$ of $\BB$ is of order cardinality and contains all the archimedean places of $F$. 

Fix a maximal order $O_{\BB_f}=\prod_{v\nmid\infty} O_{\BB_v}$ of $\BB_f=\BB\otimes_\AA\af$ invariant under the main involution of $\BB_f$. 
This gives a maximal open compact subgroup $O_{\BB_f}^\times$ of $\bfcross$. 
Let $U=\prod_{v\nmid\infty} U_v$ be an open compact subgroup of $O_{\BB_f}^\times$. 


For any $v\nmid\infty$, fix a nonzero element $\jv\in O_{\bb_v}$ orthogonal to $E_v$ such that $v(q(\mathfrak{j}_v))$ is non-negative and minimal; i.e.,
$v(q(\mathfrak{j}_v))\in\{0, 1\}$. 
Then $v(q(\mathfrak{j}_v))=1$ if and only if $\bv$ is nonsplit (and thus $E_v/F_v$ is inert by assumption). 
The existence of $\jv$ is basic and verified in \cite[\S7.2]{YZ}.
For $v\mid \infty$, fix an nonzero element $\jv\in \bb_v$ orthogonal to $E_v$.
Then we have an element $\fj=(\jv)_v$ of $\BB$, which gives an orthogonal decomposition
$\BB=E_\AA+E_\AA\fj$.

For any place $v$ nonsplit in $E$, let $B(v)$ be the nearby quaternion algebra over $F$ obtained by changing the Hasse invariant of $\BB$ at $v$. Fix an embedding $E\to B(v)$ and isomorphisms $B(v)_{w}\simeq \bb_{w}$ for all $w\neq v$, which are assumed to be compatible with the embedding $E_{\BA}\to \BB$.

The compatible isomorphism $B(v)_{w}\simeq \bb_{w}$ always exists. In fact, if $E_w$ is a field, it follows from the classical Noether--Skolem theorem. If $E_w$ is not a field, we have a a splitting $E_w\simeq F_w\oplus F_w$, whose idempotents induce splittings $B(v)_v\simeq M_2(F_w)$ and $\bb_v\simeq M_2(F_w)$.

At $v$, we also take an element $j_v\in B(v)_v$ orthogonal to $E_v$, such that $v(q(j_v))$ is non-negative and minimal as above. 
We remark that this set $\{\mathfrak j_{v'}:v'\neq v\} \cup \{j_v\}$ is not required to be the localizations of a single element of $B(v)$.

Let $\phi=\otimes_v \phi_v$ be a  Schwartz function in $\OCS(\bb\times \AA^\times)$. 
The isomorphism $B(v)_{\AA^v}\simeq \BB^v$ induces an identification 
$\OCS(B(v)_{\AA^v}\times (\AA^v)^\times)\simeq \OCS(\bb^v\times (\AA^v)^\times). $
Thus we can view $\phi^v$ as an element of $\OCS(\bb^v\times (\AA^v)^\times)$ and we will do this all the time.

We make the following basic assumptions:
\begin{itemize}
\item[(a)] If $v$ is archimedean,  $\phi_v$ is the standard Gaussian function as in 
\S\ref{sec theta eisenstein}. 
\item[(b)] Assume that 
$U_v$ is maximal at every $v\in \Sigma_f$. 
\item[(c)] Assume that $E$ is inert at every $v\in \Sigma_f$. 
\end{itemize}

We will make these basic assumptions throughout most of the paper.

\subsubsection*{Restrictive conditions} 

Here we introduce some very restrictive assumptions, which are required in local computations, and in the proof of Theorem \ref{main}.  
All our results are true under these conditions, but we will mention from time to time when these conditions are really essential in computations.

Assume that $U=O_{\BB_f}\cross$ is maximal. 
Assume that the embedding $E_\AA\to \BB$ sends $\wh O_E$ into the maximal order $O_{\BB_f}$. 
As a consequence, $U=\prod_{v\nmid\infty} U_v$ contains (the image of) $\wh O_E^\times=\prod_{v\nmid\infty} O_{E_v}^\times$. 


As for the Schwartz function $\phi=\otimes_v \phi_v$, we make the following choices:
\begin{itemize}
\item[(1)] If $v$ is archimedean, set $\phi_v$ to be the standard Gaussian function as in \S\ref{sec theta eisenstein}. This is already mentioned above.
\item[(2)] If $v$ is non-archimedean and split in $\BB$, set $\phi_v$ to be the standard characteristic function $1_{O_\bv\times\ofv\cross}$. 
\item[(3)] If $v$ is non-archimedean and nonsplit in $\BB$, set $\phi_v$ to be  the characteristic function $1_{O_\bv^\times\times\ofv\cross}$ (instead of the standard $1_{O_\bv\times\ofv\cross}$). 
\end{itemize}
By definition, $\phi$ is invariant under both the left action and the right action of $U$. 

Note that \cite[\S7.2]{YZ} assumes that there is a set $S_2$ consisting of two places of $F$ split in $E$ such that $\phi_v$ takes a specific degenerate form for $v\in S_2$. 
We do not make this assumption here, since this assumption exactly kills the terms we need  for our main theorem.


\subsection{Explicit local derivatives}
Assume the basic conditions in \S\ref{choices} for the moment.
We will assume all the conditions (at the relevant places) when we do explicit local calculations.

Recall that we have defined the Eisenstein series
\begin{eqnarray*}
E(s, g,u,\phi)
= \sum_{\gamma\in P(F)\bs \gl(F)} \delta(\gamma g)^s  r (\gamma g)\phi (0,u).
\end{eqnarray*}
Note that this Eisenstein series uses the whole Schwartz function 
$\phi\in\OCS(\bb\times\across)$ and thus have weight two, comparing to the Eisenstein series $E(s, g,u,\phi_2)$ in the definition of the derivative series at the beginning of \S\ref{sec 3.1}, which only uses $\phi_2\in\OCS(E_\AA\fj\times\across)$ and thus has weight one.
In this section, we will abbreviate
$$
E(s, g,u)
=E(s, g,u,\phi).$$

We have the usual Fourier expansion:
$$E(s, g, u)=E_0(s, g,u) + \sum_{a \in F\cross} W_a(s, g,u)$$
with
\begin{eqnarray*}
E_0(s, g,u)
&=& \delta(g)^s  r (g)\phi (0,u)  + W_0(s, g,u),\\
W_a(s, g,u)
&=&  \int_{\adele} \delta(wn(b)g)^s  r(wn(b)g)\phi(0,u) \psi(-ab) db, \ \ a\in F.
\end{eqnarray*}
We also introduce the local Whittaker function
\begin{eqnarray*}
W_{a,v}(s, g,u)
=  \int_{F_v} \delta(wn(b)g)^s  r(wn(b)g)\phi(0,u) \psi(-ab) db, \ \ a\in F_v,\ u\in F_v\cross, \
g\in\gl(F_v).
\end{eqnarray*}

The goal of this subsection is compute the derivative of this local Whittaker function at $s=0$. 
In the following, assume all the conditions on $(F,E,\BB, U,\phi)$ listed in \S \ref{choices}.

\subsubsection*{Local holomorphic projection: archimedean place}

Recall that $\phi$ at any archimedean place is the standard Gaussian as in \S\ref{sec theta eisenstein}. 

\begin{lem} \label{constant by holomorphic projection}
Let $v$ be an archimedean place. 
\begin{itemize}
\item[(1)]
For any $a\in F_v$ with $a>0$,
\begin{align*}
W_{1,v}(0,d^*(a),u)
=& -4\pi^2  ae^{-2\pi a}1_{F_{v,+}\cross}(u),\\
W_{1,v}'(0,d^*(a),u)
=& -\left(\frac{\pi}{2} e^{-2\pi a} +2\pi^2 (\log\pi +\gamma-1) ae^{-2\pi a}\right) 1_{F_{v,+}\cross}(u). 
\end{align*}
Here $\gamma$ is Euler's constant.

\item[(2)] 
The holomorphic projection
\begin{align*}
\pr'_\psi \, W_{1,v}'(0,g,u)=-\frac 12  (1+\log 4) W_{1,v}(0,g,u).
\end{align*}

\end{itemize}
\end{lem}

\begin{proof}
We first check (1).
By $wn(b)d^*(a)=d(a)wn(a^{-1}b)$, it is easy to get 
$$W_{1,v}(s,d^*(a),u)=a^{-\frac s2} W_{a,v}(s,1,a^{-1}u).$$
 It is reduced to compute 
\begin{eqnarray*}
W_{a,v}(s,1,u)
=  \int_{F_v} \delta(wn(b))^s  r(wn(b))\phi(0,u) \psi(-ab) db,\quad a>0.
\end{eqnarray*}
Assume $u>0$; otherwise, the above vanishes.

The process is parallel to \cite[Proposition 2.11]{YZZ} and also uses the technique of \cite{KRY1}. 
In fact, from the proof of \cite[Proposition 2.11]{YZZ} for the case $d=4$, 
\begin{align*}
W_{a,v}(s,1,u)
=-   \frac{2\pi^{s+2}}{\Gamma(\frac s2+2) \Gamma(\frac s2)}   e^{-2\pi a}
\int_{0}^{\infty}  e^{-2\pi t} (t+2a)^{\frac s2+1}   t^{\frac s2 -1}  dt.
\end{align*}
To see its behavior at $s=0$, write 
\begin{align*}
W_{a,v}(s,1,u)
= -e^{-2\pi a} \frac{\pi^{s+2}s}{\Gamma(\frac s2+2) \Gamma(\frac s2+1)}   
\int_{0}^{\infty}  e^{-2\pi t} (t+2a)^{\frac s2+1}   t^{\frac s2 -1}  dt.
\end{align*}
Then the product before the integral has a simple zero at $s=0$.

Since
\begin{align*}   
& \int_{0}^{\infty}  e^{-2\pi t} (t+2a)^{\frac s2+1}   t^{\frac s2 -1}  dt 
-(2a)^{\frac s2+1}  (2\pi)^{-\frac s2}\Gamma(\frac s2)   \\
=& \int_{0}^{\infty}  e^{-2\pi t}  \frac{(t+2a)^{\frac s2+1}-(2a)^{\frac s2+1}}{t}   t^{\frac s2}  dt 
=\frac{1}{2\pi}+O(s),
\end{align*}
we get
\begin{align*}
-W_{a,v}(s,1,u)
=& e^{-2\pi a} \frac{\pi^{s+2}s}{\Gamma(\frac s2+2) \Gamma(\frac s2+1)} 
\left( \frac{1}{2\pi}  
+(2a)^{\frac s2+1}  (2\pi)^{-\frac s2}\Gamma(\frac s2)  \right)+O(s^2)\\
=& e^{-2\pi a} \frac{\pi^{s+2}}{\Gamma(\frac s2+2) \Gamma(\frac s2+1)} 
\left( \frac{1}{2\pi}  s
+2(2a)^{\frac s2+1}  (2\pi)^{-\frac s2}\Gamma(\frac s2+1)  \right)+O(s^2).
\end{align*}
It follows that 
\begin{align*}
-W_{a,v}(0,1,u)
=&\ 4 \pi^2 ae^{-2\pi a}, \\
-W_{a,v}'(0,1,u)
=&\ \frac{\pi}{2} e^{-2\pi a} +2\pi^2 (\log\pi +\gamma-1) ae^{-2\pi a}
+2\pi^2 ae^{-2\pi a}\log a,\\
-W_{1,v}'(0,d^*(a),u)
=& \frac{\pi}{2} e^{-2\pi a} +2\pi^2 (\log\pi +\gamma-1) ae^{-2\pi a}. 
\end{align*}

Now we compute the holomorphic projection 
$$\pr'_\psi W_{1,v}'(0,g,u)=4\pi W^{(2)}(g)\cdot \quasilim \int_{Z(\RR)N(\RR)\bs \GL_2(\RR)} \delta(h)^s 
W_{1,v}'(0,h,u)\overline{W^{(2)}(h)}dh.
$$
By the Iwasawa decomposition,  
\begin{align*}
\pr'_\psi W_{1,v}'(0,g,u)=
&4\pi W^{(2)}(g) \quasilim \int_0^{\infty} y^s e^{-2\pi y} W_{1,v}'(0,d^*(y),u)  \frac{dy}{y} \\
=&- 4\pi W^{(2)}(g) \quasilim \int_0^{\infty} y^s e^{-2\pi y} 
\left(\frac{\pi}{2} e^{-2\pi y} +2\pi^2 (\log\pi +\gamma-1) ye^{-2\pi y}\right)
\frac{dy}{y}.
\end{align*}
The integral above is computed by
\begin{align*}
& \int_0^{\infty} y^s e^{-2\pi y} 
\left(\frac{\pi}{2} e^{-2\pi y} +2\pi^2 (\log\pi +\gamma-1) ye^{-2\pi y}\right)
\frac{dy}{y} \\
=& \frac{\pi}{2} \int_0^{\infty} y^s e^{-4\pi y} \frac{dy}{y} 
+2\pi^2 (\log\pi +\gamma-1) \int_0^{\infty} y^{s+1} e^{-4\pi y} \frac{dy}{y} \\
=& \frac{\pi}{2} (4\pi)^{-s} \Gamma(s)
+2\pi^2 (\log\pi +\gamma-1) (4\pi)^{-s-1} \Gamma(s+1).
\end{align*}
Its constant term is equal to
\begin{align*}
 \frac{\pi}{2} (-\log(4\pi)-\gamma)
+2\pi^2 (\log\pi +\gamma-1) (4\pi)^{-1} = - \frac{1}{2}\pi (1+\log 4).
\end{align*}
Hence, 
$$\pr'_\psi W_{1,v}'(0,g,u)= 2\pi^2 (1+\log 4) W^{(2)}(g).$$
This holds for $u>0$.
By the result of (1),
$$
W_{1,v}(0,g,u)=- 4\pi^2 W^{(2)}(g) 1_{F_{v,+}\cross}(u).
$$
Then (2) follows.
\end{proof}

\subsubsection*{Derivative of Whittaker functions: non-archimedean place}

As in \S\ref{sec notation}, $p_v$ denotes the maximal ideal of $\ofv$, and
$d_v\in F_v$ denotes the local different of $F$ over $\BQ$. 

Assume all the conditions in \S\ref{choices}.
Recall that for a non-archimedean place $v$, we have $\phi_v=1_{O_\bv\times\ofv\cross}$  if $v$ is split in $\BB$, and $\phi_v=1_{O_\bv^\times\times\ofv\cross}$  if $v$ is nonsplit in $\BB$.

\begin{lem} \label{local explicit}
Let $v$ be a non-archimedean place of $F$, and let $a\in F_v^\times$. 
\begin{itemize}
\item[(1)] 
Let $v$ be a non-archimedean place split in $\bb$. 
Then $W_{a,v}(s,1,u)$ is nonzero only if $u\in \ofv\cross$ and 
$v(a)\geq -v(d_v)$. 
In the case $u\in \ofv\cross$ and $a\in \ofv$,
\begin{eqnarray*}
W_{a,v}(s,1,u)
&=&  |d_v|^{s+\frac{3}{2}}\frac{(1-N_v^{-(s+2)})(1-N_v^{-(v(a)+1)(s+1)})}{1-N_v^{-(s+1)}}\\
&& +\   |d_v|^{\frac{5}{2}}  
\frac{(1-N_v^{-s})(1-|d_v|^{s-1})}{1-N_v^{-(s-1)}}.
\end{eqnarray*}
Therefore, for $u\in \ofv\cross$ and $a\in \ofv$,
\begin{eqnarray*}
&& W_{a,v}'(0,1,u)-\frac12\log|a|_v W_{a,v}(0,1,u)  \\
&=&  \left(-\zeta_v'(2)/\zeta_v(2) + \log|d_v|\right) W_{a,v}(0,1,u)\\
&& +\,
|d_v|^{\frac32}  \frac{1+N_v^{-1}}{2(1-N_v^{-1})}
\big((r+2)N_v^{-(r+1)}-rN_v^{-(r+2)}-(r+2)N_v^{-1}+r\big) \log N_v \\
&& +\,   |d_v|^{\frac{3}{2}}  
\frac{1-|d_v|}{N_v-1}\log N_v.
\end{eqnarray*}
Here $r=v(a)$. 

\item[(2)]
Let $v$ be a non-archimedean place nonsplit in $\bb$. 
Then $W_{a,v}(s,1,u)$ is nonzero only if $v(a)\geq -v(d_v)$ and $u\in \ofv\cross$, and it is constant (depending on $s$) for $(a,u)\in p_v\times \ofv\cross$.
Moreover, for any $u\in F_v\cross$,
\begin{eqnarray*}
 \int_{F_v} \left(W_{a,v}'(0,1,u)-\frac12\log|a|_v W_{a,v}(0,1,u) \right)da=0.
\end{eqnarray*}
\end{itemize}

\end{lem}

\begin{proof}
The calculation is rather involved due to the non-triviality of $d_v$. 
To simplify the calculation, we move between two different types of methods. 
We divide the process into three steps due to the ramifications of $v$ in $\BB$ and over $\BQ$. 

\medskip
\noindent\emph{Step 1. unramified case:} $v$ is split in $\BB$ and $|d_v|=1$.  
Apply the formula
\begin{eqnarray*}
W_{a,v}(s,1,u) =
\int_{F_v} \delta(wn(b))^s \ r(wn(b))\phi_{v}(0,u) \psi_v(-ab) db.
\end{eqnarray*}
Note that $\phi_v$ is invariant under the action of $\GL_2(O_{F_v})$ (as symplectic similitudes), as $\GL_2(O_{F_v})$ is generated by $w, m(a), n(b)$ with all $a\in O_{F_v}^\times, b\in O_{F_v}$. 
Thus the Iwasawa decomposition gives
$$
r(g)\phi_v(0,u)=\delta(g)^2\ 1_{O_{F_v}^\times}(u), \quad g\in \SL_2(F_v).
$$
Notice
$$
\delta(wn(b)) = \left\{ \begin{array}{rl}
 1 \ \quad &\mbox{ if $b\in O_{F_v}$,} \\
 |b|^{-1} &\mbox{ otherwise}.
       \end{array} \right.
$$
Assuming $u\in O_{F_v}^\times$ (so that $W_{a,v}(s,1,u)\neq 0$), we have
\begin{eqnarray*}
W_{a,v}(s,1,u) 
&=&
\int_{F_v} \delta(wn(b))^{s+2} \psi_v(-ab) db \\
&=& \int_{O_{F_v}} \psi_v(-ab) db+ \int_{F_v-O_{F_v}} |b|^{-(s+2)} \psi_v(-ab) db.
\end{eqnarray*}
Write the domain $F_v-O_{F_v}$ of the second integral as a disjoint union of $p_v^{-n}- p_v^{-(n-1)}$ for $n\geq1$. We have
\begin{eqnarray*}
W_{a,v}(s,1,u) 
&=&
 \int_{O_{F_v}} \psi_v(-ab) db+ \sum_{n=1}^{\infty}\int_{p_v^{-n}- p_v^{-(n-1)}} N_v^{-n(s+2)} \psi_v(-ab) db\\
 &=&
 \int_{O_{F_v}} \psi_v(-ab) db+ \sum_{n=1}^{\infty}\int_{p_v^{-n}} N_v^{-n(s+2)} \psi_v(-ab) db
 -\sum_{n=1}^{\infty}\int_{p_v^{-(n-1)}} N_v^{-n(s+2)} \psi_v(-ab) db\\
  &=&
(1-N_v^{-(s+2)}) \sum_{n=0}^{\infty} N_v^{-n(s+2)} \int_{p_v^{-n}}  \psi_v(-ab) db.
\end{eqnarray*}
It is nonzero only if $a\in O_{F_v}$. 
In that case, 
\begin{eqnarray*}
W_{a,v}(s,1,u) 
=(1-N_v^{-(s+2)}) \sum_{n=0}^{v(a)} N_v^{-n(s+2)} N_v^{n}
= \frac{(1-N_v^{-(s+2)})(1-N_v^{-(v(a)+1)(s+1)})}{1-N_v^{-(s+1)}}.
\end{eqnarray*}

It follows that 
\begin{eqnarray*}
W_{a,v}(0,1,u)  =  \frac{(1-N_v^{-2})(1-N_v^{-(v(a)+1)})}{1-N_v^{-1}}
=  (1+N_v^{-1})(1-N_v^{-(v(a)+1)})
\end{eqnarray*}
and 
\begin{eqnarray*}
&& W_{a,v}'(0,1,u)-\frac12\log|a|_v W_{a,v}(0,1,u)  \\
&=&  W_{a,v}(0,1,u)\big(\frac{W_{a,v}'(0,1,u)}{W_{a,v}(0,1,u)}-\frac12\log|a|_v \big)
\\ 
&=&  W_{a,v}(0,1,u)
\big( \frac{N_v^{-2}}{1-N_v^{-2}}+ 
\frac{(v(a)+1)N_v^{-(v(a)+1)}}{1-N_v^{-(v(a)+1)}}
-\frac{N_v^{-1}}{1-N_v^{-1}}+\frac12 v(a)
\big)\log N_v\\
&=&  \frac{N_v^{-2}\log N_v}{1-N_v^{-2}} W_{a,v}(0,1,u)+  
\frac{1-N_v^{-2}}{2(1-N_v^{-1})^2}
\big((r+2)N_v^{-(r+1)}-rN_v^{-(r+2)}-(r+2)N_v^{-1}+r\big) \log N_v. 
\end{eqnarray*}
This proves part (1) under $|d_v|=1$. 

\medskip
\noindent\emph{Step 2. A general formula:} 
By the proof of \cite[Proposition 6.10(1)]{YZZ}, we have 
\begin{eqnarray*}
W_{a,v}(s,1,u,\phi_v)
= \gamma(\bv, uq)|d_v|^{\frac{1}{2}}  (1-N_v^{-s}) \sum_{n=0}^{\infty}
N_v^{-n(s-1)} \int_{D_n(a)} \phi_{v}(x,u) d_ux,
\end{eqnarray*}
where $d_ux$ is the self-dual measure of $(\bv, uq)$, and 
$$D_n(a)=\{x \in \bb_{v}:  uq(x) \in a+p_v^nd_v^{-1}\}$$
is a subset of $\bb_{v}$. The local Weil index $\gamma(\bv,uq)=\pm 1$ coincides with the Hasse invariant of $\bv$.
Note that the quadratic space $(\BB_v, uq)$ here is different from the quadratic space $(E_v\mathfrak j_v, uq)$ in \cite[Proposition 6.10(1)]{YZZ},  
but the proof is similar. 

It is easy to see that $W_{a,v}(s,1,u)\neq 0$ only if $u\in \ofv\cross$ and 
$v(a)\geq -v(d_v)$. In the following, we always assume $u\in \ofv\cross$ and 
$v(a)\geq 0$.

\medskip
\noindent\emph{Step 3. matrix case:} $v$ is split in $\BB$ and $d_v$ is arbitrary. 
By the normalization of $\psi_v:F_v\to \BC^\times$ in \cite[\S1.6.1]{YZZ}, the characteristic function $1_{O_{F_v}}$ is not self-dual under $\psi_v$ if $|d_v|\neq 1$. Consequently, $\phi_v$ is not invariant under the action of $\GL_2(O_{F_v})$. Then the method of Step 1 does not work in this case, and we are going to use the formula in Step 2. 

We have 
\begin{eqnarray*}
W_{a,v}(s,1,u)
=  |d_v|^{\frac{1}{2}}  (1-N_v^{-s}) \sum_{n=0}^{\infty}
N_v^{-n(s-1)} \vol(D_n(a)\cap O_\bv)
\end{eqnarray*}
with 
$$D_n(a)\cap O_\bv=\{x \in O_\bv:  uq(x) \in a+p_v^nd_v^{-1}\}.$$
By $\vol(O_\bv)=|d_v|^2$, we write
\begin{eqnarray*}
W_{a,v}(s,1,u)
=  |d_v|^{\frac{5}{2}}  (1-N_v^{-s}) \sum_{n=0}^{\infty}
N_v^{-n(s-1)} \frac{\vol(D_n(a)\cap O_\bv)}{\vol(O_\bv)}.
\end{eqnarray*}
We use this expression because it holds for any Haar measure on $\bv$. 
Split the summation according to $n<v(d_v)$ and $n\geq v(d_v)$. 
It gives 
$$W_{a,v}(s,1,u)=W_{a,v}(s,1,u)_{n<v(d_v)}+ W_{a,v}(s,1,u)_{n\geq v(d_v)}$$ 
accordingly. 
In the following we compute the terms on the right-hand side separately. 

By the assumption $a\in O_{F_v}$, for any $n<v(d_v)$, we have 
$$D_n(a)=\{x \in \bb_{v}:  q(x) \in p_v^nd_v^{-1}\} \supset O_\bv.$$
It follows that 
\begin{eqnarray*}
W_{a,v}(s,1,u)_{n<v(d_v)}
=  |d_v|^{\frac{5}{2}}  (1-N_v^{-s}) \sum_{n=0}^{v(d_v)-1}
N_v^{-n(s-1)}
=  |d_v|^{\frac{5}{2}}  
\frac{(1-N_v^{-s})(1-|d_v|^{s-1})}{1-N_v^{-(s-1)}}.
\end{eqnarray*}

A direct calculation of $W_{a,v}(s,1,u)_{n\geq v(d_v)}$ is quite involved, so we compare it with the unramified case instead. For clarification, denote 
$$D_n(a)^\circ=\{x \in \bv:  uq(x) \in a+p_v^n\},$$
which is equal to the set $D_n(a)$ in the unramified case in Step 1. 
For $n\geq v(d_v)$, the substitution $n\mapsto n+v(d_v)$ gives
\begin{eqnarray*}
W_{a,v}(s,1,u)_{n\geq v(d_v)}
=  |d_v|^{\frac{5}{2}} |d_v|^{s-1} (1-N_v^{-s}) \sum_{n=0}^{\infty}
N_v^{-n(s-1)} \frac{\vol(D_n(a)^\circ \cap O_\bv)}{\vol(O_\bv)}.
\end{eqnarray*}
This is equal to $W_{a,v}(s,1,u)$ in the case $|d_v|=1$ considered in Step 1. In other words, the result of Step 1 gives the combinatorial equality 
\begin{eqnarray*}
(1-N_v^{-s}) \sum_{n=0}^{\infty}
N_v^{-n(s-1)} \frac{\vol(D_n(a)^\circ \cap O_\bv)}{\vol(O_\bv)}
=\frac{(1-N_v^{-(s+2)})(1-N_v^{-(v(a)+1)(s+1)})}{1-N_v^{-(s+1)}}
\end{eqnarray*}
in the current setting. 
Hence, we have 
\begin{eqnarray*}
W_{a,v}(s,1,u)_{n\geq v(d_v)}
&=&  |d_v|^{s+\frac{3}{2}}\frac{(1-N_v^{-(s+2)})(1-N_v^{-(v(a)+1)(s+1)})}{1-N_v^{-(s+1)}}.
\end{eqnarray*}

Now we have a formula for $W_{a,v}(s,1,u)$, and some elementary computations
finish the proof of part (1) of the lemma. 

\medskip
\noindent\emph{Step 4. division case:} $v$ is nonsplit in $\BB$.
Then the formula in Step 2 becomes
\begin{eqnarray*}
W_{a,v}(s,1,u)
= -|d_v|^{\frac{1}{2}}  (1-N_v^{-s}) \sum_{n=0}^{\infty}
N_v^{-n(s-1)} \vol(D_n(a)\cap O_\bv^\times),
\end{eqnarray*}
where 
$$D_n(a)\cap O_\bv^\times=\{x \in O_\bv^\times:  uq(x) \in a+p_v^nd_v^{-1}\}.$$
Here we have assumed $u\in \ofv\cross$; otherwise, $W_{a,v}(s,1,u)=0$. 
Assume $v(a)\geq -v(d_v)$ by the same reason. 

If $v(a)>0$, the condition $uq(x) \in a+p_v^nd_v^{-1}$ is equivalent to 
$n\leq v(d_v)$.
In this case $D_n(a)\cap O_\bv^\times=O_\bv^\times$. 
Therefore,
\begin{eqnarray*}
W_{a,v}(s,1,u)
&=& -|d_v|^{\frac{1}{2}}  (1-N_v^{-s}) \sum_{n=0}^{v(d_v)}
N_v^{-n(s-1)} \vol(O_\bv^\times)\\
&=& -|d_v|^{\frac{1}{2}}  (1-N_v^{-s}) 
\frac{1-N_v^{(v(d_v)+1)(1-s)}}{1-N_v^{1-s}} \vol(O_\bv^\times).
\end{eqnarray*} 
This proves the first assertion in (2). 

It remains to verify the formula 
\begin{eqnarray*}
 \int_{F_v} \left(W_{a,v}'(0,1,u)-\frac12\log|a|_v W_{a,v}(0,1,u) \right)da=0.
\end{eqnarray*}

We first check that 
$$
\log|a|_v W_{a,v}(0,1,u)=0, \quad\ \forall a\in F_v^\times. 
$$
In fact, it suffices to check $W_{a,v}(0,1,u)=0$ if $v(a)\neq 0$. 
This is an easy consequence of the local Siegel--Weil formula in 
\cite[Theorem 2.2]{YZZ} or \cite[Proposition 2.9]{YZZ}. 
Alternatively, we can verify it by the type of computation here. 
Since we already know the vanishing for $v(a)>0$ from the above computation, 
it remains to check the case 
$-v(d_v) \leq v(a)< 0$. 
In this case, the condition $uq(x)-a \in p_v^nd_v^{-1}$ with $x\in O_\bv^\times$ is equivalent to $a \in p_v^nd_v^{-1}$. It holds only if $n\leq v(ad_v)$. Under this condition $D_n(a)\cap O_\bv^\times=O_\bv^\times$. 
It follows that 
\begin{eqnarray*}
W_{a,v}(s,1,u)
&=& -|d_v|^{\frac{1}{2}}  (1-N_v^{-s}) \sum_{n=0}^{v(ad_v)}
N_v^{-n(s-1)} \vol(O_\bv^\times)\\
&=& -|d_v|^{\frac{1}{2}} (1-N_v^{-s}) \frac{1-N_v^{(v(ad_v)+1)(1-s)}}{1-N_v^{1-s}}\vol(O_\bv^\times).
\end{eqnarray*}
Thus $W_{a,v}(0,1,u)=0$.

It is reduced to prove
\begin{eqnarray*}
 \int_{F_v} W_{a,v}'(0,1,u) da=0.
\end{eqnarray*}
We are going to compute 
\begin{eqnarray*}
 \int_{F_v} W_{a,v}(s,1,u) da
= -|d_v|^{\frac{1}{2}}  (1-N_v^{-s}) \sum_{n=0}^{\infty}
N_v^{-n(s-1)}  \int_{F_v} \vol(D_n(a)\cap O_\bv^\times) da.
\end{eqnarray*}
Use a Fuibini type of result to change the order of the last integral. We have
$$
\int_{F_v} \vol(D_n(a)\cap O_\bv^\times) da
=\iint_{uq(x)-a \in p_v^nd_v^{-1}}  dx da
=\int_{O_\bv^\times} \int_{uq(x)+p_v^nd_v^{-1}} da dx
=\vol(O_\bv^\times) |d_v|^{-\frac{1}{2}}N_v^{-n}.
$$
Hence,
\begin{eqnarray*}
 \int_{F_v} W_{a,v}(s,1,u) da
= -|d_v|^{\frac{1}{2}}  (1-N_v^{-s}) \sum_{n=0}^{\infty}
N_v^{-n(s-1)}  \vol(O_\bv^\times) |d_v|^{-\frac{1}{2}}N_v^{-n}
= -  \vol(O_\bv^\times) .
\end{eqnarray*}
Taking derivative at $s=0$, we get the desired result. 
The proof is complete. 
\end{proof}

\subsubsection*{Derivative of intertwining operators}

Recall from Theorem \ref{analytic series}(4), under the relation $\phi_{v}=\phi_{1,v}\otimes \phi_{2,v}$, we have defined
$$c_{\phi_v}(g,y,u)=r(g)\phi_{1,v}(y,u) W_{0,v}^{\circ}\ '(0,g,u, \phi_{2,v}) 
+ \log \delta(g_v)r(g)\phi_v(y,u), $$
where the normalization
\begin{eqnarray*}
W_{0,v}^{\circ}(s,g,u,\phi_{2,v})
= \gamma_{u,v}^{-1}|D_v|^{-\frac{1}{2}}|d_v|^{-\frac{1}{2}}\frac{L(s+1,\eta_v)}{L(s,\eta_v)}
 W_{0,v}(s,g,u,\phi_{2,v}).
\end{eqnarray*}
Here $\gamma_{u,v}$ is the Weil index of $(E_v\jv, uq)$. 

Assume all the conditions in \S\ref{choices}.
The following result is a variant of \cite[Lemma 7.6]{YZ}.

\begin{lem} \label{derivative of intertwining}
Denote $w=\matrixx{}{1}{-1}{}$.
For any non-archimedean place $v$ nonsplit in $\BB$, 
\begin{eqnarray*}
r(w)\phi_v(0,u)
=- |d_v|^2 N_v^{-1} (1-N_v^{-2}) 1_{\ofv\cross}(u),
\end{eqnarray*}
and
\begin{eqnarray*}
c_{\phi_v}(w,0,u)
=\begin{cases}
\displaystyle 
\frac{(1-N_v)\log N_v}{1+N_v}
r(w)\phi_{v}(0,u)  & \text{ if } 2\mid v(d_v),\\
0 & \text{ if }  2\nmid v(d_v).
\end{cases}
\end{eqnarray*}
\end{lem}

\begin{proof}
Note that $v$ is inert in $E$ by assumption, and 
$$
\phi_v=1_{O_{\bv}^\times \times \ofv\cross}, \quad\
\phi_{1,v}=1_{\oev^\times \times \ofv\cross}, \quad\
\phi_{2,v}=1_{\oev \jv \times \ofv\cross}.
$$
We need to compute
$$c_{\phi_v}(w,0,u)=r(w)\phi_{1,v}(0,u) W_{0,v}^{\circ}\ '(0,w,u, \phi_{2,v}).$$
It is easy to have 
$$
r(w)\phi_{1,v}(0,u)= \gamma(E_{v},uq)1_{\ofv\cross}(u)\int_{\oev^\times}dx
=\gamma(E_{v},uq) |d_v| (1-N_v^{-2}) 1_{\ofv\cross}(u),
$$
$$
r(w)\phi_{2,v}(x_2,u)= \gamma(E_{v}\jv,uq)
 |d_vq(\jv)| \cdot 1_{d_v^{-1}q(\jv)^{-1}\oev\jv}(x_2)1_{\ofv\cross}(u).
$$
This proves the first result, as $v(q(\jv))=1$ by assumption, and
$$
\gamma(E_{v},uq)\gamma(E_{v}\jv,uq)
=\gamma(\BB_{v},uq)=-1.
$$


Now we prove the second identity. 
From the definition
\begin{eqnarray*}
W_{0,v}(s, g,u, \phi_{2,v})= 
\int_{F_v} \delta(wn(b)g)^s  r(wn(b)g)\phi_{2,v}(0,u) db,
\end{eqnarray*}
we have 
\begin{eqnarray*}
W_{0,v}(s,w,u, \phi_{2,v})= W_{0,v}(s, 1,u, r(w)\phi_{2,v}).
\end{eqnarray*}
Its computation is similar to that of \cite[Lemma 7.6]{YZ}.
In fact, we still have 
\begin{eqnarray*}
W_{0,v}(s,1,u, r(w)\phi_{2,v})
= \gamma(E_{v}\jv,uq) |d_v|^{\frac{1}{2}} (1-N_v^{-s}) \sum_{n=0}^{\infty} N_v^{-ns+n}
\int_{D_n} r(w)\phi_{2,v}(x_2,u) d_ux_2,
\end{eqnarray*}
where
$$
D_n =\{x_2 \in E_v\jv:  uq(x_2) \in p_v^nd_v^{-1}\}
$$
and the measure $d_ux_2$ gives $\vol(\oev\jv)= |d_vuq(\jv)|$.
It follows that 
\begin{eqnarray*}
W_{0,v}^{\circ}(s,w,u,\phi_{2,v})
= \frac {1-N_v^{-2s} }{1+N_v^{-(s+1)}} \sum_{n=0}^{\infty} N_v^{-ns+n}
\int_{D_n} r(w)\phi_{2,v}(x_2,u) d_ux_2.
\end{eqnarray*}

Apply the formula of $r(w)\phi_{2,v}(x_2,u)$ in the above. 
Assume $u\in\ofv\cross$ in the following. 
Note that for any $n\geq 0$,
$$
D_{n} =p_v^{[\frac {n-v(d_v)}{2}]}\oev\jv\ \subset \
d_v^{-1}q(\jv)^{-1}\oev\jv.
$$
We have 
\begin{eqnarray*}
W_{0,v}^{\circ}(s,w,u,\phi_{2,v})
= \gamma(E_{v}\jv,uq)
 |d_vq(\jv)| \frac {1-N_v^{-2s} }{1+N_v^{-(s+1)}}
   \sum_{n=0}^{\infty} N_v^{-ns+n}
\vol(D_n).
\end{eqnarray*}
The summation is equal to 
\begin{eqnarray*}
\sum_{n=0}^{\infty} N_v^{-ns+n}
N_v^{-2[\frac {n-v(d_v)}{2}]} |d_vq(\jv)|
=\begin{cases}
N_v^{-1} (1+N_v^{1-s})(1-N_v^{-2s})^{-1} & \text{ if } 2\mid v(d_v),\\
(1+N_v^{-1-s})(1-N_v^{-2s})^{-1}& \text{ if }  2\nmid v(d_v).
\end{cases}
\end{eqnarray*}
Hence, 
\begin{eqnarray*}
W_{0,v}^{\circ}(s,w,u, \phi_{2,v})
=\gamma(E_{v}\jv,uq) |d_vq(\jv)| \cdot \begin{cases}
N_v^{-1} (1+N_v^{1-s})(1+N_v^{-1-s})^{-1} & \text{ if } 2\mid v(d_v),\\
1 & \text{ if }  2\nmid v(d_v).
\end{cases}
\end{eqnarray*}
Then 
\begin{eqnarray*}
W_{0,v}^{\circ}\ '(0,w,u, \phi_{2,v})= \gamma(E_{v}\jv,uq) |d_vq(\jv)|\cdot
\begin{cases}
\displaystyle
\frac{(1-N_v)\log N_v}{1+N_v}  & \text{ if } 2\mid v(d_v),\\
0 & \text{ if }  2\nmid v(d_v).
\end{cases}
\end{eqnarray*}
This finishes the proof. 
 
\end{proof}

\subsubsection*{An average formula}

Assume all the conditions in \S\ref{choices}.

Let $v$ be a non-archimedean place nonsplit in $\BB$.
Consider the difference
$$
\bar k_{\phi_v}(y,u)= k_{\phi_v}(1,y,u)-m_{\phi_v}(y,u)\log N_v, \quad (y,u)\in (B(v)_v-E_v)\times F_v^\times.
$$
Here $k_{\phi_v}(1,y,u)$ is defined in Theorem \ref{analytic series}(2).
On the other hand, the function $m_{\phi_v}(y,u)$ is defined in \cite[\S8.2, Notation 8.3]{YZZ}, which works for the settings of both \cite[\S8.2]{YZZ} (for $v$ nonsplit in $E$ but split $\BB$) and 
\cite[\S 8.3]{YZZ} (for $v$ nonsplit in both $E$ and $\BB$).
The definition is based on local intersection numbers on local integral models of the Shimura curve $X_U$ above $O_{F_v}$, and we refer to \cite[\S 8.2-8.3]{YZZ} and  \cite[\S 8.2]{YZ} for how to form a pseudo-theta series using $m_{\phi_v}(y,u)$ and $\phi^v$ to express arithmetic intersection numbers of CM points above $v$.

For our purpose, $\bar k_{\phi_v}(y,u)$
extends to a Schwartz function in $\ol\CS(B(v)_v\times F_v\cross)$, as a combination of \cite[Lemma 7.4, Lemma 8.7]{YZ}. This fits into the setting of nonsingular pseudo-theta series.
In the following, we compute the action of $w=\matrixx{}{1}{-1}{}$ on this Schwartz function. 

\begin{lem} \label{average k}
For any non-archimedean place $v$ nonsplit in $\BB$, 
$$
r(w) \bar k_{\phi_v}(0,u) =-r(w) \phi_{v}(0,u)
 \cdot
\begin{cases}
\displaystyle \left(\frac{N_v}{N_v+1}+\frac{v(d_v)}{2}\right) \log N_v , & 2\mid v(d_v); \\ 
\displaystyle   \frac{v(d_v)+1}{2}  \log N_v , & 2\nmid v(d_v).
\end{cases}
$$
\end{lem}

\begin{proof}
Write $y=y_1+y_2$ according to $B(v)_v=E_v+E_vj_v$ as usual. 
Note that $E_v$ is unramified over $F_v$, and $v(q(j_v))=0$.  
By \cite[Lemma 8.7]{YZ}, 
$$
 m_{\phi _v}(y, u)=
  \phi_{v}(y_1,u) 1_{\oev j_v}(y_2)
\cdot \frac{1}{2} v(q(y_2)).
$$
The function 
\begin{equation*}
k_{\phi_v}(1,y,u)= \frac{L(1,\eta_v)}{\mathrm{vol}(E_v^1)} 
\phi_{1,v}(y_1,u) {W_{uq(y_2),v}^\circ}'(0,1,u, \phi_{2,v}).  
\end{equation*}
is computed in \cite[Lemma 7.4]{YZ}. 
Here $\mathrm{vol}(E_v^1)=|d_v|^{\frac 12}$ in the current case. 
From the proof of \cite[Lemma 7.4]{YZ} (written in \cite[p. 596]{YZ}), which has also computed ${W_{a,v}^\circ}'(0,1,u)$ for $-v(d_v) \leq v(a) <0$, we have 
\begin{equation*}
k_{\phi_v}(1,y,u)= (\log N_v)\phi_{1,v}(y_1,u) \cdot 
\begin{cases}
\ 0,  & v(q(y_2))<-v(d_v); \\
\displaystyle \frac{N_v|q(y_2)|^{-1}-|d_v|}{N_v^2-1},  & -v(d_v) \leq v(q(y_2))<0;\\
\displaystyle \left(\frac{N_v-|d_v|}{N_v^2-1} + \frac{1}{2} v(q(y_2)) \right), \quad
 &  v(q(y_2))\geq 0.
\end{cases}
\end{equation*}
It follows that 
\begin{equation*}
\bar k_{\phi_v}(y,u)= (\log N_v)\phi_{1,v}(y_1,u) \phi_{2,v}'(y_2,u),
\end{equation*}
where $\phi_{2,v}'\in \ol\CS(E_vj_v \times F_v^\times)$ is given by
\begin{equation*}
\phi_{2,v}'(y_2,u) = 
 1_{\ofv\cross}(u) \cdot
\begin{cases}
\ 0,  & v(q(y_2))<-v(d_v); \\
\displaystyle \frac{N_v|q(y_2)|^{-1}-|d_v|}{N_v^2-1},  & -v(d_v) \leq v(q(y_2))<0;\\
 \displaystyle \frac{N_v-|d_v|}{N_v^2-1}  , \quad
 &  v(q(y_2))\geq 0.
\end{cases}
\end{equation*}

Now we compute
$$
r(w) \phi_{2,v}'(0,u)=
\gamma(E_vj_v, uq)
\int_{E_vj_v} \phi_{2,v}'(y_2,u)d_uy_2.
$$
Assume $u\in \ofv\cross$. 
Note that $E_v$ is unramified over $F_v$, $v(q(j_v))=0$, 
and $\vol(\oev j_v)=|d_v|$.  
Writing $v(q(y_2))=2i$ for $i\in \ZZ$, we have 
$$
r(w) \phi_{2,v}'(0,u)=
\gamma(E_vj_v, uq) 
\left(\sum_{-v(d_v)\leq 2i<0} \frac{N_v N_v^{2i}-|d_v|}{N_v^2-1} \vol(\varpi_v^i\oev^\times j_v)
+ \frac{N_v-|d_v|}{N_v^2-1} \vol(\oev j_v) \right),
$$
where $\varpi_v\in\ofv$ is a uniformizer. 
It follows that
$$
r(w) \phi_{2,v}'(0,u)=
\gamma(E_vj_v, uq) |d_v|
\left(\sum_{-v(d_v)\leq 2i<0} \frac{N_v N_v^{2i}-|d_v|}{N_v^2-1} (1-N_v^{-2}) N_v^{-2i} 
+ \frac{N_v-|d_v|}{N_v^2-1} \right).
$$
An elementary computation gives
$$
r(w) \phi_{2,v}'(0,u)=
\gamma(E_vj_v, uq)\cdot |d_v|
\cdot 1_{\ofv\cross}(u) \cdot
\begin{cases}
\displaystyle \frac{1}{N_v+1}+\frac{v(d_v)}{2N_v}, & 2\mid v(d_v); \\ 
\displaystyle   \frac{v(d_v)+1}{2N_v} , & 2\nmid v(d_v).
\end{cases}
$$

Note that
$$
r(w) \phi_{2,v}(0,u)=
\gamma(E_v\jv, uq)\cdot N_v^{-1}|d_v|
\cdot 1_{\ofv\cross}(u).
$$
It remains to check
$\gamma(E_vj_v, uq)=-\gamma(E_v\jv, uq)$ for the Weil indexes. 
This follows from 
$$
\gamma(E_v, uq)\gamma(E_v\jv, uq)
=\gamma(\bv, uq)=\gamma(\bv, q)
=-1,
$$
$$
\gamma(E_v, uq)\gamma(E_vj_v, uq)
=\gamma(B(v)_v, uq)=\gamma(B(v)_v, q)
=1.
$$
\end{proof}

\section{Height series} \label{sec height series}

The goal of this section is to decompose the height series $Z(g, (t_1, t_2))_U$ into a sum of pseudo-theta series and pseudo-Eisenstein series, and compute some related terms. This is mainly treated in \cite{YZZ,YZ}, but we do need to compute some extra terms for the purpose here.

\subsection{Shimura curve and height series} \label{sec shimura curve}

The goal of this subsection is to recall the basics of Shimura curves and the height series in \cite{YZZ,YZ}.

\subsubsection*{Shimura curve} 
 
Let $F$ be a totally real number field, and $\BB=\prod_v'\BB_v$ be a totally definite incoherent quaternion algebra over $F$ with ramification set $\Sigma$. 
Then $\Sigma$ is a finite set of places of $F$ of odd cardinality and containing all the archimedean places.

For any open compact subgroup $U$ of $\bfcross$, 
we have a Shimura curve $X_U$, which is a projective and smooth curve over 
$F$. For any embedding $\sigma: F\hookrightarrow \BC$, 
it has the usual complex uniformization
$$X_{U, \sigma}(\BC) = B(\sigma)^\times \bs \CH^\pm\times
\bb_f^\times / U \cup \{\rm cusps\}.$$
Here $B(\sigma)$ denotes the nearby quaternion algebra over $F$ with ramification set $\Sigma \setminus \{\sigma\}$, 
the action of $B(\sigma)^\times$ on $\CH^\pm$ is via a fixed isomorphism $B(\sigma)_\sigma\simeq M_2(\RR)$, 
and the action of $B(\sigma)^\times$ on $\bfcross$ is via a fixed 
 isomorphism $B(\sigma)_{\af}\to \BB_f$ over $\af$ as in \S\ref{choices}.

\subsubsection*{Integral models}

The exposition here agrees with that of \cite[\S4.2]{YZ}. 

Fix a maximal order $O_{\BB_f}=\prod_{v\nmid\infty} O_{\BB_v}$ of $\BB_f$. 
This gives a maximal open compact subgroup $O_{\BB_f}^\times$ of $\bfcross$. 
Let $U=\prod_v U_v$ be an open compact subgroup of $O_{\BB_f}^\times$ such that 
$U_v$ is maximal at every $v\in \Sigma_f$. 
Here we briefly recall the integral model $\CX_U$ of $X_U$ over $O_F$.

For any positive integer $N$, denote by $U(N)=(1+N O_{\BB_f})^\times$ the open compact subgroup of $O_{\BB_f}^\times$. 
Its local component is $U(N)_v=(1+N O_{\BB_v})^\times$. 

By the works of Deligne--Rapoport \cite{DR}, 
Carayol \cite{Ca}, and Boutot--Zink \cite{BZ}, 
if $U$ is sufficiently small in the sense that it is contained in $U(N)$ for some $N\geq 3$,  
there is a canonical integral model $\CX_{U}$ of $X_U$ over 
$O_{F}$, which is a projective and flat regular scheme over $O_{F}$.
It is further smooth at every place $v\notin\Sigma_f$ such that $U_v$ is maximal, and semistable at every place $v\in\Sigma_f$ (where $U_v$ is maximal by assumption).

For general $U$ (still with $U_v$ maximal for all $v\in \Sigma_f$), the
canonical integral model $\CX_{U}$ of $X_U$ over 
$O_{F}$ is constructed as follows. 
Take a sufficiently small open compact subgroup $U'$ of $\bb_f^\times$
such that $U_v'$ maximal for all $v\in \Sigma_f$ and that $U'$ is a normal subgroup of $U$. 
From the above paragraph, we have a canonical integral model $\CX_{U'}$
of $X_{U'}$ over $O_{F}$.
The finite group $U/U'$ acts on $X_{U'}$ via the right Hecke translation, and the quotient is exactly $X_U$. 
The action extends to an action of $U/U'$ on $\CX_{U'}$, so the quotient  scheme
$$\CX_{U}:=\CX_{U'}/(U/U')$$
is an integral model of $X_U$ over $O_F$. 

As in \cite[Corollary 4.6]{YZ},  $\CX_U$ is a normal and $\QQ$-factorial integral scheme, projective and flat over $O_{F}$. 
Here ``$\QQ$-factorial'' means that any Weil divisor has a positive integer multiple which is a Cartier divisor.
By Raynaud's abstract result in \cite[Appendice]{Ra}, $\CX_U$ is actually semistable over $O_F$, and smooth over every non-archimedean place $v$ split in $\BB$. 
The smoothness can be found at the top of page 195 of the loc. cit..
As a consequence, for any finite extension $H$ of $F$, 
the base change $\CX_U\times_{O_F}O_H$ is also flat and semistable, so it is still 
normal and $\QQ$-factorial.  
 
The integral model $\CX_U$
does not depend on the choice of $N$.
We remark that $\CX_U$ is also $\QQ$-Gorentein in the sense that
 there is a natural notion of dualizing sheaf $\omega_{\CX_{U}/O_F}$ as a  $\QQ$-line bundle on $\CX_U$. 
 In fact, denote by $\CX_U^{\mathrm{reg}}$ the regular locus of $\CX_U$. 
 The singular locus $\CX_U\setminus \CX_U^{\mathrm{reg}}$ is finite since $\CX_U$ is normal. 
Then the relative dualizing sheaf $\omega_{\CX_U^{\mathrm{reg}}/O_F}$ is a line bundle over $\CX_U^{\mathrm{reg}}$. 
Take a divisor $\CD$ on $\CX_U^{\mathrm{reg}}$ linearly equivalent to $\omega_{\CX_U^{\mathrm{reg}}/O_F}$, and take $\CD'$ to be the Zariski closure of $\CD$ in $\CX_U$. 
By $\QQ$-factoriality, some multiple of $\CD'$ is a Cartier divisor, and thus the same multiple of $\omega_{\CX_U^{\mathrm{reg}}/O_F}$ extends to a line bundle on $\CX_U$.

\subsubsection*{Arithmetic Hodge bundle}

The exposition here follows from that of \cite[\S4.2]{YZ}, which introduces the {canonical arithmetic Hodge bundle} $\CL_U$ over $\CX_U$. 

Recall that  the {Hodge bundle} of $X_U$,  as $\QQ$-line bundle over $X_U$, is defined by
$$
L_{U}= \omega_{X_{U}/F} \otimes\CO_{X_U}\Big(\sum_{Q\in X_U(\overline F)} (1-e_Q^{-1}) Q\Big).
$$
We refer to \S\ref{sec intro shimura curve} for an explanation of the ramification indexes.

As in \S\ref{sec intro shimura curve},  denote 
$$S(U)=\Spec\, O_F\setminus \{v: U_v \text{ is not maximal}\}.$$ 
Then $\CL_U$ is a $\QQ$-line bundle on $\CX_U$ constructed as follows. If $U$ is sufficiently small in the sense that $U$ is contained in $U(N)=(1+N O_{\bb_f})^\times$ for some integer $N\geq 3$, then over the open subscheme 
$\CX_{U, S(U)}=\CX_{U}\times_{\Spec\, O_F} S(U)$,
we have 
$$
\CL_U|_{\CX_{U, S(U)}}= \omega_{\CX_{U}/O_F} \otimes \CO_{\CX_{U, S(U)}}\Big(\sum_{Q\in X_U} (1-e_Q^{-1}) \CQ\Big).
$$
Here $\omega_{\CX_{U}/O_F}$ is the relative dualizing sheaf, the summation is through closed points $Q$ of $X_U$, $\CQ$ is the Zariski closure of $Q$ in $\CX_{U,S(U)}$, and $e_Q$ is the ramification index of any point of $X_U(\ol F)$ corresponding to $Q$.
If $U$ is a maximal open compact subgroup of $\bfcross$, for any sufficiently small normal open compact subgroup $U'$ of $\bfcross$ contained in $U$, we have
$$
\CX_{U, S(U')}= \CX_{U',S(U')}/(U/U'), \quad
\CL_{U}|_{\CX_{U, S(U')}}= \deg(\pi)^{-1}N_{\pi} (\CL_{U'}|_{\CX_{U', S(U')}}),
$$
where
$N_{\pi}:\Pic(\CX_{U', S(U')})\to \Pic(\CX_{U, S(U')})$ is the norm map with respect to the natural map $\pi: \CX_{U', S(U')}\to \CX_{U, S(U')}$. 
Varying $U'$, we glue $\{\CL_{U}|_{\CX_{U, S(U')}}\}_{U'}$ together to form the $\QQ$-line bundle $\CL_U$ over $\CX_U$ for maximal $U$.
For general $U$, we take an embedding $U\subset U_0$ into a maximal $U_0$, and than define $\CL_U$ to be the pull-back of $\CL_{U_0}$ via the natural map $\CX_U \to \CX_{U_0}$.

At any archimedean place $\sigma:F\to \BC$, 
 the {Petersson metric} of $\CL_U$ is given by 
$$\|f(\tau)d\tau\|_{\mathrm{Pet}}=2\, \Im(\tau) |f(\tau)|,$$
where $\tau$ is the standard coordinate function on $\CH$, and $f(\tau)$ is any meromorphic modular form of weight 2 over $X_{U,\sigma}(\CC)$. 
Thus we have {the arithmetic Hodge bundle}
$$\ol\CL_U=(\CL_U, \{\|\cdot\|_\sigma\}_\sigma).$$
It is a hermitian $\QQ$-line bundle over $\CX_U$. 

As in \cite[Theorem 4.7]{YZ}, the system $\{\ol\CL_U\}_U$ is compatible with pull-back morphisms. 
In other words, if $U'\subset U$ are two open compact subgroups, then the natural morphism $\pi:\CX_{U'}\to \CX_U$ gives $\pi^*\ol\CL_U= \ol\CL_{U'}$ in $\wh\Pic(\CX_{U'})_\QQ$. 
Moreover, the system is also compatible with the norm map in the sense that 
$N_{\pi}(\ol\CL_{U'})=\deg(\pi)\ol\CL_{U}$ in $\wh\Pic(\CX_{U})_\QQ$. This is used as the construction in \cite[Theorem 4.7]{YZ}, but can also be deduced by the identity $N_{\pi}\circ \pi^*=\deg(\pi)$ as operators in $\wh\Pic(\CX_{U})_\QQ$.

Note that $X_U$ is connected by not geometrically connected over $F$. 
For any connected component $X_\alpha$ of $X_{U, \overline F}$, denote by 
$L_\alpha=L_{U,\alpha}$ the pull-back of $L_U$ to $X_\alpha$. 
Denote by $\kappa_U^\circ$ the degree of $L_\alpha$ over $X_\alpha$, which is independent of $\alpha$ by the Galois action. 
Denote 
$$\xi_U=(\kappa_U^\circ)^{-1}L_U, \quad
\hat\xi_U=(\kappa_U^\circ)^{-1}\ol\CL_U, \quad
\xi_\alpha=(\kappa_U^\circ)^{-1}L_\alpha,$$
which are respectively a $\QQ$-line bundle over $X_U$, a hermitian $\QQ$-line bundle over $\CX_U$, and a  $\QQ$-line bundle over $X_\alpha$.

\subsection*{Quotient at split places}

If $F=\QQ$ and $U$ is maximal, by \cite[Lem. 2.1]{Yu2}, we have
$$
\CL_U= \omega_{\CX_{U}/O_F} \otimes \CO_{\CX_{U}}\Big(\sum_{Q\in X_U} (1-e_Q^{-1}) \CQ\Big),
$$
where $\CQ$ is the Zariski closure of $Q$ in $\CX_{U}$.
This gives an equivalent definition of $\CL_U$ in this case.

If $F$ is totally real, the result does not hold over $O_F$, which is the reason why we have to take the above intricate definition of $\CL_U$. However, the following result asserts that the result still holds at places $v$ split in $\BB$. 
It will not be used in elsewhere of this paper, but we include it here for completeness.

\begin{prop}
Let $U'\subset U$ be a normal open subgroup which is sufficiently small as an open compact subgroup of $\bb_f^\times$. 
Let $v$ be a non-archimedean place of $F$ split in $\BB$ such that $U$ and $U'$ are maximal at $v$. Denote by $k$ the residue field of $v$ in $O_F$.
Then the following hold.
\begin{itemize}
\item[(1)] 
The fibers $\CX_{U',k}$ and $\CX_{U, k}$ are smooth over $k$, and the morphism $\CX_{U',k}\to \CX_{U, k}$ induces an isomorphism $\CX_{U',k}/(U/U')\to \CX_{U, k}$. 
\item[(2)] 
The quotient morphism $\CX_{U'}\to \CX_{U}$ is \'etale at every generic point of the fiber $\CX_{U, k}$ of $\CX_{U}$ above $v$.
\item[(3)] Over $S=\Spec\, O_F\setminus \Sigma_f$, there is a canonical isomorphism
$$
\CL_{U, S}= \omega_{\CX_{U,S}/S} \otimes \CO_{\CX_{U,S}}\Big(\sum_{Q\in X_U} (1-e_Q^{-1}) \CQ\Big),
$$
where $\CQ$ is the Zariski closure of $Q$ in $\CX_{U,S}$.
\end{itemize}

\end{prop}

\begin{proof}

By the uniformization, an element $g\in U$ acts trivially on $X_{U'}$ if and only if 
$g$ lies in $(F^\times U')\cap U=(F^\times \cap U)U'=O_F^\times U'$.   It follows that the action induces a faithful action of the finite group 
$G=U/(O_F^\times U')$ on $X_{U'}$.

We need the basic theory of quotient of schemes by finite groups.
We refer to \cite[p. 66, Thm]{Mu} for the case of varieties, where the 
current case is similar.
In particular, the map $\CX_{U'}\to \CX_{U}$ is a topological quotient in the sense that the underlying topological space of $\CX_U$ is exactly the quotient of the underlying topological space of $\CX_{U'}$ by $G$. 
Moreover, if $\CX_{U'}$ is a union of affine open subschemes $\mathrm{Spec}(A)$ stabilized by $G$, then the quotient $\CX_{U}$ is the union of affine open subschemes 
$\mathrm{Spec}(A^G)$. 
As $\CX_{U'}$ is regular, the affine ring $A$ is normal, so $A^G$ is normal by definition. 
This implies that $\CX_U$ is normal.

By the action of $G$ on $\CX_{U',k}$, we also have a quotient $\CX_{U',k}/G$. 
There is a natural morphism $\CX_{U',k}/G \to \CX_{U,k}$, induced by the homomorphism
$A^G\otimes_{O_{F,(v)}} k \to (A\otimes_{O_{F,(v)}} k)^G$ from the above description of the quotient process. 
The morphism $\CX_{U',k}/G \to \CX_{U,k}$ is a homeomorphism, since the quotients are topological quotients. 
As $\CX_{U',k}$ is smooth over $k$, the quotient $\CX_{U',k}/G$ is also smooth over $k$. 
By \cite[Appendice, p. 195]{Ra}, $\CX_{U,k}$ is also smooth. 
This forces the homeomorphism $\CX_{U',k}/G \to \CX_{U,k}$ to be an isomorphism.
This proves (1). 

Now we see that (2) implies (3). 
By varying $U'$, it suffices to prove the isomorphism of (3) over the local ring $O_{F,(v)}$ of $O_{F}$ at $v$.
Denote $\CX_{U,(v)}=\CX_U\times_{O_F}O_{F,(v)}$ and $\CX_{U',(v)}=\CX_{U'}\times_{O_F}O_{F,(v)}$ for convenience.
As above, we already know that $\CX_{U,(v)}$ is smooth over $O_{F,(v)}$. 
Once (2) holds, the ramification divisor of $\pi:\CX_{U',(v)}\to \CX_{U,(v)}$ is horizontal, thus it is just the Zariski closure of the ramification divisor of $\CX_{U'}\to \CX_{U}$. 
This proves that in the case $|\Sigma_f|>1$, 
$$
N_\pi(\omega_{\CX_{U',(v)}/O_{F,(v)}})
\simeq \omega_{\CX_{U,(v)}/O_{F,(v)}} \otimes \CO_{\CX_{U,(v)}}\Big(\sum_{Q\in X_U} (1-e_Q^{-1}) \CQ_{O_{F,(v)}}\Big);
$$
in the case $|\Sigma_f|=1$, the formula holds after suitably modifying by cusps. 
This proves (3). 

Now we prove (2), which is the essential part of the theorem. 
Let $\xi$ be a generic point of (an irreducible component of) $\CX_{U',k}$, and denote by $k(\xi)$ the residue field of $\xi$. 
Denote the stabilizers
$$G_\xi^+=\{g\in G:g(\xi)=\xi\}, \quad
G_\xi=\{g\in G_\xi^+: g \text{\ acts trivially on } k(\xi)\}.$$ 
We claim that $G$ acts freely at $\xi$ in the sense that $G_\xi$ is trivial. 

We first see that the claim implies that $\pi:\CX_{U'}\to \CX_{U}$ is \'etale at $\xi$. 
In fact, denote by $\eta$ the image of $\xi$ in $\CX_U$.
Denote by $G\xi\subset \CX_{U'}$ the orbit of $\xi$, and denote $k(G\xi)=\oplus_{\xi'\in G\xi} k(\xi')$. 
Note that there is a homomorphism $k(\eta)\to k(G \xi)^G\simeq k(\xi)^{G_\xi^+}$. 
By assumption, $G_\xi^+$ acts freely on $k(\xi)$, so the index $[k(G \xi):k(\xi)^{G_\xi^+}]=|G|$. 
It follows that $[k(G \xi):k(\eta)]\geq |G|$.
On the other hand, as $\CX_U$ is normal, the local ring of $\CX_{U}$ at $\eta$ is a discrete valuation ring, and thus it is easy to see that $\pi:\CX_{U'}\to \CX_{U}$ is flat over $\eta$. As a consequence, we have $\deg(\pi^{-1}(\eta)/\eta)=|G|$. 
Note that $\mathrm{Spec}\, k(G \xi)$ is the reduced structure of $\pi^{-1}(\eta)$.
The bound of the degree forces $\pi^{-1}(\eta)=\mathrm{Spec}\, k(G \xi)$,   
$k(\eta)\simeq k(\xi)^{G_\xi^+}\simeq k(G \xi)^G$, and that $k(\xi)$ is Galois over $k(\eta)$. 
As a consequence, $\pi^{-1}(\eta)$ is unramified over $\eta$, and thus $\CX_{U'}\to \CX_{U}$ is \'etale at $\xi$. 

It remains to prove that $G$ acts freely at every generic point of $\CX_{U',k}$.
It suffices to prove that $G$ acts freely at every generic point of  
$\CX_{U',\bar k}$.

Our proof is based on the theory of Carayol \cite{Ca}. 
In particular, \cite[\S 6]{Ca} introduces a universal $p$-divisible group $E_\infty$ of height two over $\CX_{U'}$. 
As in \cite[\S 6.7]{Ca}, every geometric point of $\CX_{U',\bar k}$ is either supersingular or ordinary. 
By \cite[\S 11.2]{Ca}, the set of supersingular points of $\CX_{U',\bar k}$ can be expressed as
$$
\CX_{U',\bar k}^{\rm ss}\simeq B^\times \backslash  \ZZ\times B_{\BA_f^v}^\times /U'^v.
$$
Here $B=B(v)$ is the nearby quaternion algebra over $F$ by switching the invariant of $\BB$ at $v$, and $B^\times$ acts on $\ZZ$ via the canonical isomorphism 
$\tilde v:B_v^\times/O_{B_v}^\times\to \ZZ$ via the isomorphism, where $O_{B_v}$ is the unique maximal order of the division algebra $B_v$. 
 
Fix a supersingular point $x$ of $\CX_{U',\bar k}$, and assume that $x$ lies on an irreducible component $C$ of $\CX_{U',\bar k}$. 
By \cite[\S 6.6]{Ca}, the deformation of $E_\infty|_x$ over the category of local Artinian rings with residue field $\bar k$
is represented by the completion 
$\widehat \CO_{\CX', x}$ of the local ring of 
$\CX'=\CX_{U',O_{F_v}^{\rm ur}}$ at $x$. 
As $\widehat \CO_{C, x}=\widehat \CO_{\CX', x}/\wp$ for the maximal ideal of $\wp$ of $O_{F_v}^{\rm ur}$, the deformation of $E_\infty|_x$ over the category of local Artinian $\bar k$-algebras with residue field $\bar k$
is represented by  
$\widehat \CO_{C, x}$. 
Moreover, the base change $E_\infty|_{\widehat \CO_{C, x}}$ is the universal $p$-divisible group. 

By definition, the deformation induces an action of 
$\Aut(E_\infty|_x)\simeq O_{B_v}^\times$ on $\widehat \CO_{C, x}$, whose kernel is exactly $\Aut(E_\infty|_{\widehat \CO_{C, x}})$. 
A key property is that $\Aut(E_\infty|_{\widehat \CO_{C, x}})\simeq O_{F_v}^\times$. 
In fact, denote by $M$ the algebraic closure of the fraction field of $\widehat \CO_{C, x}$. 
Note that $\End(E_\infty|_{\widehat \CO_{C, x}})$ has natural injections into both 
$\End(E_\infty|_M)$
and $\End(E_\infty|_x)\simeq O_{B_v}$ as $O_{F_v}$-algebras. 
The $p$-divisible group $E_\infty|_M$ is ordinary, since the supersingular locus 
$\CX_{U',\bar k}^{\rm ss}$ is a finite set, and thus  $\End(E_\infty|_M)\simeq O_{F_v}\oplus O_{F_v}$ by splitting $E_\infty|_M$ into a direct sum of a local part and an \'etale part. 
These force $\End(E_\infty|_{\widehat \CO_{C, x}})=O_{F_v}$, and thus 
$\Aut(E_\infty|_{\widehat \CO_{C, x}})\simeq O_{F_v}^\times$. 
Hence, we conclude that the kernel of $O_{B_v}^\times \to \Aut(\widehat \CO_{C, x})$ is 
$O_{F_v}^\times$.

Now we are ready to prove that $G$ acts freely at every generic point of $\CX_{U',k}$. 
Let $C$ be an irreducible component of $\CX_{U',\bar k}$, and let $\xi_C$ be its generic point.
Let $x$ be a supersingular point of $C$.
To prove that $G$ acts freely on $\xi_C$, it suffices to prove that $G_x$ acts freely on the local ring $\CO_{C, x}$ of $x\in C$, so it suffices to prove that 
$G_x$ acts freely on the completion $\widehat \CO_{C, x}$ of $\CO_{C, x}$. 
Assume that $x$ is represented by $(m_x,b_x)\in \ZZ\times B_{\BA_f^v}^\times$ in the above coset expression of $\CX_{U',\bar k}^{\rm ss}$. 
An element $g\in U$ fixes $x$ if and only if 
$$B^\times (m_x,b_x U'^v) 
=B^\times (m_x,b_x (g^v)^{-1} U'^v),$$
which holds if and only if 
there is $\gamma\in B^\times$ such that 
$$
\tilde v(\gamma)=0,\quad
g^v\in b_x^{-1} \gamma b_x  U'^v.
$$
In this case, $g$ acts on $\widehat \CO_{C, x}$ via $\gamma\in O_{B_v}^\times$, and the action of $O_{B_v}^\times$ on $\widehat \CO_{C, x}$ is the one described above in terms of the deformation theory. 
In particular, the kernel of  $O_{B_v}^\times \to \Aut(\widehat \CO_{C, x})$ is 
$O_{F_v}^\times$. 

Assume that $g$ acts trivially on $\widehat \CO_{C, x}$. This implies 
$\gamma\in O_{F_v}^\times$ by the kernel identity. 
Recall that $\gamma\in B^\times$ is a global element, so $\gamma\in F^\times$.
Then the above relation of $\gamma$ and $g$ becomes
$$
\tilde v(\gamma)=0,\quad
g^v\in \gamma U'^v.
$$
It follows that $g\in F^\times U'$, and thus 
$g$ lies in $(F^\times U')\cap U=O_F^\times U'$. 
Then the image of $g$ in
$G=U/(O_F^\times U')$ is trivial. 
This proves that $G_x$ acts freely on $\widehat \CO_{C, x}$, and thus $G$ acts freely on $\xi_C$.
The proof is complete.
\end{proof}

\subsubsection*{CM points}
Let $E/F$ be a totally imaginary quadratic extension, 
with a fixed embedding $E_\BA\hookrightarrow \BB$ over $\BA=\BA_F$.
As in \S\ref{choices},
for any archimedean place $\sigma:F\hookrightarrow \CC$, we also fix an embedding $E\hookrightarrow B(\sigma)$ compatible with 
$E_\BA\hookrightarrow \BB$ and $B(\sigma)_{\af}\to \BB_f$. 
Then $E^\times$ acts on $\CH^\pm$ via $E\hookrightarrow B(\sigma)$.

Fix an embedding $\sigma_0:F\hookrightarrow \CC$. 
Let  $z_0\in\CH$ be the unique fixed point of $E\cross$ on $\CH$, where the action is via $E\to B(\sigma_0)$.
Via the complex uniformization, we have a CM point 
$$
[\beta]_U=[z_0,\beta]_U
$$ 
on $X_{U, \sigma_0}(\BC)$ for any $\beta\in \bfcross$. 
By definition of the canonical model, $[\beta]_U$ descends to an algebraic 
point of $X_U$ defined over $E^\ab$. 
Here $E^\ab$ is the maximal abelian extension of $E$, endowed with embeddings 
$F\to E\to E^\ab\to\ol F\to \CC$ refining $\sigma_0$. 
We may also abbreviate $[\beta]_U$ as $[\beta]$ or just $\beta$.

In particular, we have the CM point 
$$P=P_U=[z_0,1]_U.$$ 
We will view the finite abelian group
$$C_U=E\cross\bs E^\times_{\af} / (U\cap E^\times_{\af})$$
as a set of CM points via
the natural injection 
$$
C_U \lra X_U(E^\ab), \quad
t\longmapsto [t]_U.
$$
Note that the point $[t]_U$ depends only on the class of $t$ in $C_U$.

For any CM point $[\beta]_U$ with $\beta\in \bfcross$, assume that it lies in the connected component $X_\alpha$ of $X_{U,\overline F}$. 
Then the divisor 
$$[\beta]_U^\circ:=[\beta]_U-\xi_\alpha$$
is a divisor of degree 0 on $X_\alpha$, 
also viewed as a divisor on $X_{U,\overline F}$ by push-forward. 
We usually abbreviate $\beta^\circ=[\beta]_U^\circ$. 
By abuse of notations, we may also write $\xi_\alpha$ as $\xi_\beta$. 

Note the set of connected components of $X_{U, \overline F}$ is isomorphic to the group $F_+^\times\bs\afcross/q(U)$, and under this isomorphism we have $\alpha=q(\beta)$
(assuming $[\beta]_U$ lies in $X_\alpha$).

\subsubsection*{Height series}

Finally, we recall the generating series in \cite[\S3.4.5]{YZZ}. 
Let $\phi\in \ol\CS (\BB\times \BA^\times)$ be a Schwartz function invariant under $K=U\times U$; i.e., for any $(h_1,h_2)\in U\times U$, the Weil representation $r(h_1,h_2)\phi=\phi$ in the sense of \S\ref{sec theta eisenstein}.
More precisely, this means that 
$$
\phi(h_1^{-1}xh_2 , q(h_1)q(h_2)^{-1}u)=\phi(x,u), \quad 
(x,u)\in \BB\times \across.
$$
Then the generating series is defined by
\begin{equation*}  
Z(g,\phi)_U=Z_0(g,\phi)_U+Z_*(g,\phi)_U, \quad g\in \gla,
\end{equation*}
where
\begin{eqnarray*}
Z_0(g,\phi)_U&=& -  \sum_{\alpha\in F_+\cross\bs \afcross/q(U)} \sum_{u\in
\mu_U^2 \bs F\cross}
  E_0(\alpha^{-1}u,r(g)\phi)\ L_{K,\alpha},\\
  Z_*(g,\phi)_U&=& w_U \sum_{a\in F^\times} \sum_{x\in U\bs \bb_f^\times/U} r(g)\phi(x,aq(x)^{-1})\ Z(x)_U.
\end{eqnarray*}
Here $\mu_U=F^\times \cap U$, and $w_U=|\{1,-1\}\cap U|$ is equal to 1
or 2.
For any $x\in \bfcross$, 
$Z(x)_U$ denotes the Hecke correspondence on $X_U$ determined by the double coset $UxU$, which is also viewed as a divisor on $X_U\times X_U$. 
By \cite[Theorem 3.17]{YZZ}, $Z(g,\phi)_U$
is an automorphic form of $g\in\gla$ with coefficients in $\Pic(X_U\times X_U)_\CC$.

Let $E/F$ be a totally imaginary quadratic extension, with a fixed embedding
$E_\BA\hookrightarrow \BB$ over $\BA$. 
Recall from \cite[\S3.5.1, \S5.1.2]{YZZ} and \cite[\S8.1]{YZ} that we have a height series
\begin{eqnarray*}
Z(g, (t_1, t_2), \phi)_U
= \pair{Z(g, \phi)_U\ t_1^{\circ},\ t_2^{\circ} }_{\NT},
\quad t_1,t_2\in E^\times(\af).
\end{eqnarray*}
Here $Z(g, \phi)_U$ acts on $t_1^{\circ}$ as a correspondence, and the Neron--Tate height  over $F$ is defined as in \cite[\S7.1.2]{YZZ}. 
By the modularity, $Z(g, (t_1, t_2), \phi)_U$ is an automorphic form in $g\in \gla$. 
By \cite[Lemma 3.19]{YZZ}, it is actually a cusp form. 
In particular, the constant term $Z_0(g,\phi)$ of the generating function plays no role here.

Throughout this section, we will assume the basic conditions in \S \ref{choices}. In particular, $\phi$ is standard at every archimedean place $v$. 
It follows that $Z(g, \phi)_U$ and $Z(g, (t_1, t_2), \phi)_U$ are  automorphic forms in $g\in \gla$, holomorphic of parallel weight two at archimedean places. 
We will assume the restrictive conditions when it comes to explicit calculations.

\subsection{Weakly admissible extensions} \label{sec admissible}

In order to decompose the height series in terms of the arithmetic Hodge index theorem of Faltings--Hriljac, the notion of admissible extensions are used in \cite{YZ, YZZ}.
However,  there is a minor mistake involving misconceptions about admissible extensions in \cite{YZ, YZZ}. 
In fact, the Green's function is not admissible, but only weakly admissible in the current sense. 
As we will see, this mistake does not affect the main results of \cite{YZ, YZZ}, but it does affect the results here. 
In the following, we review the admissibility notion as described in \cite[\S 7.1-7.2]{YZZ}, introduce the weak admissibility notion in the mean time, and then point out the mistake and the correction.

\subsubsection*{Terminology for arithmetic intersection theory}

Here we review and modify some terminology of 
\cite[\S7.1.3-7.1.5]{YZZ} and make some additional definitions in the following.
Our setting is slightly more general than the loc. cit. by allowing the arithmetic surface to be $\QQ$-factorial (instead of being integral and semistable). 
Note that the integral model $\CX_U$ of the Shimura curve is $\QQ$-factorial as in the last subsection, and it remains so by reasonable base changes.

Let $X$ be a projective and smooth curve over a number field $F$.
By taking the definitions over every connected component, we can assume that $X$ is connected, but we do not assume that it is geometrically connected. 
Denote by $F'$ the algebraic closure of $F$ in the function field of $X$, so that $X$ is geometrically connected over $F'$.

Let $\CX$ be a projective, flat, normal and $\QQ$-factorial integral model of $X$ over $O_F$. By $\QQ$-factoriality, local intersection multiplicities of properly intersecting Weil divisors can be defined as rational numbers, so we can still consider arithmetic intersection theory over $\CX$. In the following, a divisor on $\CX$ means a Weil divisor on $\CX$, and the finite part of an arithmetic divisor on $\CX$ is allowed to be a Weil divisor. By linear combination of rational coefficients, we also have the notion of $\QQ$-divisors and arithmetic $\QQ$-divisors.

Let $\wh \CD_1=(\CD_1, g_1)$ and $\wh \CD_2=(\CD_2, g_2)$ be arithmetic 
 $\QQ$-divisors on $\CX$. Here as a convention, $\CD_i$ is a $\QQ$-divisor on 
 $\CX$, and $g_i$ is a Green's function of $\CD_i(\CC)$ on $X(\CC)=\coprod_{\sigma:F\to \CC}  {X_{\sigma}(\CC)}$.
Note that $X_{\sigma}(\CC)$ is not connected unless $F'=F$, but this does not affect our exposition.
Recall that the arithmetic intersection number is defined as 
$$
\wh \CD_1\cdot \wh \CD_2
=\wh \CD_1\cdot  \CD_2
+  \sum_{\sigma:F\to \CC}  \int_{X_{\sigma}(\CC)} g_{2,\sigma}\, c_1(\CD_{1,\sigma}, g_{1,\sigma}), 
$$
where the extra subscripts $\sigma$ indicate base change or restriction to the compact Riemann surface $X_{\sigma}(\CC)$, and
$c_1(\CD_{1,\sigma}, g_{1,\sigma})$ denotes the Chern form of the hermitian line bundle associated to $(\CD_{1,\sigma}, g_{1,\sigma})$ on $X_{\sigma}(\CC)$. 

It remains to explain the intersection number $\wh \CD_1\cdot \CD_2$, which is of  independent importance. 
If $\CD_1$ intersects $\CD_2$ properly on $\CX$,  we have 
$$
\wh \CD_1\cdot \CD_2=(\CD_1\cdot \CD_2)+ \sum_{\sigma:F\to \CC} 
g_{1,\sigma}(\CD_{2,\sigma}(\CC)), 
$$
where the finite part 
$$
\CD_1\cdot \CD_2
=\sum_{v\nmid\infty} (\CD_1\cdot \CD_2)_{v}
$$ 
is the usual intersection number on $\CX$ decomposed in terms of non-archimedean places $v$ of $F$, and the infinite part 
$g_{1,\sigma}(\CD_{2,\sigma}(\CC))$ is understood as $\sum_i a_ig_{1,\sigma}(z_i)$ if $\CD_{2,\sigma}(\CC)=\sum_i a_i z_i$ as a divisor on $X_{\sigma}(\CC)$. 
In general, there is an arithmetic $\QQ$-divisor $\wh \CD_1'=(\CD_1', g_1')$ on $\CX$ linearly equivalent to $\wh \CD_1$ such that  
 $\CD_1'$ intersects  $\CD_2$ properly, and then we set 
$\wh \CD_1\cdot \CD_2=\wh \CD_1'\cdot \CD_2.$
The result is independent of the choice of $\wh \CD_1'$.

Note that $\CX$ is actually a scheme over $\Spec\, O_{F'}$ since it is normal.
We can also group the above intersection numbers in terms of places over $F'$. For example, if $\CD_1$ intersects $\CD_2$ properly on $\CX$ as above,  then 
$$
\wh \CD_1\cdot \CD_2=(\CD_1\cdot \CD_2)+ \sum_{\sigma':F'\to \CC} 
g_{1,\sigma'}(\CD_{2,\sigma'}(\CC)), 
$$
where the finite part 
$$
\CD_1\cdot \CD_2
=\sum_{v'\nmid\infty} (\CD_1\cdot \CD_2)_{v'}
$$ 
is a sum over non-archimedean places $v'$ of $F'$.

\subsubsection*{Flat arithmetic divisors}

An arithmetic $\QQ$-divisor $\wh \CD=(\CD,g_\CD)$ on $\CX$ is said to be \emph{flat} if its intersection number with any vertical arithmetic divisor is 0, or equivalently the following two
conditions hold:
\begin{itemize}
\item[(a)] The Chern form  $c_1(\CD_{\sigma}(\CC), g_{\CD,\sigma})$ on $X_{\sigma}(\BC)$ is 0 for any
embedding $\sigma:F\to \CC$;
\item[(b)] The intersection number $\CD \cdot \CV=0$ for any irreducible component $\CV$ of any closed fiber of $\CX$ above $\Spec\,O_{F}$. 
\end{itemize}
These imply that $\CD$ has degree zero on $X$ over $F$. 
The notion ``flat'' depends only on the arithmetic divisor class of $\wh \CD$, so it is naturally a property for arithmetic divisor classes or hermitian line bundles.

The Hodge index theorem of Faltings and Hriljac (cf. \cite[Theorem 7.4]{YZZ}) also holds here. Namely, let $\wh \CD_1=(\CD_1, g_1)$ and $\wh \CD_2=(\CD_2, g_2)$ be two flat arithmetic $\QQ$-divisors on $\CX$, then 
$$
\wh\CD_1\cdot \wh\CD_2=-\pair{\CD_{1,F},\CD_{2,F}}_\NT.
$$
In fact, let $\CX'\to \CX$ be a resolution of singularity. The pull-back of $\wh\CD_1$ and 
$\wh\CD_2$ are still flat over $\CX'$, so the result follows from that on $\CX'$.

\subsubsection*{Weakly admissible extensions}

Resume the above notations for $(X, \CX, F, F')$. 
Now we introduce our key admissibility notion. 

Fix an arithmetic divisor
class $\hat\xi \in \wh{\Pic}(\CX)_\QQ$ whose generic fiber has degree 1 on $X$ over $F'$.
Let $\wh \CD=(\CD,g_\CD)$ be an arithmetic $\QQ$-divisor on $\CX$. We can always write
$\CD=\CH+\CV$ where $\CH$ is the horizontal part of $\CD$, and $\CV$ is the vertical part of
$\CD$. The arithmetic
divisor $\wh \CD$ is called \emph{$\hat\xi$-admissible} if the following
conditions hold:
\begin{itemize}
\item[(1)] The difference $\wh \CD-\deg_{F'} (\CD_{F'})\cdot \hat\xi$ is flat over $\CX$;
\item[(2)] The intersection number $(\CV\cdot \hat\xi)_{v'}=0$ for any non-archimedean place $v'$
of $F'$;
\item[(3)] The integral $\ds\int_{X_{\sigma'}(\BC)} g_\CD c_1(\hat\xi)=0$  for any
embedding $\sigma':F'\to \CC$.
\end{itemize}
The arithmetic divisor $\wh \CD$ is called \emph{weakly $\hat\xi$-admissible} if it satisfies conditions (1) and (2). 

A \emph{$\hat\xi$-admissible extension} (resp.  \emph{weakly $\hat\xi$-admissible extension}) of a $\QQ$-divisor $D_0$ over $X$ is a {$\hat\xi$-admissible} (resp.  \emph{weakly $\hat\xi$-admissible}) arithmetic $\QQ$-divisor 
$\wh \CD=(\CD,g_\CD)$ over $\CX$ such that the generic fiber $\CD_{F}=D_0$ over $X$.

The notion \emph{$\hat\xi$-admissible} is introduced in \cite[\S7.1.5]{YZZ}, while the notion 
\emph{weakly $\hat\xi$-admissible} is a new one added here.
Note that the {$\hat\xi$-admissible} extension exists and is unique.
On the other hand, without condition (3), condition (1) only determines the Green's function up to constant functions over $X_{\sigma'}(\CC)$. 

Nonetheless, in our calculation over the Shimura curve, we do have a fixed choice of Green's functions as follows. 
To illustrate the idea, we will specify a symmetric and smooth function $g: X(\CC)\times X(\CC)\setminus \Delta\to \RR$ such that for any $P\in X(\CC)$, the 1-variable function $g(P,\cdot)$ is a Green's function for the divisor $P$ over $X(\CC)$ with curvature form equal to $c_1(\hat\xi)$. 
Then for any divisor $D$ over $X$, we take the Green's function $g_D=g(D,\cdot)$. 

Let $D_1, D_2$ be two $\QQ$-divisors over $X_{\ol F}$. 
Then $D_1, D_2$ are realized as $\QQ$-divisors over $X_{L}$ for a finite extension $L$ of $F$. 
Assume that $\CX_{O_L}$ is still $\QQ$-factorial. 
By abuse of notations, still denote the pull-back of $\hat\xi$ to $\CX_{O_L}$ by $\hat\xi$.

For $i=1,2$, let 
$\wh D_i=(\ol D_i+\CV_i, g_i)$ be a weakly $\hat\xi$-admissible extension of $D_i$ over $\CX_{O_L}$.
Here $\ol D_i$ is the Zariski closure of $D_i$ in $\CX_{O_L}$, $\CV_i$ is the (uniquely determined) vertical $\QQ$-divisor over $\CX_{O_L}$, and $g_i=g_{D_i}$ is a Green's function over 
$X_L(\CC)$ determined by the 2-variable function $g$ above. 
As in \cite[\S7.1.6]{YZZ}, it will be convenient to denote
$$
\pair{D_1, D_2}:=-\frac{1}{[L:F]}\wh D_1\cdot \wh D_2. 
$$
The definition is independent of the choice of $L$.
We will have a decomposition $\pair{\cdot, \cdot}=-i-j$ in the following. 

We first have equalities
$$
\wh D_1\cdot \wh D_2
=(\ol D_1, g_1)\cdot (\ol D_2+\CV_2, g_2)
=(\ol D_1, g_1)\cdot \ol D_2+\ol D_1\cdot \CV_2+ \int_{X_L(\CC)} g_2\, c_1(D_1,g_1).
$$
Here the first equality holds by $\CV_1\cdot \wh D_2=0$, a consequence of 
condition (1) for $\wh D_2$ and condition (2) for $\CV_1$, and the intersection number 
$(\ol D_1, g_1)\cdot \ol D_2$ is explained above.

So we can write 
$$
\pair{D_1, D_2}=-i(D_1,D_2)-j(D_1,D_2)
$$
with 
$$
i(D_1,D_2):=\frac{1}{[L:F]} (\ol D_1, g_1)\cdot \ol D_2
$$
and 
$$
j(D_1,D_2):=\frac{1}{[L:F]} \ol D_1\cdot \CV_2+\frac{1}{[L:F]} \int_{X_L(\CC)} g_2\, c_1(D_1,g_1).
$$
We  further have a decomposition according to places $v$ of $F$ by
$$j(D_1,D_2)=\sum_{v} j_v(D_1,D_2)\log N_v$$
with
$$
j_v(D_1,D_2):=
\begin{cases}
\ds\frac{1}{[L:F]}  (\ol D_1\cdot \CV_2)_v, & v\nmid\infty,\\
\ds \frac{1}{[L:F]} \int_{X_v(\CC)} g_2\, c_1(D_1,g_1), & v\mid\infty.
\end{cases}
$$
Here we take the convention $\log N_v=1$ for archimedean $v$. 
The local intersection numbers make sense by viewing $\CX_{O_L}$ as a scheme over $O_F$. 

The pairing $j_v(D_1,D_2)$ is symmetric in $D_1,D_2$ for any place $v$. For example, if $v$ is non-archimedean, this is because
$$
(\ol D_1\cdot \CV_2)_v-(\ol D_2\cdot \CV_1)_v
=((\ol D_1+\CV_1)\cdot \CV_2)_v+((\ol D_2+\CV_2)\cdot \CV_1)_v
=(\deg_{F'}(D_1)\hat\xi\cdot \CV_2)_v+(\deg_{F'}(D_2)\hat\xi\cdot \CV_1)_v
=0.
$$

Note that $j_v(D_1,D_2)$ for archimedean $v$ does not necessarily vanish if $\wh D_2$ is not $\hat\xi$-admissible but only weakly $\hat\xi$-admissible.
This is different from \cite[\S7.1.7]{YZZ}, where it considers the $\hat\xi$-admissible case, and thus $j_v(D_1,D_2)=0$ for archimedean $v$.

If $D_1, D_2$ have disjoint supports over $X_{\ol F}$, we can also decompose
$$i(D_1,D_2)=\sum_{v} i_v(D_1,D_2)\log N_v$$
with 
$$
i_v(D_1,D_2):=
\begin{cases}
\ds\frac{1}{[L:F]} (\ol D_1\cdot \ol D_2)_v, & v\nmid\infty,\\
\ds\frac{1}{[L:F]}  g_1(D_{2,v}(\CC)), & v\mid\infty.
\end{cases}
$$

Each of the pairings $i,j,i_v,j_v$ is  symmetric as long as it is defined. In fact, for non-archimedean $v$, this is automatic for $i_v$, and this holds for $j_v$ by $(\ol D_1\cdot \CV_2)_v=-(\CV_1\cdot \CV_2)_v$. 
For archimedean $v$, $g_1(D_{2,v}(\CC))=g_2(D_{1,v}(\CC))$ as they come from the same symmetric 2-variable function $g$, so $i_v$ is symmetric, which implies the symmetry of $j_v$ by Stokes' formula.

As in \cite[\S7.1.7]{YZZ}, we can also introduce the pairings $i_{\bar v}$ and $j_{\bar v}$, and write $i_{ v}$ and $j_{ v}$ respectively as averages of $i_{\bar v}$ and $j_{\bar v}$ over the Galois group $\Gal(\ol F/F)$.

The mistake in \cite{YZZ,YZ} is that the arithmetic extensions used to compute the height pairing are not 
 $\hat\xi$-admissible, but only weakly $\hat\xi$-admissible extension. This will incur  $j_v(D_1,D_2)$ for archimedean $v$.
In the following, we first review the Green's function, compute this extra term, and then decompose the height series by taking into account of the integration term.
We will see that the extra term does not affect the main results of  \cite{YZZ,YZ}, but do affect the main result of this paper.

\subsubsection*{Local nature}

The pairings $i_v, j_v$ for a non-archimedean $v$ has a local nature in that it can be defined in terms of the integral model $\CX_{O_{F_v}}$ instead of $\CX$. Here we explore the situation slightly. 

Let $R$ be a discrete valuation ring with function field $K$. 
Let $Y$ be a projective and smooth curve over $K$. 
For simplicity, assume that $Y$ is geometrically connected, and we omit the extension of the notion to the more general setting. 
Let $\CY$ be a projective, flat, normal and $\QQ$-factorial integral model of $Y$ over $R$.
By a divisor (or $\QQ$-divisor), we mean a Weil divisor (or Weil $\QQ$-divisor).  
 
Fix a $\QQ$-divisor
class $\xi \in {\Pic}(\CY)_\QQ$ whose generic fiber has degree 1 on $Y$.
A $\QQ$-divisor $\CD=\CH+\CV$ on $\CY$, with a horizontal part $\CH$ and a vertical part $\CV$, 
 is called \emph{$\xi$-admissible} if the following
conditions hold:
\begin{itemize}
\item[(1)] The difference $\CD-\deg(\CD_{K})\cdot \xi$ is \emph{flat} in that 
its intersection with any irreducible component of the special fiber of $\CY$ is 0;
\item[(2)] The intersection number $\CV\cdot \xi=0$.
\end{itemize}
Then we have the notion of a \emph{$\xi$-admissible extension} of a $\QQ$-divisor over $Y$.

Let $D_1, D_2$ be two $\QQ$-divisors over $Y$. 
For $i=1,2$, let 
$\wh D_i=\ol D_i+\CV_i$ be a $\xi$-admissible extension of $D_i$ over $\CY$.
Here $\ol D_i$ is the Zariski closure of $D_i$ in $\CY$, and $\CV_i$ is the (uniquely determined) vertical $\QQ$-divisor over $\CY$. 
Assume that the supports of $D_1$ and $D_2$ on $Y$ are disjoint. 
Then we have 
$$
\wh D_1\cdot \wh D_2
=\ol D_1\cdot \ol D_2+\ol D_1\cdot \CV_2.
$$
This gives a decomposition 
$$
\pair{D_1, D_2}=-i_R(D_1,D_2)-j_R(D_1,D_2)
$$
with 
$$
\pair{D_1, D_2}_R=-\wh D_1\cdot \wh D_2,\quad
i_R(D_1,D_2)=\ol D_1\cdot \ol D_2, \quad
j_R(D_1,D_2)= \ol D_1\cdot \CV_2.
$$

Note that $i_R(D_1,D_2)$ is only defined if the supports of $D_1$ and $D_2$ on $Y$ are disjoint, but $j_R(D_1,D_2)$ can actually be defined by the same formula without this condition.

It is worth noting that if the special fiber of $\CY$ is irreducible, then $\CV_1=\CV_2=0$ and $j_R(D_1,D_2)=0$ identically. 

Finally, we mention a projection formula for $i_R$ and $j_R$. 
With the pair $(Y, \CY, \xi)$ over $R$, let $(Y',\CY',\xi')$ be another such pair. 
Assume that there is finite morphism $\pi:Y'\to Y$ 
extending to an $R$-morphism $\tilde\pi:\CY'\to \CY$ such that $\tilde\pi^*\xi=\xi'$. 
Let $D$ and $D'$ be $\QQ$-divisors over $Y$ and $Y'$ respectively. 
Then we have 
$$
j_R(D',\pi^* D)=j_R(\pi_*D', D).
$$
If the supports of $\pi_*D'$ and $D$ on $Y$ are disjoint, we also have
$$
i_R(D',\pi^* D)=i_R(\pi_*D', D).
$$

\subsection{Integral of the Green's function}

Return to the situation that $F$ is totally real, and $X_U$ is a Shimura curve over $F$. 
Fix an archimedean place $v$ of $F$.
The Green's function $g$ over $X_{U,v}(\BC)$ is defined in \cite[\S8.1.1]{YZZ} following the original idea of Gross--Zagier \cite{GZ}.
As noted above, the Green's function is not $\hat\xi$-admissible but only weakly $\hat\xi$-admissible. 
The goal here is to compute the integral of the Green's function, which measures the failure of admissibility.  
The results of this subsection work under the basic conditions in \S\ref{choices}.

Let us first recall the Green's function briefly. 
For any two points $z_1, z_2\in \CH$, the hyperbolic cosine of the
hyperbolic distance between them is given by
$$d(z_1,z_2)=1+\frac{|z_1-z_2|^2}{2\Im(z_1)\Im(z_2)}.$$
It is invariant under the action of $\gl(\RR)$. For any $s\in \BC$
with $\Re(s)>0$, denote
$$m_s(z_1,z_2)=Q_s(d(z_1,z_2)),$$
where
$$Q_s(t)=\int_0^{\infty} \left(t+\sqrt{t^2-1}\cosh u\right)^{-1-s}du $$
is the Legendre function of the second kind. 

Denote by $B=B(v)$ the nearby quaternion algebra.
For any two distinct points of
$$X_{U,v}(\BC)=B_+\cross\bs \CH\times B^\times(\BA_f) /U$$
represented by $(z_1, \beta_1), (z_2, \beta_2)\in \CH\times
B^\times_{\af} $, we denote
$$
g_s((z_1, \beta_1), (z_2, \beta_2)): =\sum_{\gamma\in \mu_U
\backslash B_+\cross} m_s(z_1,\gamma z_2)\
 1_{U}(\beta_1^{-1}\gamma\beta_2).
$$
It converges for $\Re(s)>0$ and has meromorphic continuation to $s=0$ with a simple pole.

The Green's function $g: X_{U,v}(\BC)^2\setminus \Delta\to \RR$
is defined by
$$g((z_1, \beta_1), (z_2, \beta_2)):=\quasilim \ g_s((z_1,
\beta_1), (z_2, \beta_2)).$$
Here $\quasilim$ denotes the constant term at $s=0$ of the Laurent expansion of $g_s((z_1,
\beta_1), (z_2, \beta_2))$.
In particular, for a fixed point 
$P=(z_1, \beta_1)\in X_{U,v}(\BC)$, we can view $g(P,\cdot)$ as a function over $X_{U,v}(\BC)$ with logarithmic singularity at $P$. 

The first part of the following result is a classical one in the computation of Selberg's trace formula, which is a special case of \cite[Proposition 6.3.1(3)]{OT}. 
The second part of the following result is essentially a special case of \cite[Proposition 3.1.2]{OT}, and our proof is a variant of that of the loc. cit.. 

\begin{lem}\label{int of Green}
Let $v$ be an archimedean place of $F$ and $P\in X_{U,v}(\BC)$ be a point.
\begin{enumerate}[(1)]
\item
The residue $\Res_{s=0}g_s\, (P,Q)$ is nonzero only if $Q$ lies in the same connected component as $P$.
In that case,
$$
\Res_{s=0}g_s\, (P,Q) =\frac{1}{\kappa_U^\circ},
$$
where $\kappa_U^\circ$ denotes the degree of $L_U$ on a connected component of $X_{U,v}(\BC)$. 
\item
The integral
$$
\int_{X_{U,v}(\BC)} g(P,\cdot) c_1(\ol\CL_U)=-1.
$$
\end{enumerate}

\end{lem}

\begin{proof}
We will prove that for $\Re(s)>0$, 
$$
\int_{X_{U,v}(\BC)} g_s(P,\cdot) c_1(\ol\CL_U)=
\frac{1}{s(s+1)}.
$$
This implies (2) by taking the constant term. 
It also implies (1).
In fact, the differential equation of the Legendre function transfers to a functional equation
$$
\Delta g_s(P,\cdot)=s(s+1) g_s(P,\cdot).
$$
This implies $\Delta (\Res_{s=0} g_s(P,\cdot))=0$, since $g_s(P,\cdot)$ has at most a simple pole at $s=0$.
It follows that $\Res_{s=0} g_s(P,\cdot)$ is constant on the connected component of $P$. 
Then the integration of $g_s(P,\cdot)$ determines the constant. 

Now prove the formula for the integration of $g_s(P,\cdot)$.
Denote $P=(z_1, \beta_1)$ as above.
As in \cite[\S8.1.1]{YZZ}, the function $g_s(P,\cdot)$ is nonzero only over the connected component of $X_{U,v}(\BC)$ containing $P$.
This connected component is isomorphic to 
$\Gamma\bs\CH$ with $\Gamma= B_+\cross\cap \beta_1 U \beta_1^{-1}$, and the embedding $\Gamma\bs\CH\to X_{U,v}(\BC)$ is given by $z\mapsto (z,\beta_1)$.
Then the induced function $g_s(P,\cdot)$ on $\Gamma\bs\CH$ is given by 
$$
g_s(P,z)=g_s((z_1, \beta_1),(z,\beta_1))
=\sum_{\gamma\in \mu_U
\backslash B_+\cross} m_s(z_1,\gamma z)\
 1_{U}(\beta_1^{-1}\gamma\beta_1)
=\sum_{\gamma\in \mu_U
\backslash \Gamma} m_s(z_1,\gamma z).
$$
It follows that 
$$
\int_{X_{U,v}(\BC)} g_s(P,\cdot) c_1(\ol\CL_U)=
\int_{\Gamma\bs\CH} \sum_{\gamma\in \mu_U
\backslash \Gamma} m_s(z_1,\gamma z) c_1(\ol\CL_U).
$$

Note that the stabilizer of $\CH$ in $\Gamma$ is exactly $\Gamma\cap
F\cross=\mu_U$. Moreover, as in the proof of \cite[Lemma 3.1]{YZZ}, 
$c_1(\ol\CL_U)$ is represented by the standard volume form $\ds\frac{dx\wedge dy}{2\pi y^2}$ over $\CH$. 
Therefore, the integral is further equal to 
$$
\int_{\CH} m_s(z_1, z) \frac{dx dy}{2\pi y^2},
$$
where $z=x+iy$ is as in the convention.

Note that $m_s(\gamma z_1, \gamma z)=m_s(z_1, z)$ for any $\gamma\in \SL_2(\RR)$, the above integral is independent of $z_1$. 
It follows that we can assume $z_1=i$. 
This gives 
$$m_s(i,z)=Q_s(d(i,z))$$
with
$$d(i,z)=1+\frac{|i-z|^2}{2\Im(i)\Im(z)}
=\frac{x^2+y^2+1}{2y}
.$$
Then the integral becomes 
$$
\int_{\CH} Q_s\big(\frac{x^2+y^2+1}{2y}\big) \frac{dx dy}{2\pi y^2}.
$$
We need to prove that this integral is equal to $\ds\frac{1}{s(s+1)}$.
The remaining part is purely analysis.

Denote by $\mathbb D=\{z'\in \CC:|z'|<1\}$ the standard open unit disc.
Under the standard isomorphism $\CH\to \BD$ given by $\ds z'=\frac{z-i}{z+i}$ and $\ds z=i\frac{1+z'}{1-z'}$, the integral becomes
$$
\int_{\BD} Q_s\big(\frac{1+|z'|^2}{1-|z'|^2}\big) \frac{4dx' dy'}{2\pi (1-|z'|^2)^2}.
$$
Here $z'=x'+iy'$ as usual. 
In terms of the polar coordinate $z'=re^{i\theta}$, the integral becomes
$$
\int_0^1 Q_s\big(\frac{1+r^2}{1-r^2}\big) \frac{4rdr}{(1-r^2)^2}
=\int_1^\infty Q_s(t) dt.
$$

Recall from \cite[\S II.2]{GZ} that the Legendre function $Q_s(t)$ satisfies the differential equation
$$
\left((1-t^2)\frac{d^2}{dt^2} -2t\frac{d}{dt} + s(s+1) \right)Q_s=0.
$$
This gives
$$
s(s+1) Q_s= \frac{d}{dt} \left((t^2-1)\frac{d}{dt} Q_s \right).
$$
As a consequence, the original integral is equal to 
$$
\frac{1}{s(s+1)} \left.\left( (t^2-1)\frac{d}{dt} Q_s \right)\right|_1^\infty.
$$

By \cite[II, (2.6)]{GZ}, we can express $Q_s$ by
$$
Q_s(t)=\frac{2^s\Gamma(s+1)^2}{\Gamma(2s+2)}\big(\frac{1}{t+1}\big)^{s+1} F\big(s+1,s+1;2s+2;\frac{2}{t+1}\big).
$$
Here the hypergeometric function 
$$
F(a,b;c;t)=\sum_{n=0}^\infty \frac{(a)_n(b)_n}{(c)_n} \frac{t^n}{n!},
$$
where $(a)_n=a(a+1)\cdots (a+n-1).$
If $c$ is not a negative integer, the series $F(a,b;c;t)$ is absolutely convergent for $|t|<1$, and satisfies the functional equation 
$$
\frac{d}{dt} F(a,b;c;t)=\frac{ab}{c}F(a+1,b+1;c+1;t).
$$
This gives
\begin{align*}
(t^2-1)\frac{d}{dt} Q_s(t)
=& -(s+1)\frac{2^s\Gamma(s+1)^2}{\Gamma(2s+2)} \frac{t-1}{(t+1)^{s+1}} F\big(s+1,s+1;2s+2;\frac{2}{t+1}\big)\\
& -\frac{2^{s+1}\Gamma(s+2)^2}{\Gamma(2s+3)} \frac{t-1}{(t+1)^{s+2}} F\big(s+2,s+2;2s+3;\frac{2}{t+1}\big).
\end{align*}
It follows that for $\Re(s)>0$,  the function $\ds (t^2-1)\frac{d}{dt} Q_s(t)$ converges to 0 as $t\to \infty$. 

By \cite[Theorem 2.1.3]{AAR},  as $t\to 1^+$, 
$$
\big(1-\frac{2}{t+1}\big) F\big(s+1,s+1;2s+2;\frac{2}{t+1}\big)\lra 0$$
and
$$
\big(1-\frac{2}{t+1}\big) F\big(s+2,s+2;2s+3;\frac{2}{t+1}\big)\lra \frac{\Gamma(2s+3)}{\Gamma(s+2)^2}.
$$
Therefore,
$$
\lim_{t\to 1^+} (t^2-1)\frac{d}{dt} Q_s(t)=-\frac{2^{s+1}\Gamma(s+2)^2}{\Gamma(2s+3)} \frac{1}{2^{s+1}} \frac{\Gamma(2s+3)}{\Gamma(s+2)^2}
=-1.
$$ 
This finishes the proof.
\end{proof}

\subsection{Decomposition of the height series}

Assume all the conditions in \S\ref{choices}. 
Note that we particularly need the condition that $U$ contains $\wh O_E^\times$. 

The goal of this subsection to decompose the height series 
\begin{eqnarray*}
Z(g, (t_1, t_2), \phi)_U
= \pair{Z(g, \phi)_U\ t_1^{\circ},\ t_2^{\circ} }_{\NT},
\quad t_1,t_2\in E^\times(\af).
\end{eqnarray*}
The series only depends on the classes of $t_1, t_2\in C_U$ as introduced in \S\ref{sec shimura curve}. 
The height series was treated in \cite[\S 7.1-7.2]{YZZ} in terms of the arithmetic Hodge index theorem and admissible extensions.
But as mentioned above, there is a minor mistake caused by the fact that the Green's function is only weakly admissible, so we will present the correct result here. 
We will still follow the idea of \cite{YZZ, YZ}, but we will also take into account the extra term caused by weak admissibility.

By \ref{sec shimura curve}, we have a canonical arithmetic model $(\CX_U,\ol\CL_U)$ of $(X_U,L_U)$ over $O_F$. Note that $\CX_U$ is a projective and flat normal integral scheme over $O_F$. 
Moreover, the base change $\CX_{U, O_{H}}$ is $\QQ$-factorial for any finite extension $H$ of $F$, so intersection theory is still well-defined for Weil divisors over $\CX_{U, O_{H}}$. 

Recall that $\kappa_U^\circ$ is the degree of $L_U$ over any connected component of $X_{U,\ol F}$. 
Denote 
$$\xi=(\kappa_U^\circ)^{-1} L_U\in \Pic(X_U)_\QQ,\quad
\hat\xi=(\kappa_U^\circ)^{-1} \ol\CL_U\in \wh\Pic(\CX_U)_\QQ.$$ 
For any finite extension $M$ of $F$ unramified over $\Sigma_f$, we can pull $\hat\xi$  back to the base change $\CX_{U, O_{M}}$. Still denote the pull-back by $\hat\xi$ by abuse of notations. 
Then we have the notion of weakly $\hat\xi$-admissible extensions of divisors over $\CX_{U, O_{M}}$.

In particular, for any CM point $[\beta]_U$ represented by $\beta\in \bfcross$, it is defined over the abelian extension $H(\beta)$ of $E$ determined by the open compact subgroup $\beta U \beta^{-1}\cap E^\times(\af)$ of $E^\times(\af)$ via the class field theory. 
By assumption, $U_v$ is maximal at any $v\in \Sigma_f$, so 
$\beta_v U_v \beta_v^{-1}\cap E_v^\times= U_v \cap E_v^\times=O_{E_v}\cross$. 
It follows that the extension $H(\beta)$ of $E$ is unramified over $\Sigma_f$. 
By this, we obtain a weakly $\hat\xi$-admissible extension
$$
\hat\beta=(\ol P_\beta+V_\beta, g(P_\beta,\cdot))
$$ 
of $P_\beta$ over
$\CX_{U, O_{H(\beta)}}$. 
Here $P_\beta$ is the point of $X_U(H(\beta))$ corresponding to $[\beta]$, 
$\ol P_\beta$ is the Zariski closure in $\CX_{U, O_{H(\beta)}}$,
 $V_\beta$ is a vertical divisor over $\CX_{U, O_{H(\beta)}}$,
and $g(P_\beta,\cdot)$ is the Green's function reviewed above. 

Note that the weakly $\hat\xi$-admissible extension $\hat\beta$ is unique, as the Green's function is already chosen. Moreover, the base change of $\hat\beta$ by any extension $M$ of $H(\beta)$ unramified over $\Sigma_f$ is still weakly a $\hat\xi$-admissible extension, which we still denote by $\hat\beta$ by abuse of notations. 

Finally, consider
$$
Z(g, (t_1, t_2), \phi)_U
= \pair{Z_*(g, \phi)_U(t_1-\xi_{t_1}),\ t_2-\xi_{t_2} }_{\NT},
\quad t_1,t_2\in E^\times(\af).
$$
Then the arithmetic Hodge index theorem of Faltings and Hriljac (cf. \cite[Theorem 7.4]{YZZ}) gives
$$
Z(g, (t_1, t_2), \phi)_U
=-( (Z_*(g, \phi)_Ut_1)^\wedge- (Z_*(g, \phi)_U\xi_{t_1})^\wedge )\cdot (\hat t_2-\hat \xi_{t_2} ).
$$
We understand that the arithmetic intersection on the right-hand side involves base changes by finite extensions of $F$ unramified over $\Sigma_f$ to realize $t_1$ as a rational point, and the intersection numbers should be normalized by the degrees of the base changes.
The extension $\hat\xi_{t_2}$ of $\xi_{t_2}$ is given by the corresponding connected component of $\hat\xi$ (over suitable base changes of $X_U$).
The extension $(Z_*(g, \phi)_U\xi_{t_1})^\wedge$ of $ {Z_*(g, \phi)_U \xi_{t_1}}$ is defined similarly, as 
$ {Z_*(g, \phi)_U \xi_{t_1}}$ is a linear combination of connected components of $\xi$. 
The weakly $\hat \xi$-admissible extension $\hat t_2$ of $t_2$ is introduced above. 
The weakly $\hat \xi$-admissible extension $(Z_*(g, \phi)_Ut_1)^\wedge$ of $Z_*(g, \phi)_Ut_1$ is defined similarly, as $Z_*(g, \phi)_Ut_1$ is a linear combination of CM points of the form $[\beta]\in \CMU$.

Take the notational convention 
$$
\pair{D,  D'}:=-\wh D\cdot \wh D',
$$
where $D, D'$ are the divisor classes involved above, and $\wh D, \wh D'$ are the arithmetic extensions introduced above. The right-hand side involves a normalizing factor again if a base change is taken. 
Then the decomposition is written as
$$Z(g, (t_1, t_2))_U
=\pair{Z_*(g,\phi)_U t_1,  t_2}-\pair{Z_*(g,\phi)_U t_1, \xi_{t_2}}
-\pair{Z_*(g,\phi)_U \xi_{t_1}, t_2} +\pair{Z_*(g,\phi)_U\xi_{t_1}, \xi_{t_2}}.$$
Now we summarize the result term by term in the following. 

\begin{thm} \label{height series}
For any $t_1,t_2\in C_U$,
$$Z(g, (t_1, t_2))_U
=\pair{Z_*(g,\phi)_U t_1,  t_2}-\pair{Z_*(g,\phi)_U t_1, \xi_{t_2}}
-\pair{Z_*(g,\phi)_U \xi_{t_1}, t_2} +\pair{Z_*(g,\phi)_U\xi_{t_1}, \xi_{t_2}},$$
where the first term on the right-hand side has the expression 
\begin{align*}
\pair{Z_*(g,\phi)_U t_1,  t_2}
=& -\sum_{v\ \nonsplit} (\log N_v) \barint_{C_U}
\CM^{(v)}_{\phi}(g,(tt_1,tt_2)) dt\\
& -\sum_{v\nmid\infty} \CN^{(v)}_{\phi}(g,(t_1,t_2)) \log N_v-\sum_{v\nmid\infty} j_v(Z_*(g,\phi)t_1,t_2)\log N_v\\
& -\frac{i_0(t_2,t_2)}{[E\cross\cap U:\mu_U]} \Omega_\phi(g,(t_1,t_2))\\
&-\frac12 [F:\QQ]  E_*(0,g,r(t_1,t_2)\phi).
\end{align*}
Here $E_*(0,g,\phi)$ is the non-constant part of the Eisenstein series introduced right before Proposition \ref{analytic series extra},
and the first three lines on the right hand side are the same as the formula of
$Z(g, (t_1,t_2),\phi))_U$ in \cite[Theorem 8.6]{YZ}. 
Namely, they are explained in the following. 
\begin{itemize}
\item[(1)]
The modified arithmetic self-intersection number
$$i_0(t_2,t_2)=i(t_2,t_2)- \sum_{v} i_v(t_2,t_2) \log N_v,$$
where the local term 
$$
i_v(t_2,t_2)=\barint_{C_U}  i_{\bar v}(tt_2,tt_2) dt.
$$
Here the term $i_{\bar v}$ is defined in \cite[\S8.2]{YZ} by case-by-case formulas according to the type of the place $v$, which extends the notion of the $i_{\bar v}$-part in \S\ref{sec admissible} to the current setting of non-proper intersection. 
\item[(2)]
The pseudo-theta series
\begin{eqnarray*}
\Omega_\phi(g,(t_1,t_2))
= \sumu \sum_{y \in
E\cross } r(g,(t_1,t_2))\phi(y,u).
\end{eqnarray*}

\item[(3)]
For any place $v$ nonsplit in $E$,
\begin{eqnarray*}
\CM^{(v)}_{\phi}(g,(t_1,t_2))
&=&w_U \sum_{a\in F\cross}  \quasilim \sum_{y \in \mu_U \backslash (B_+\cross-E\cross)}
r(g,(t_1,t_2))\phi(y)_a   m_s(y), \quad v|\infty, \\
\CM^{(v)}_\phi(g,(t_1,t_2))
&=&\sum_{u\in \mu_U^2\bs F\cross}\sum _{y\in B-E}  r(g,(t_1,t_2))\phi^v (y, u)\  
m_{r(g,(t_1,t_2))\phi _v}(y, u),\quad v\nmid \infty.
\end{eqnarray*}
Here $m_s(y)$ is introduced in \cite[\S8.2]{YZ}, 
and $m_{\phi _v}(y, u)$ is introduced in \cite[Proposition 8.3]{YZ} (and originally in \cite[\S8.2, Notation 8.3]{YZZ}). 

\item[(4)]
For any non-archimedean $v$,
\begin{eqnarray*}
\CN^{(v)}_\phi(g,(t_1,t_2))
&=&\sum_{u\in \mu_U^2\bs F\cross}\sum _{y \in E\cross}  r(g,(t_1,t_2))\phi^v (y, u)\ r(t_1,t_2)n_{r(g)\phi _v}(y, u).
\end{eqnarray*}
Here $n_{\phi _v}(y, u)$ is introduced in \cite[Proposition 8.3]{YZ} for $v$ nonsplit in $E$
and in \cite[Proposition 8.5]{YZ} for $v$ split in $E$. 
\end{itemize}
\end{thm}

\begin{proof}

This is computed in \cite[Theorem 8.6]{YZ}, except that we will have an extra term coming from the weak admissibility. 
In fact, we first write 
$$\pair{Z_*(g) t_1,  t_2}=-i({Z_*(g) t_1,  t_2})-j({Z_*(g) t_1,  t_2}).$$
Then we write 
$$
j({Z_*(g) t_1,  t_2})
=\sum_v j_v({Z_*(g) t_1,  t_2}),
$$
where the sum is over all places $v$ of $F$ instead of just non-archimedean places.
The extra terms are $j_v({Z_*(g) t_1,  t_2})$ for archimedean $v$, while the other terms are computed in the proof of \cite[Theorem 8.6]{YZ}.

If $v$ is archimedean, by definition
$$
j_v({Z_*(g) t_1,  t_2})
=\int_{X_{U,v}(\CC)} g(t_2,\cdot) c_1((Z_*(g) t_1)^\wedge).
$$
Note that only the part of $c_1((Z_*(g) t_1)^\wedge)$ supported on the connected components of $t_2$ contributes to the integral. 
Recall the terminology for the connected components of $Z_*(g)$ in \cite[\S4.3.1]{YZZ}. 
Then we only need to consider the component $Z_*(g)_{q(t_1^{-1}t_2)}$ of $Z_*(g)$.
By the weak admissibility, 
$$
c_1((Z_*(g)_{q(t_1^{-1}t_2)} t_1)^\wedge)
=\deg(Z_*(g)_{q(t_1^{-1}t_2)} )\, c_1(\hat \xi_{t_2}). 
$$
By \cite[Proposition 4.2]{YZZ}, 
$$
\deg(Z_*(g)_{q(t_1^{-1}t_2)} )
=-\frac12 \kappa_U^\circ E_*(0,g,r(t_1,t_2)\phi). 
$$
It follows that 
$$
j_v({Z_*(g) t_1,  t_2})
=-\frac12 \kappa_U^\circ E_*(0,g,r(t_1,t_2)\phi)\int_{X_{U,v}(\CC)} g(t_2,\cdot) c_1(\hat \xi_{t_2}).
$$
By Lemma \ref{int of Green},
$$
\int_{X_{U,v}(\BC)} g(t_2,\cdot) c_1(\ol\CL_U)=-1.
$$
Hence,
$$
j_v({Z_*(g) t_1,  t_2})
=\frac12  E_*(0,g,r(t_1,t_2)\phi).
$$
This finishes the proof.
\end{proof}

\begin{remark} \label{erratum}
The extra term $\ds -\frac12 [F:\QQ]  E_*(0,g,r(t_1,t_2)\phi)$ in Theorem \ref{height series} appears due to the weak admissibility. This term is a priori missed in \cite{YZZ,YZ}. 
However, it does not affect the main results of \cite{YZZ,YZ}, since both articles assume \cite[Assumption 5.4]{YZZ}, under which the extra term vanishes. 
\end{remark}

\subsection{Comparison at archimedean places}

Assume the basic conditions in \S\ref{choices}. 
Let $v$ be an archimedean place of $F$. Recall that in Theorem \ref{height series}, $Z(g,(t_1,t_2))_U$ has a $v$-component 
\begin{eqnarray*}
\CM^{(v)}_{\phi}(g,(t_1,t_2))
=w_U \sum_{a\in F\cross}  \quasilim \sum_{y \in \mu_U \backslash (B_+\cross-E\cross)}
r(g,(t_1,t_2))\phi(y)_a   m_s(y).
\end{eqnarray*}

On the other hand, recall that in Theorem \ref{analytic series}, $\pr' I'(0,g,\phi)_U$ has a $v$-component 
$$
\overline \CK^{(v)}_{\phi}(g,(t_1,t_2)) = w_U \sum_{a\in F\cross}
\quasilim
\sum_{y\in \mu_U\backslash (B(v)_+\cross -E\cross)}  r(g,(t_1,t_2)) \phi(y)_a\
k_{v, s}(y), $$
where
$$
k_{v, s}(y) = \frac{\Gamma(s+1)}{2(4\pi)^{s}}
\int_1^{\infty} \frac{1}{t(1-\lambda(y)t)^{s+1}}  dt,
\qquad \lambda(y)=q(y_2)/q(y).
$$

The goal of this subsection is to compute their difference. The final result is as follows.

\begin{prop} \label{archimedean comparison}
For any $t_1,t_2\in C_U$,
$$
\CK^{(v)}_{\phi}(g,(t_1,t_2))-\overline \CM^{(v)}_{\phi}(g,(t_1,t_2))
=\frac12 (\gamma+\log(4\pi)-1) E_*(0,g,r(t_1,t_2)\phi). 
$$
Here $\gamma$ is Euler's constant.
\end{prop}

\begin{proof}

This is computed as in \cite[Proposition 8.1]{YZZ}, but we need some extra work to take care of contribution from the residue of the Green's function, which is missed in the loc. cit..

As in the calculation of Gross--Zagier \cite{GZ}, 
$$\int_1^{\infty} \frac{1}{t(1-\lambda t)^{s+1}}dt=
\frac{2\Gamma(2s+2)}{\Gamma(s+1)\Gamma(s+2)}Q_s(1-2\lambda)+O_s(|\lambda|^{-s-2}).$$
Moreover, the error term $O_s(|\lambda|^{-s-2})$ vanishes at $s=0$. 
This is a combination of the equations in the first line and in the 12th line of \cite[p. 304]{GZ}, by noting that the left-hand sides of those two equations are equal. 
It follows that 
$$
k_{v, s}(y) = 
\frac{\Gamma(2s+2)}{(4\pi)^{s}\Gamma(s+2)}Q_s(1-2\lambda(y))+O_s(|\lambda(y)|^{-s-2}).
$$

Denote 
\begin{eqnarray*}
\CM^{(v)}_{\phi}(s,g,(t_1,t_2))
=w_U \sum_{a\in F\cross}   \sum_{y \in \mu_U \backslash (B_+\cross-E\cross)}
r(g,(t_1,t_2))\phi(y)_a   m_s(y).
\end{eqnarray*}
Then we have 
$$
\overline \CM^{(v)}_{\phi}(g,(t_1,t_2)) 
=\quasilim \overline \CM^{(v)}_{\phi}(s,g,(t_1,t_2))
$$
and
$$
\overline \CK^{(v)}_{\phi}(g,(t_1,t_2)) 
= \quasilim \frac{\Gamma(2s+2)}{(4\pi)^{s}\Gamma(s+2)}
 \overline \CM^{(v)}_{\phi}(s,g,(t_1,t_2)).
$$
Then
$$
\overline \CK^{(v)}_{\phi}(g,(t_1,t_2)) 
-\overline \CM^{(v)}_{\phi}(g,(t_1,t_2)) 
= \quasilim \left(\frac{\Gamma(2s+2)}{(4\pi)^{s}\Gamma(s+2)}-1\right)
 \overline \CM^{(v)}_{\phi}(s,g,(t_1,t_2)).
$$
Note that $\ds\frac{\Gamma(2s+2)}{(4\pi)^{s}\Gamma(s+2)}-1$ vanishes at $s=0$, and its derivative at $s=0$ is given by
$$
\Gamma'(2)-\log(4\pi)=1-\gamma-\log(4\pi).
$$
We will see that the series $\overline \CM^{(v)}_{\phi}(s,g,(t_1,t_2))$ has a simple pole at $s=0$, coming from the pole of $g_s$ as in Lemma \ref{int of Green}. 
Hence,
$$
\overline \CK^{(v)}_{\phi}(g,(t_1,t_2)) 
-\overline \CM^{(v)}_{\phi}(g,(t_1,t_2)) 
= (1-\gamma-\log(4\pi))\, \Res_{s=0}
 \overline \CM^{(v)}_{\phi}(s,g,(t_1,t_2)).
$$

We can see the simple pole of $\overline \CM^{(v)}_{\phi}(s,g,(t_1,t_2))$ and compute the residue as follows. 
By a simple transformation as in \cite[Proposition 8.1]{YZZ},
we have
\begin{eqnarray*}
\CM^{(v)}_{\phi}(s,g,(t_1,t_2))
=\sum_{a\in F^\times} \sum_{x\in \bb_f^\times/U} r(g)\phi(x)_a\ g_s(t_1x,t_2)
=g_s(Z_*(g,\phi)t_1, t_2).
\end{eqnarray*}
By Lemma \ref{int of Green}, we have
\begin{eqnarray*}
\Res_{s=0} \CM^{(v)}_{\phi}(s,g,(t_1,t_2))
=\frac{1}{\kappa_U^\circ} \deg(Z_*(g,\phi)_{q(t_1^{-1}t_2)}).
\end{eqnarray*}
Here $Z_*(g,\phi)_{q(t_1^{-1}t_2)} t_1$ is the part of $Z_*(g,\phi) t_1$ that lies in the same connected component as $t_2$.
See \cite[\S4.3.1]{YZZ} for the connected components of $Z_*(g,\phi)$. 
In particular, \cite[Proposition 4.2]{YZZ} gives
$$
\deg(Z_*(g,\phi)_{q(t_1^{-1}t_2)} )
=-\frac12 \kappa_U^\circ E_*(0,g,r(t_1,t_2)\phi). 
$$
This finishes the proof.
\end{proof}

\begin{remark} \label{erratum2}
Note that \cite[Proposition 8.1]{YZZ} asserts 
$\CM^{(v)}_{\phi}(g,(t_1,t_2))=\overline \CK^{(v)}_{\phi}(g,(t_1,t_2))$, which is a priori wrong by Proposition \ref{archimedean comparison}.
However, it holds under \cite[Assumption 5.4]{YZZ}, so the correction does not affect the main results of \cite{YZZ,YZ}. 
The situation is similar to Remark \ref{erratum}. 
\end{remark}

\subsection{The $j$-part from bad reduction}

Assume only the basic conditions in \S\ref{choices}. 

If $v$ is a non-archimedean place of $F$ split in $\BB$, then the $j$-part $j_{v}(Z_*(g,\phi)t_1,t_2)=0$ automatically. 
This is a trivial consequence of the fact that every connected component of the fiber $\CX_U$ above $v$ is integral, since it is the quotient of an integral model smooth above $v$. 

Let $v$ be a non-archimedean place nonsplit in $\BB$ (and thus inert in $E$).  
The $j$-part $j_{v}(Z_*(g,\phi)t_1,t_2)$ is treated briefly in \cite{YZZ}
and \cite[Lemma 8.9]{YZ}. 
For the purpose here, we need some precise information. 

Note that $U_v$ is maximal by the basic conditions. 
We further assume that $\phi_v=1_{O_\bv^\times\times\ofv\cross}$, which is part of the restrictive conditions in \S\ref{choices}, but we do not assume that $U$ is maximal. 
This gives us flexibility to vary $U$, which is essential in our proof of the following result. 

\begin{prop} \label{averaged j}
Let $v$ be a non-archimedean place nonsplit in $\BB$ and inert in $E$.  
Then the $j$-part $j_{\bar v}(Z_*(g,\phi)_Ut_1,t_2)$ is a non-singular pseudo-theta series of the form 
$$
\sum_{u\in \mu_U^2\backslash F\cross} \sum_{y\in B(v)-\{0\}}
  r(g,(t_1,t_2)) \phi^v(y,u)\ r(t_1,t_2)l_{r(g)\phi_v}(y,u).
$$
Furthermore, 
$$
\int_{B(v)_v} l_{\phi_v}(y,u)dy=\frac14
|d_v|^2N_v^{-1} (1-N_v^{-1})^2
\cdot 1_{O_{F_v}^\times}(u),
$$
and thus
$$
r(w)l_{\phi_v}(0,u)=- \frac{N_v-1}{4(N_v+1)} r(w)\phi_v(0,u).
$$
\end{prop}

Note that the first part of the proposition is exactly \cite[Lemma 8.9]{YZ}. 
The goal of this subsection is to prove the proposition. 
The idea is still to compute the intersection numbers in terms of the $p$-adic uniformization, but to make a clear picture of the quotient process, we will pass to sufficiently small level structure.

For simplicity, we denote by $B=B(v)$ the nearby quaternion algebra in the following.
Recall that \S\ref{choices} chooses an isomorphism $B(v)_{\AA^v}\simeq \BB^v$ and an embedding $E\to B(v)$, which are compatible with the embedding $E_\AA\to \BB$.  
Here we further fix an isomorphism $B(v)_v\simeq M_2(F_v)$.


\subsubsection*{Pass to smaller level}

Let $U'$ be an open compact subgroup of $U$ with $U'_w=U_w=O_{\BB_w}^\times$ for all non-archimedean places $w$ nonsplit in $\BB$. 
We first check that Proposition \ref{averaged j} for $U'$ implies that for $U$. 
This essentially follows from the local nature and the projection formula of the $j$-part as explained in \S\ref{sec admissible}.

Denote by $\pi:X_{U'}\to X_{U}$ the natural map, and denote by 
$\tilde\pi:\CX_{U'}\to \CX_{U}$ the induced map between the integral models over $O_F$. 
Note the compatibility $\tilde\pi^*\ol\CL_U=\ol\CL_{U'}.$
By the projection formula explained in \S\ref{sec admissible}, we have 
$$j_{\bar v}\big( \pi^*\big(Z_*(g,\phi)_U [t_1]_U\big),[t_2]_{U'}\big)
=j_{\bar v}(Z_*(g,\phi)_U[t_1]_U,[t_2]_U).$$
Here we have used the easy fact $\pi_*([t_2]_{U'})=[t_2]_U$.

We claim that 
$$
\pi^*\big(Z_*(g,\phi)_U [t_1]_U\big)=[\mu_{U'}^2:\mu_U^2]\cdot Z_*(g,\phi)_{U'} [t_1]_{U'}
$$
as a divisor on $X_{U'}$.
In fact, by \cite[Lemma 3.2]{YZZ}, 
$$
\pi^*\big(Z_*(g,\phi)_U [t_1]_U\big)=\big((\pi\times\pi)^*Z_*(g,\phi)_{U}\big) [t_1]_{U'}.
$$
By \cite[Lemma 3.18]{YZZ}, 
$$
(\pi\times\pi)^*Z_*(g,\phi)_{U}=[\mu_{U}^2:\mu_{U'}^2]^{-1} Z_*(g,\phi)_{U'}.
$$
This gives the claim. 
As a consequence, the projection formula gives
$$ j_{\bar v}\big( Z_*(g,\phi)_{U'} [t_1]_{U'},[t_2]_{U'}\big)
=[\mu_{U}^2:\mu_{U'}^2] \, j_{\bar v}\big( Z_*(g,\phi)_{U} [t_1]_U,[t_2]_U\big).$$

Assume that the proposition holds for $U'$. 
Then we have 
$$j_{\bar v}\big( Z_*(g,\phi)_{U'} [t_1]_{U'},[t_2]_{U'}\big)
=\sum_{u\in \mu_{U'}^2\backslash F\cross} \sum_{y\in B-\{0\}}
  r(g,(t_1,t_2)) \phi^v(y,u)\ r(t_1,t_2)l_{r(g)\phi_v}(y,u)_{U'}.
$$
Here the last subscript of $l_{r(g)\phi_v}(y,u)_{U'}$ indicates its dependence on $U'$. 
We assume the slight extra condition that $l_{r(g)\phi_v}(y,u)_{U'}$ as a function of $u$ is invariant under the action of $\ofv\cross$, which will be seen later by its precise expression when $U'$ is sufficiently small. 
Then the proposition holds for $U$ by setting 
$$l_{r(g)\phi_v}(y,u)_{U}=l_{r(g)\phi_v}(y,u)_{U'}.$$

\subsubsection*{Sufficiently small level}

As in \S\ref{sec shimura curve}, 
for any positive integer $N$, denote by $U(N)=(1+N O_{\BB_f})^\times$ the open compact subgroup of $O_{\BB_f}^\times$. 
We say that $U$ (with $U_w$ maximal for all $w\in\Sigma_f$) is \emph{sufficiently small} if $U\subset U(N)$ for 
some integer $N\geq 3$ coprime to the places in $\Sigma_f$.




Denote by $k_v$ the residue field of $v$, and $\bar k_v$ its algebraic closure. 
Denote by $\CX_{U, \bar k_v}=\CX_U\times_{O_F} \bar k_v$ the geometric special fiber.
The goal here is to prove that if $U$ is {sufficiently small}, then $\CX_{U, \bar k_v}$ is a (reduced) semistable curve over the residue field $\bar k_v$, in which every irreducible component $C$ is isomorphic to $\BP^1_{\bar k_v}$, and intersects exactly $N_v+1$ other irreducible components respectively at $N_v+1$ distinct points of $C(\bar k_v)$. 
Moreover, the set of irreducible components can be written as a disjoint union $S_0\cup S_1$ of two subsets, such that any two distinct components in the same subset do not intersect.

The proof is an application of the $p$-adic uniformization of \v Cerednik--Drinfe'ld (cf. \cite{BC}) over $\QQ$ and that of Boutot--Zink \cite{BZ} over a totally real field. 
Recall that the uniformization gives an isomorphism  
$$\wh\CX_{{U}}\times_{\Spf\, O_{F_v}} \Spf\, O_{F_v^\ur} =B^{\times}\bs 
(\wh\Omega \times_{\Spf\, O_{F_v}} \Spf\, O_{F_v^\ur})\times \BB_f\cross/{U}.$$
Here $\wh\CX_{U}$ denotes the formal completion of $\CX_{U}$ along the special fiber above $v$, $F_v^\ur$ denotes the completion of the maximal unramified extension of $F_v$, and $\wh \Omega$ is Deligne's integral model of Drinfe'ld (rigid-analytic) upper half plane $\Omega$ over $O_{F_v}$. 
The group $B_v\cross\cong \GL
_2(F_v)$ acts on $\wh\Omega$ by the fractional linear transformation, and on 
$\bb_v^\times/U_v\cong \ZZ$ via translation by $v\circ q=v\circ\det$.  

Denote a subgroup
$$
B^\dagger=\{\gamma\in B: v(q(\gamma)) =0\}.
$$
Then the uniformization is equivalent to 
$$\wh\CX_{U,O_{F_v^\ur}} =B^\dagger\bs 
\wh\Omega_{O_{F_v^\ur}}\times (\BB_f^v)\cross/{U^v}
=(\mu_{U}\bs B^\dagger)\bs 
\wh\Omega_{O_{F_v^\ur}}\times (\BB_f^v)\cross/{U^v}.$$
Here we denote $\wh\CX_{U, O_{F_v^\ur}}=\wh\CX_{{U}}\times_{\Spf\, O_{F_v}} \Spf\, O_{F_v^\ur}$
and
$\wh\Omega_{O_{F_v^\ur}}=\wh\Omega \times_{\Spf\, O_{F_v}} \Spf\, O_{F_v^\ur}$.
Note that $\wh\Omega_{O_{F_v^\ur}}\times  (\BB_f^v)\cross/{U^v}$
is a disjoint union of countably many copies of $\wh\Omega_{O_{F_v^\ur}}$.

Recall that the dual graph of the special fiber (or equivalently the underlying scheme) of $\wh\Omega$ is just the Bruhat--Tits tree of $\gl$ over $F_v$, whose vertices are homethety classes of lattices in $F_v^2$. 
Choose an irreducible component $C_0\simeq \BP^1$ of the special fiber of $\wh \Omega$.  Apply the transitive action of $\gl(F_v)$ on the set of all irreducible components.
Note that the stabilizer of $\gl(F_v)$ on $C_0$ is $F_v^\times \gl(\ofv)$.
It follows that the irreducible components of the special fiber of $\wh \Omega$ are indexed by 
$\gl(F_v)/F_v^\times \gl(\ofv).$

Denote by $\wt S_0$ (resp. $\wt S_1$)
the set of irreducible components of the special fiber of 
$\wh\Omega_{O_{F_v^\ur}}\times (\BB_f^v)\cross/{U^v}$ represented by 
$\alpha_v F_v^\times \gl(\ofv) \times \beta^v U^v$ 
such that $2|v(q(\alpha_v))$ (resp. $2\nmid v(q(\alpha_v))$). 
We say that two irreducible components $C_1, C_2$ have the \emph{same parity}, if
they belong to the same $\wt S_i$ for some $i=0,1$; otherwise, we say that they have  \emph{different parities}.  

By definition, any $\gamma\in B^\dagger$ satisfies $v(q(\gamma))=0$, so the action of $B^\dagger$ stabilizes $\wt S_0$ and $\wt S_1$ respectively. 
In other words, the action does not change the parity of an irreducible component.
Denote by $S_0$ and $S_1$ respectively the set of irreducible components of the special fiber of $\CX_{U, O_{F_v}^\ur}$ coming from $\wt S_0$ and $\wt S_1$ respectively via the quotient process.
Then $S_0$ and $S_1$ are still disjoint. 
So we can also talk about parities of irreducible components of the special fiber of $\CX_{U, O_{F_v}^\ur}$. 

By the property of Bruhat--Tits tree, if $C$ and $C'$ intersect, then they correspond to adjacent lattices of $F_v^2$, so $C$ and $C'$ should have different parities.
In other words, components of the same parity do not intersect. 
In fact, for two adjacent lattices $\Lambda, \Lambda'$ of $F_v^2$ (corresponding to $C$ and $C'$), the relation 
$p_v \Lambda\subset \Lambda' \subset \Lambda$
implies that $v(\det(\Lambda))$ and $v(\det(\Lambda'))$ have different parities. 
Here 
$\det(\Lambda)$ denotes the transition matrix of the determinant of an $O_{F_v}$-basis of $\Lambda$ to the standard basis of $F_v^2$.

Return to the uniformization 
$$\wh\CX_{U,O_{F_v^\ur}} 
=(\mu_{U}\bs B^\dagger)\bs 
\wh\Omega_{O_{F_v^\ur}}\times (\BB_f^v)\cross/{U^v}.$$
Then the irreducible components of the special fiber of $\CX_{U,O_{F_v^\ur}}$ are indexed by 
$$
(\mu_{U}\bs B^\dagger)\bs  B_{\af}\cross/(F_v^\times \gl(\ofv) U^v).
$$
We claim that for sufficiently small $U^v$, 
the action of $\mu_{U}\bs B^\dagger$ on $B_{\af}\cross/(F_v^\times \gl(\ofv) U^v)$ is free.
 
To prove the claim, take $b\in B_{\af}\cross$, and assume that $\gamma\in B^\dagger$
stabilizes $b(F_v^\times \gl(\ofv) U^v)$, 
and we need to prove $\gamma\in \mu_{U}$. 
The stabilizing condition gives  
$$
\gamma b(F_v^\times \gl(\ofv) U^v) =b(F_v^\times \gl(\ofv) U^v), \quad 
\gamma \in b(F_v^\times \gl(\ofv) U^v)b^{-1}.
$$
By the definition of $B^\dagger$, we have $v(q(\gamma))=0$, so
$$\gamma \in b(\gl(\ofv) U^v)b^{-1}.$$ 

It suffices to prove $\gamma \in F^\times$. 
The proof is very similar to that in 
\cite[Proposition 4.1]{YZ}, but we reproduce it here for convenience of readers. 
Assume the contrary that $\gamma\notin F^\times$. Then $E'=F(\gamma)=F+F\gamma$ is a quadratic CM extension of $F$ contained in $B$. 
Moreover, $\gamma$ lies in  $b(\gl(\ofv) U^v)b^{-1}\cap E_{\af}'^\times$, which is an open and compact subgroup of $E_{\af}'^\times$.
Note that $\wh O_{E'}^\times$ is the unique maximal open compact subgroup of $E_{\af}'^\times$.  
It follows that $\gamma\in \wh O_{E'}^\times$, and thus $\gamma\in O_{E'}^\times$ is a unit. 

On the other hand, by the assumption $U\subset U(N)$, 
we have $1-\gamma\in N b(M_2(\ofv) O_{\BB^v})b^{-1}$. 
The intersection $b(M_2(\ofv) O_{\BB^v})b^{-1}\cap E_{\af}'$ is a compact subring of $E_{\af}'$, so it is contained in the unique maximal compact subring $\wh O_{E'}$ of $E_{\af}'$. It follows that $1-\gamma$ lies in $N\wh O_{E'}\cap E'=N O_{E'}$. 

With the condition $\gamma\in  O_{E'}^\times$ and 
$\gamma\in 1+ N O_{E'}$, it is easy to have a contradiction.
In fact, the element $\zeta=\gamma/ \bar\gamma$ lies in $O_F^\times$ and has absolute value 1 at all archimedean places, so it must be a root of unity. 
By $\gamma\in 1+ N O_{E'}$, we have $\zeta\in 1+N O_{E'}$ and $(1-\zeta)/N\in O_{E'}$. 
Then $(1-\zeta)/N\in \ZZ[\zeta]$ since $\ZZ[\zeta]$ is integrally closed.  
It follows that
we have $\ZZ[\zeta]/N \ZZ[\zeta] \simeq \ZZ/N\ZZ$. 
Note that $\ZZ[\zeta]$ is a free $\ZZ$-module, so it has rank 1 over $\ZZ$. 
It follows that $\zeta=\pm 1$. 
As  $\zeta\in 1+N O_{E'}$, we have $\zeta=1$ and $\gamma\in F^\times.$

In summary, we have just proved that for sufficiently small $U$, 
the action of $\mu_{U}\bs B^\dagger$ on the set of irreducible components of the special fiber of 
$\wh\Omega_{O_{F_v^\ur}}\times (\BB_f^v)\cross/{U^v}$
is free.
Now we prove that for any closed point $x$ on the special fiber, the action of $\mu_{U}\bs B^\dagger$ is free at $x$. 

If $x$ lies on a single irreducible component $C$, then the stabilizer of $x$ is contained in the stabilizer of $C$ and thus must be trivial. 
If $x$ is a node lying in two distinct irreducible components $C, C'$, assume that $\gamma\in \mu_{U}\bs B^\dagger$ is an element fixing $x$, and we need to prove $\gamma=1$.
Note that the action of $\gamma$ either switches $C$ and $C'$, or stabilizes each of $C$ and $C'$.
As $C$ and $C'$ intersect, they must have different parities. 
Then $\gamma$ cannot switch them because $\gamma$ does not change parities.
As a consequence $\gamma$ stabilizes each of $C$ and $C'$.
Then $\gamma=1$ as proved above. 

Therefore, the action of $\mu_{U}\bs B^\dagger$ on the special fiber of 
$\wh\Omega_{O_{F_v^\ur}}\times (\BB_f^v)\cross/{U^v}$
is free.
Then the quotient map is actually \'etale, and
any irreducible component $C$ is mapped birationally to its image $\overline C$ in the special fiber of $\CX_{U,O_{F_v^\ur}}$.
To prove that $C\to \overline C$ is actually an isomorphism, we need to further check that 
the nodes of $C$ is injective into $\overline C$. 
In other words, we need to prove that there is no $\gamma\in \mu_{U}\bs B^\dagger$ mapping a node $x_1\in C$ to another node $x_2\in C$. 
Assume that $\gamma$ exists. 
Denote by $C_1$ (resp. $C_2$) the irreducible component other than $C$ passing through $x_1$ (resp. $x_2$). 
By the intersection configuration,  $C_1$ and $C_2$ have the same parity, which is different from that of $C$. 
Note $\gamma$ maps the set $\{C, C_1\}$ to the set $\{C, C_2\}$. 
It follows that $\gamma$ maps $C$ to $C$, as it does not change parity. 
Then $\gamma$ is trivial again, and we have a contradiction.

Via the quotient process, we conclude that 
 the special fiber $\CX_{U, \bar k_v}$ of $\CX_{U,O_{F_v^\ur}}$  is a semistable curve over 
 $\bar k_v$, in which every irreducible component $C$ is isomorphic to $\BP^1_{k_v}$, and intersects exactly $N_v+1$ other irreducible components respectively at $N_v+1$ points in $C$. 
Moreover, the set of irreducible components can be written as a disjoint union $S_0\cup S_1$ of two subsets, such that any two distinct components in the same subset do not intersect.

In the following, we assume that $U$ is sufficiently small, so that the quotient process is free, and the special fiber has the above nice properties.

\subsubsection*{Pass to the uniformization}

Note that the first part of the proposition is exactly \cite[Lemma 8.9]{YZ}. 
In the following, we first recall the formula of $l_{\phi_v}(g,y,u)$ in \cite[Lemma 8.9]{YZ}, and then compute its average by a more careful analysis of the $p$-adic uniformization. 

By the basic conditions in \S\ref{choices}, $O_{E_v}^\times\subset O_{\BB_v}\cross$, so the CM point $[t]_U$ for any $t\in E_{\af}^\times$ is defined over $F_v^\ur$. 
Therefore, to compute $j_{\overline v}$, by choosing an embedding 
$F_v^\ur\to (\ol F)_{\bar v}$, it suffices to compute it over the integral model $\CX_{U,O_{F_v^\ur}}$. Note that the $j$-part can be defined locally as in \S\ref{sec admissible}. 

Recall the uniformization 
$$\wh\CX_{U, O_{F_v^\ur}} =B^{\times}\bs 
\wh\Omega_{O_{F_v^\ur}} \times \BB_f\cross/U.$$
By definition,
$$j_{\overline v}(Z_*(g)t_1, t_2)= \overline{Z_*(g)t_1}\cdot
V_{t_2}.$$
Here $\overline{Z_*(g)t_1}$ is the Zariski closure in $\CX_{U,O_{F_v^\ur}}$, and $V_{t_2}$ is the unique vertical divisor on $\CX_{U,O_{F_v^\ur}}$, supported on the geometrically connected component of $t_2$ in $\CX_{U,O_{F_v^\ur}}$,
 satisfying the following properties:
\begin{itemize}
\item[(1)] $(V_{t_2}+\bar t_2)\cdot C=\hat\xi \cdot C$ for any vertical divisor $C$ of $\CX_{U,O_{F_v^\ur}}$;
\item[(2)] $V_{t_2} \cdot \hat\xi =0$.
\end{itemize}

Write $V_{1}=\sum_i a_i W_i$ (for $t_2=1$), where $\{W_i\}_i$ is the set of irreducible components of the special fiber of $\CX_{U,O_{F_v^\ur}}$ lying in the same connected component as $1$. 
Let $\wt W_i$ be an irreducible component of the special fiber of $\wh\Omega _{O_{F_v^\ur}}$ lifting $W_i$.
Write $\wt V=\sum_i a_i \wt W_i$, viewed as a vertical divisor of 
$\wh\Omega_{O_{F_v^\ur}}$. 

Via the $p$-adic uniformization, the proof of \cite[Lemma 8.9]{YZ} actually gives a non-singular pseudo-theta series
$$j_{\overline v}(Z_*(g)t_1, t_2)= \sumu \sum_{\gamma \in B^\times}
 r(g,(t_1,t_2))\phi^v(\gamma, u)\ r(t_1,t_2) l_{\phi_v}(g,\gamma,u),
$$
where for $g\in \gl(F_v), \ \gamma\in B_v^\times,\ u\in F_v^\times$,
$$l_{\phi_v}(g,\gamma,u)=
\sum_{x \in\bb_v^\times /U_v} r(g)\phi_v(x, uq(\gamma)/q(x))
1_{O_{F_v}^\times}(q(x)/q(\gamma))\ (\gamma^{-1} z_0\cdot \wt V).
$$
Here $z_0\in \Omega (E_v)$ is a point in $\Omega (\BC_v)$
fixed by $E_v^\times$, viewed as a point of $\Omega (F_v^\ur)$ via an embedding $E_v\to F_v^\ur$ fixed by us once for all. 
Moreover, we also have
$$l_{\phi_v}(1,\gamma,u)= 
(\gamma^{-1}\bar z_0\cdot \wt V)\cdot 1_{O_{F_v}^\times}(q(\gamma)) \cdot 1_{O_{F_v}^\times}(u).$$
Here $\bar z_0$ denotes the section of $\wh\Omega_{O_{F_v^\ur}}$ corresponding to $z_0$. 
The non-singularity of the pseudo-theta series is implied by the fact that 
$l_{\phi_v}(1,\gamma,u)$ is actually a compactly-supported function of $(\gamma, u)\in B_v^\times \times F_v^\times$.
To make the right-hand nonzero, assume that $u\in O_{F_v}^\times$ in the following. 

Note that the proof of \cite[Lemma 8.9]{YZ} claims that there is a unique point of $\Omega (F_v^\ur)$ fixed by $E_v^\times$. This is wrong. In fact, there are exactly two points of $\Omega (\BC_v)=\CC_v-F_v$ fixed by $E_v^\times$, both of which are defined over $E_v$. 
To see this, by the Noether--Skolem theorem, the embedding 
$E_v\to M_2(F_v)$ is a conjugate of the embedding
$$
a+b\sqrt{D_v} \longmapsto \matrixx{a}{bD_v}{b}{a},
$$
where $D_v\in F_v^\times$ is the relative discriminant of $E_v$ over $F_v$.
This explicit case gives fixed points $\pm\sqrt{D_v}$.
This mistake of the loc. cit. does not affect the other results, since any fixed point gives the same results for our purpose.

\subsubsection*{Compute the average}

Let $\widetilde W^0$ to be the unique irreducible component of the special fiber of $\wh \Omega_{O_{F_v}^\ur}$ 
intersecting the section $\bar z_0$.
Recall that the irreducible components of the special fiber of $\wh \Omega_{O_{F_v}^\ur}$ are indexed by 
$\gl(F_v)/F_v^\times \gl(\ofv).$
We can further adjust the indexing such that $\widetilde W^0$ corresponds to the trivial coset. 

Denote by $\alpha_i F_v^\times \gl(\ofv)$ the coset representing the component $\wt W_i$.
Then we have  
$$
(\gamma^{-1} z_0\cdot \wt W_i)
=1_{\alpha_iF_v^\times \gl(\ofv)}(\gamma^{-1})
=1_{F_v^\times \gl(\ofv) \alpha_i^{-1}}(\gamma).
$$
By $\wt V=\sum_i a_i \wt W_i$, we have
$$\int_{B_v} l_{\phi_v}(1,\gamma,u) d\gamma 
= \sum_i a_i \int_{B_v^\dagger} 
(\gamma^{-1} z_0\cdot \wt W_i) d\gamma
=\sum_i a_i\ \vol(B_v^\dagger\cap F_v^\times \gl(\ofv)\alpha_i^{-1}).$$
Here 
$$
B_v^\dagger=\{\gamma\in B_{v}: q(\gamma) \in O_{F_v}^\times\}.
$$

For any $\alpha\in B_v^\times$, it is easy to have 
$$B_v^\dagger\cap F_v^\times \gl(\ofv) \alpha^{-1}
=\begin{cases}
\emptyset, & 2\nmid v(q(\alpha));  \\
 \gl(\ofv) \varpi_v^{v(q(\alpha))/2} \alpha^{-1},
& 2\mid v(q(\alpha)).
\end{cases}$$
Here $\varpi_v$ is a uniformizer of $F_v$. 
Note that the self-dual measure on $B_v=M_2(F_v)$ gives
$$\vol(\GL_2(O_{F_v}))=|\GL_2(O_{F_v}/p_v)|\cdot\vol(1+p_vM_2(O_{F_v}))=|d_v|^2(1-N_v^{-1})(1-N_v^{-2}).$$
As a consequence 
$$\int_{B_v} l_{\phi_v}(1,\gamma,u) d\gamma 
=|d_v|^2(1-N_v^{-1})(1-N_v^{-2}) \sum_{i:\, 2\mid v(q(\alpha_i))} a_i.$$

Recall that $S_0$ (resp. $S_1$) denotes the set of $i$ such that $2\mid v(q(\alpha_i))$ (resp. $2\nmid v(q(\alpha_i))$).
Then $S=S_0\cup S_1$ is the set of all indexes $i$. 
Denote 
$$
A_0=\sum_{i\in S_0} a_i, \qquad
A_1=\sum_{i\in S_1} a_i. 
$$
We need to compute $A_0$.
We are going to prove the following equation:
$$
A_0+A_1=0, \qquad
A_0-A_1=\frac{1}{2(N_v+1)}. 
$$
The relations give 
$$A_0=\frac{1}{4(N_v+1)}, \qquad
\int_{B_v} l_{\phi_v}(1,\gamma,u) d\gamma 
=\frac14 |d_v|^2N_v^{-1}(1-N_v^{-1})^2.
$$ 
Then the last equality of the proposition follows from Lemma \ref{derivative of intertwining}.

It remains to prove the two equations of $A_0$ and $A_1$. 
We need the following intersection results:
\begin{enumerate}[(1)]
\item The orders $|S_0|$ and $|S_1|$ are equal. In fact, for any $x_v\in \BB_v\cross$ with $2\nmid v(q(x_v))$, the Hecke correspondence $Z(x_v)_U$ corresponding to $U_vx_vU_v=x_vU_v$ is an automorphism of $\CX_{U,O_{F_v^\ur}}$ and switches $S_0$ with $S_1$. Denote $n=|S_0|=|S_1|$ in the following, so $2n=|S|$.

\item $W_i\cdot W_i=-(N_v+1)$ for any $i\in S$. In fact, by the above result (as we have assumed that $U$ is sufficiently small), any irreducible component of the special fiber of $\CX_{U, O_{F_v^\ur}}$ is isomorphic to $\BP^1$, and intersects with exactly $N_v+1$ other components. As a consequence, 
$$
W_i\cdot W_i=-W_i\cdot \sum_{j\in S, \ j\neq i} W_j=-(N_v+1).
$$

\item Fix $r=0,1$. Then $W_i\cdot W_{j}=0$ for $i,j\in S_r$ with $i\neq j$. 
This is just the above fact that components of the same parity do not intersect. 

\item $W_i\cdot \hat\xi=1/(2n)$ for any $i\in S$.
Note that $\hat\xi \cdot\sum_{i\in S} W_i=1$, since $\xi$ has degree one on every  connected component of $X_{U,F_v^\ur}$. 
Then it suffices to prove that $W_i\cdot \hat\xi$ is independent of $i$, or equivalently 
$W_i\cdot \CL_U$ is independent of $i$.  
As $U$ is sufficiently small, the complex uniformization for $X_U$ is also a free quotient. 
This fact is contained in the proof of \cite[Proposition 4.1]{YZ}.
Moreover, $X_U$ has no cusp since we are in the case that $v$ is nonsplit in $\BB$.
As a consequence, the Hodge bundle $\CL_U$ is just the relative dualizing sheaf $\omega_{\CX_U/O_F}$ of the regular scheme $\CX_U$ over $O_F$.  
Apply the adjunction formula
$$
2g(W_i)-2= W_i\cdot \CL+W_i\cdot W_i.
$$
As 
$g(W_i)=0$ and $W_i\cdot W_i=-(N_v+1)$. This gives $W_i\cdot \CL=N_v-1$. 
\end{enumerate}

Now we are ready to establish the equations for $A_0$ and $A_1$. 
By the definition of $V_1=\sum_i a_i W_i$, 
we have $V_1 \cdot \hat\xi =0$. 
This is just $A_0+A_1=0$ by (4). 

On the other hand, the definition of $V_1$ also gives 
$(\bar 1+V_1)\cdot C=\hat\xi\cdot C$ for any vertical divisor $C$ of $\CX_{U,O_{F_v^\ur}}$. 
Take $C=\sum_{j\in S_1} W_j$. 
It does not intersect the Zariski closure of $1$. 
Furthermore, by (2) and (3), $W_i\cdot C=N_v+1$ for $i\in S_0$ and 
$W_i\cdot C=-(N_v+1)$ for $i\in S_1$.
By (1) and (4), we have $\hat\xi\cdot C=1/2$.
Then the identity $(\bar 1+V_1)\cdot C=\hat\xi\cdot C$ becomes
$$
\sum_{i\in S_0} a_i (N_v+1)
-\sum_{i\in S_1} a_i (N_v+1)=\frac12.
$$
This gives our second equation. 
The proof of Proposition \ref{averaged j} is complete.

\subsection{Hecke action on arithmetic Hodge classes} \label{sec 4.3}

In last subsection, we have the decomposition
$$Z(g, (t_1, t_2))
=\pair{Z_*(g,\phi) t_1,  t_2} -\pair{Z_*(g,\phi) \xi_{t_1}, t_2}
+\pair{Z_*(g,\phi)\xi_{t_1}, \xi_{t_2}} -\pair{Z_*(g,\phi) t_1, \xi_{t_2}}.$$
We have also considered a decomposition of the first term on the right-hand side.
In this subsection we consider the remaining three terms. 
The treatment here is an enhanced version of \cite[\S7.3]{YZZ}. 
We still assume all the conditions in \S\ref{choices} in the following, unless otherwise instructed.

\subsubsection*{Two easy terms}

Recall that $\kappa_U^\circ$ is the degree of $L_U$ on a connected component of $X_{U,\ol F}$. The following result works under the basic conditions of \S\ref{choices}. 

\begin{pro} \label{geometric series extra1}
$$\pair{Z_*(g,\phi) \xi_{t_1}, t_2}=-\frac 12 \kappa_U^\circ\ E_*(0,g,r(t_1,t_2) \phi)_U\cdot \pair{\xi_{t_2}, t_2},$$
$$\pair{Z_*(g,\phi)\xi_{t_1}, \xi_{t_2}}=-\frac 12 \kappa_U^\circ\ E_*(0,g,r(t_1,t_2) \phi)_U\cdot \pair{\xi_{t_2}, \xi_{t_2}}.$$
\end{pro}

\begin{proof}

We first compute $\pair{Z_*(g,\phi) \xi_{t_1}, t_2}$. 
By definition, $Z_*(g,\phi) \xi_{t_1}$ is a linear combination of 
$Z(x) \xi_{t_1}$. 
By construction, the correspondence $Z(x)$ keeps the
canonical bundle up to a multiple under pull-back and push-forward.
More precisely, one has
$$Z(x) \xi_{t_1} =(\deg Z(x)) \xi_{t_1 x},\quad\forall x\in\bb_f\cross.$$
Note that $\pair{Z(x) \xi_{t_1}, t_2}$ is nonzero only if $\xi_{t_1 x}$ and $t_2$ lie in the same geometrically connected component of $X_U$. 
It follows that 
$$\pair{Z_*(g,\phi) \xi_{t_1}, t_2}
=\pair{Z_*(g,\phi)_{U,q(t_1^{-1}t_2)} \xi_{t_1}, t_2}
=\deg Z_*(g,\phi)_{U,q(t_1^{-1}t_2)} \cdot \pair{\xi_{t_2}, t_2}.$$ 
Here $Z_*(g,\phi)_{U,q(t_1^{-1}t_2)}$ consists of the $q(t_1^{-1}t_2)$-component of $Z_*(g,\phi)_{U}$ as introduced in \cite[\S 4.2.4]{YZZ}. 
By \cite[Proposition 4.2]{YZZ},
\begin{eqnarray*}
\deg Z(g, \phi)_{U,q(t_1^{-1}t_2)}=-\frac 12 \kappa_U^\circ\ E(0,g,r(t_1,t_2) \phi)_U.
\end{eqnarray*}
This gives the formula for $\pair{Z_*(g,\phi) \xi_{t_1}, t_2}$. 
The same method also proves the formula for $\pair{Z_*(g,\phi) \xi_{t_1}, \xi_{t_2}}$. 
\end{proof}

\subsubsection*{Almost eigenvector}

It remains to consider $\pair{Z_*(g) t_1, \xi_{t_2}}$.
We follow the treatment of \cite[\S7.3.2]{YZZ} with some modification to fit the current setting. 

For any $x\in \bb_f\cross$, let $\CZ(x)$ be the Zariski closure of $Z(x)$ in $\CX_U\times_{O_F} \CX_U$.
Note that $U$ is maximal by assumption. The following are true: 
\begin{itemize}
\item[(1)] $\CZ(x_1)$ commutes with $\CZ(x_2)$ for any $x_1, x_2\in
\bb_f\cross$;
\item[(2)] $\CZ(x)=\prod_{v\nmid\infty} \CZ(x_v)$ for any $x\in \bb_f\cross$;
\item[(3)] for any $x\in \bb_f\cross$, both structure projections from
$\CZ(x)$ to $\CX_U$ are finite.
\end{itemize}

In the proof of Proposition \ref{geometric series extra1}, we already see that 
$$
Z(x) \xi= (\deg Z(x))\ \xi.
$$
In other words, $\xi$ is an eigenvector of $Z(x)$ over $X_U$.
For the arithmetic version, we will see that $\hat\xi$ generally fails to be an eigenvector of 
$\CZ(x)$, but the failure is explicitly computable. 

Define an arithmetic class $D(x)$ on $\CX_U$ by
$$D(x):=\CZ(x) \hat\xi- (\deg Z(x))\ \hat\xi.$$
Then $D(x)$ is a vertical arithmetic $\BQ$-divisor since it is zero on the generic fiber.

If $x\in \bb_v^\times$ for some non-archimedean place $v$ nonsplit in $\BB$, then we have $\deg Z(x)=1$ and $D(x)=0$. 
In fact, since $U_v=O_{\BB_v}^\times$, the double coset $U_vxU_v=xU_v$ is a single coset depending only on $v(q(x))$. 
As a consequence, $Z(x)$ is just an automorphism of $X_U$, and thus $\CZ(x)$ is an automorphism of $\CX_U$.
For any subgroup $U'^v \subset U^v$, we have a similar automorphism on $\CX_{U_vU'^v}$ determined by $x$, and this automorphism does not change the relative dualizing sheaf of $\CX_{U_vU'^v}$. By the compatibility of the arithmetic Hodge bundle with the norm map reviewed in \S\ref{sec shimura curve}, we see that $\CZ(x)$ fixes the arithmetic class $\hat\xi$. 

If $x\in \bb^\Sigma$, then $D(x)$ is a \emph{constant $\BQ$-divisor}, i.e., the pull-back of an arithmetic $\BQ$-divisor from $\Spec (O_{F'})$, where $F'$ is the algebraic closure of $F$ in the functor field of $X_U$.
Note that $F'$ is the abelian extension of $F$ with Galois group
$\pi_0(X_{U,\overline F})=F_+^\times\bs \afcross/q(U)$ via the class field theory.
See the reason for the constancy of $D(x)$ in \cite[\S7.3.2]{YZZ}.

Hence, for all $x\in \bb_f^\times$, $D(x)$ is a constant $\BQ$-divisor, i.e., the pull back of an arithmetic $\BQ$-divisor from 
$\Spec (O_{F'})$. 
By abuse of notation, we also denote by $D(x)$ the 
arithmetic degree of the arithmetic $\BQ$-divisor on 
$\Spec (O_{F'})$. 
Hence we get a number $D(x)\in \BR$. 
It is more convenient to introduce 
$$
D_0(x):= \frac{1}{\deg Z(x)} D(x).
$$
By definition, $D_0(x)$ is \emph{additive} in that
$$
D_0(x)=\sum_{v\nmid\infty} D_0(x_v).
$$
The sum has only finitely many nonzero terms.

Now we have the following basic result. 

\begin{lem} \label{derivation basic}
For any $t\in C_U$, 
\begin{align*}
\pair{Z(x) t,  \xi}
= \deg Z(x)\, \pair{t, \xi} - \deg Z(x)\, \sum_{v\nmid \infty}D_0(x_v) .
\end{align*}
\end{lem}
\begin{proof}
This is a direct consequence of \cite[Lemma 7.7]{YZZ}, which asserts
\begin{align*}
\pair{Z(x) D,  \xi}
= \deg Z(x)\, \pair{D, \xi} -  \deg(D) D(x), \quad D\in \Div(X_{U,\ol F}), \ x\in \bvcross.
\end{align*}
There is a gap in the proof of the loc. cit. due to the extra term caused by the weak admissibility, but the conclusion still holds. 
In fact, the loc. cit. proves that  
\begin{align*}
\pair{\overline{Z(x) D},  \hat\xi}
= \deg Z(x)\, \pair{\ol D, \hat\xi} - \deg(D) D(x).
\end{align*}
On the other hand, by Lemma \ref{int of Green},
$$
\pair{D,  \xi}=\pair{\overline{D},  \hat\xi}+(\kappa_U^\circ)^{-1} \deg(D),
$$
and 
$$
\pair{Z(x) D,  \xi}=\pair{\overline{Z(x) D},  \hat\xi}+(\kappa_U^\circ)^{-1} \deg Z(x) \deg(D).
$$
This implies the original statement. 
\end{proof}

\subsubsection*{The last term}

Now we are ready to compute $\pair{Z_*(g,\phi) t_1, \xi_{t_2}}$.

\begin{pro} \label{geometric series extra2}
\begin{eqnarray*}
\pair{Z_*(g,\phi) t_1, \xi_{t_2}}
= -\frac 12 \kappa_U^\circ\ E_*(0,g,r(t_1,t_2) \phi)_U\ \pair{[1], \xi} 
 +\frac 12  \sum_{v\notin \Sigma} \CF^{(v)}_\phi(g, (t_1,t_2)),
\end{eqnarray*}
where
$$
\CF^{(v)}_\phi(g, (t_1,t_2))
= \sumu\sum_{a\in F^\times}  W_{a}^v(0,g,u,r(t_1,t_2)\phi) \, 
f_{\phi_v, a}(g,(t_1,t_2),u)$$
with
$$
f_{\phi_v, a}(g,(t_1,t_2),u)=(1-N_v^{-2})|d_v|^{\frac32} |au^{-1}|_v\, \kappa_U^\circ
\sum_{y\in \bb_v(au^{-1})/U_v^1} r(g, (t_1,t_2))\phi_v(y,u)  D_0(t_{1,v}^{-1}y t_{2,v}).
$$
Here
$\bb_v(a)=\{x\in \bb_v: q(x)=a \}.$
\end{pro}
\begin{proof}
Denote $t=t_1t_2^{-1}$.
We have 
$$\pair{Z_*(g,\phi) t_1, \xi_{t_2}}
=\pair{Z_*(g,\phi)_{q(1/t)} t_1, \xi_{t_2}}
=\pair{Z_*(g,\phi)_{q(1/t)} t_1, \xi}
=\pair{Z_*(g,\phi)_{q(1/t)} [1], \xi}.$$
Here the first equality holds as in the proof of 
Proposition \ref{geometric series extra1}, the second equality holds by a similar reason of geometrically connected components, and the third equality holds by the Galois action associated to $t_1$.  

Recall that from \cite[\S 4.2.4]{YZZ} we have 
\begin{eqnarray*}
Z_*(g)_{q(1/t)}
&=& w_U \sum_{u\in \mu_U'\bs F^\times} \sum_{a\in F_+\cross}
\sum_{y\in K^t\bs \bb_f(a)}
 r(g, (t,1))\phi(y,u)Z(t^{-1}y).
\end{eqnarray*}
Here
$$K^t=
\GSpin(\BB_f,q)\cap tKt^{-1}
= \{(h_1, h_2)\in (tUt^{-1}) \times U: q(h_1)=q(h_2) \}$$
acts on 
$$\bb_f(a)=\{x\in \bb_f: q(x)=a \}$$
by $(h_1,h_2):x\mapsto h_1xh_2^{-1}$. 

Hence, Lemma \ref{derivation basic} gives 
\begin{eqnarray*}
\pair{Z_*(g,\phi)_{q(1/t)} [1], \xi}
= \deg(Z_*(g,\phi)_{q(1/t)}) \pair{[1], \xi} 
 +\sum_{v\nmid\infty} \CF^{(v)}_\phi(g, (t_1,t_2)),
\end{eqnarray*}
where
$$
\CF^{(v)}_\phi(g, (t_1,t_2))
=-w_U \sum_{u\in \mu_U'\bs F^\times} \sum_{a\in F_+\cross}
\sum_{y\in K^t\bs \bb_f(a)}
 r(g, (t,1))\phi(y,u)\cdot \deg Z(t^{-1}y)\cdot D_0(t_v^{-1}y_v).
$$
As in the proof of Proposition \ref{geometric series extra1}, we already have 
\begin{eqnarray*}
\deg Z(g, \phi)_{U,q(t^{-1})}=-\frac 12 \kappa_U^\circ\ E(0,g,r(t_1,t_2) \phi)_U.
\end{eqnarray*}
It remains to convert the above expression of $\CF^{(v)}_\phi(g, (t_1,t_2))$ to the form in the proposition. 

Consider the last summation
\begin{eqnarray*}
 \sum_{y\in K^t\bs \bb_f(a)}
 r(g, (t,1))\phi(y,u) \deg Z(t^{-1}y)\cdot D_0(t_v^{-1}y_v).
\end{eqnarray*}
By $\deg Z(t^{-1}y)=|Ut^{-1}yU/U|$, the summation is equal to  
\begin{eqnarray*}
 \sum_{y\in K^t\bs \bb_f(a)} \ \sum_{x\in Ut^{-1}yU/U}
 r(g)\phi(x,q(t)u)  D_0(x_v).
\end{eqnarray*}
Note that 
$$Ut^{-1}yU/U=Kt^{-1}y/U=t^{-1}(tKt^{-1}y/U)=t^{-1}(K^ty/U^1).$$
The summation becomes  
\begin{eqnarray*}
 && \sum_{y\in \bb_f(a)/U^1}
 r(g)\phi(t^{-1}y,q(t)u)  D_0(t_v^{-1}y_v) \\
&=& \left( \sum_{y\in \bb_f^v(a)/(U^v)^1}
 r(g, (t,1))\phi^v(y,u)  \right)
 \cdot \left( \sum_{y_v\in \bb_v(a)/U_v^1}
 r(g, (t,1))\phi_v(y,u)  D_0(t_v^{-1}y_v)\right).
\end{eqnarray*}
We assume that $v$ is split in $\BB$; otherwise $D_0(t_v^{-1}y_v)=0$ identically. 
It suffices to convert the first summation on the right-hand side in this case. 
The proof is similar to the proof of \cite[Proposition 4.2]{YZZ}, except that we do not convert the second summation on the right-hand side.

In fact, by \cite[Proposition 2.9]{YZZ}, 
\begin{eqnarray*}
 \sum_{y\in \bb_f^v(a)/(U^v)^1}
 r(g, (t,1))\phi^v(y,u)  
 = -\frac{|a|_v}{\vol((U^v)^1)\vol(\bb_\infty^1)}W_{au}^v(0,g,u,r(h)\phi).
\end{eqnarray*}
The negative sign comes from the Weil index of $\BB^v$, which is $-1$ since $\BB_v$ is a matrix algebra. 

Finally, apply equation (4.3.2) in the proof of 
\cite[Proposition 4.2]{YZZ}. 
Note that $\vol(U_v^1)=(1-N_v^{-2})|d_v|^{\frac32}$ by the normalization in \cite[\S1.6.2]{YZZ}.
It remains to check 
$$
f_{\phi_v, a}(g,(t_1,t_2),u)=f_{\phi_v, a}(g,(t,1),u).
$$
This can be obtained by writing the sum over $\bb_v(au^{-1})/U_v^1$ as an integral over $\bb_v(au^{-1})$. 
\end{proof}

For simplicity, write 
$$
\CF^{(v)}_\phi(g)
= \CF^{(v)}_\phi(g, (1,1)), \quad
f_{\phi_v, a}(g,u)=f_{\phi_v, a}(g,(1,1),u), \quad
f_{\phi_v, a}(1,u)=f_{\phi_v, a}(1,(1,1),u).
$$

\begin{pro} \label{geometric series extra3}
For any non-archimedean place $v$ nonsplit in $\BB$,
$f_{\phi_v, a}(1, u)\neq 0$ only if $a\in \ofv$ and $u\in\ofv\cross$. 
In that case, 
$$
f_{\phi_v, a}(1,u)=  
|d_v|^{\frac32} 
\frac{1+N_v^{-1}}{1-N_v^{-1}}
\big((r+2)N_v^{-(r+1)}-rN_v^{-(r+2)}-(r+2)N_v^{-1}+r\big) \log N_v.
$$
Here $r=v(a)$. 
\end{pro}
\begin{proof}
This is essentially computed in \cite{Zh1}. 
By definition,
$$
f_{\phi_v, a}(1,u)=(1-N_v^{-2})|d_v|^{\frac32} |au^{-1}|_v \, \kappa_U^\circ
\sum_{y\in \bb_v(au^{-1})/U_v^1} \phi_v(y,u)  D_0(y).
$$
It is nonzero only if $u\in \ofv\cross$ and $a\in \ofv$, which we assume in the following. 
Identify $\bb_2=M_2(F_v)$ and $O_{\BB_v}=M_2(\ofv)$. 
Note $r=v(a)\geq 0$. Denote
$$
M_2(\ofv)_r=\{y \in M_2(\ofv)_r: v(\det(y))=r\}.
$$
Then the summation equals 
$$
\sum_{y\in (\bb_v(au^{-1}) \cap O_{\BB_v})/U_v^1}  D_0(y)
=\sum_{y\in M_2(\ofv)_r/\GL_2(\ofv)} D_0(y)
=\sum_{y\in \GL_2(\ofv)\bs M_2(\ofv)_r/\GL_2(\ofv)} D(y).
$$
Note that the double coset in the last summation corresponds exactly to the classical Hecke correspondence $T(p_v^r)$. 
Hence, the above further equals 
$$
T(p_v^r) \hat\xi- \deg(T(p_v^r)) \hat\xi.
$$
Here $\deg(T(p_v^r))=\sigma_1(p_v^r)=1+N_v+\cdots+ N_v^r$. 

By \cite[Proposition 4.3.2]{Zh1}, 
$$
T(p_v^r) \overline\CL- \sigma_1(p_v^r)\ \overline\CL
= -2\sum_{i=0}^r i N_v^{r-i}\log N_v+ \log (N_v^{r\sigma_1(p_v^r)})
$$
In fact, the proposition considers a morphism 
$$
T(p_v^r) \CL\lra  \CL^{\otimes \deg(T(p_v^r))}
$$
and computes the norms of this morphism at non-archimedean places in part 1 and at archimedean places in part 2. Note that the result in part 2 of the proposition should be $N(m)^{\sigma_1(m)}$ instead of $N(m)^{2\sigma_1(m)}$.
The sum of the logarithms of these norms gives the formula.  

An elementary computation gives
$$
T(p_v^r) \overline\CL- \sigma_1(p_v^r)\ \overline\CL
= \frac{(r+2)N_v^{-(r+1)}-rN_v^{-(r+2)}-(r+2)N_v^{-1}+r}{(1-N_v^{-1})^2}
N_v^r\log N_v.
$$
The result follows by $\overline\CL=\kappa_U^\circ \cdot \hat\xi$. 
\end{proof}

The expression of $f_{\phi_v, a}(1,u)$ in the above lemma
happens to be very close to that of $W_{a,v}'(0,1,u)-\frac12\log|a|_v W_{a,v}(0,1,u)$ in  Lemma \ref{local explicit} (1). 
They will give great cancelation in our matching of the derivative series and the height series.

\section{Comparison of the two series} \label{sec comparison}

In this section, we will combine results in the last two sections to prove Theorem \ref{main}.
The upshot is to apply Lemma \ref{pseudo} to the difference
$$\mathcal D(g,\phi)=\pr I'(0, g, \phi)_U- 2 Z(g,(1,1))_U.  $$
Here we take $t_1=t_2=1$ for the CM points.
We refer to \S\ref{sec pseudo} for the notion of pseudo-Eisenstein series, and to \S\ref{sec key lemma} for a quick review of the notion of pseudo-theta series introduced in \cite[\S6.2]{YZ}. 

We will see that $\mathcal D(g,\phi)$ is a sum of finitely many non-singular pseudo-Eisenstein series and non-singular pseudo-theta series. 
Then Lemma \ref{pseudo} will imply that $\mathcal D(g,\phi)$ is the sum of the corresponding Eisenstein series and theta series. 
Since $\mathcal D(g,\phi)$ is cuspidal, its constant must be zero. 
This implies that the sum of the constant terms of the corresponding Eisenstein series and theta series is zero, which gives an equality involving the modular height of $X_U$. 
After computing all other terms, we get a formula of the modular height. 

To start with, let $(F,E,\BB, U,\phi)$ be as in \S \ref{choices}. We assume all the conditions of \S \ref{choices} throughout this section.
By Theorem \ref{analytic series}, 
$$
\pr I'(0, g,\phi)_U
=\pr' I'(0, g,\phi)_U-\pr'\CJ'(0, g, \phi)_U,
$$
By Theorem \ref{height series}, 
$$Z(g, (1, 1))_U
=\pair{Z_*(g,\phi)_U 1,  1}-\pair{Z_*(g,\phi)_U 1, \xi_{1}}
-\pair{Z_*(g,\phi)_U \xi_{1}, 1} +\pair{Z_*(g,\phi)_U\xi_{1}, \xi_{1}}.$$
Then the difference
\begin{eqnarray*}
\mathcal D(g,\phi)
&=& \pr' I'(0, g,\phi)_U-2\pair{Z_*(g,\phi)_U 1,  1}\\
&&-\pr'\CJ'(0, g, \phi)_U+ 2\pair{Z_*(g,\phi)_U 1, \xi_{1}}\\
&& +2\pair{Z_*(g,\phi)_U \xi_{1}, 1} -2\pair{Z_*(g,\phi)_U\xi_{1}, \xi_{1}}.
\end{eqnarray*}

In the following, for each of the three lines on the right-hand side of the above expression of $\mathcal D(g,\phi)$, we will describe the computational result, check that it is non-singular in the pseudo sense, and give its contribution in the equality after applying Lemma \ref{pseudo}. 

\subsubsection*{Third line}

Start with the third line, which has the simplest expression. By Proposition \ref{geometric series extra1}, 
$$\pair{Z_*(g,\phi)_U \xi_{1}, 1}=-\frac 12 \kappa_U^\circ\ E_*(0,g,\phi)_U\cdot \pair{\xi_{1}, 1},$$
$$\pair{Z_*(g,\phi)_U \xi_{1}, \xi_{1}}=-\frac 12 \kappa_U^\circ\ E_*(0,g,\phi)_U\cdot \pair{\xi_{1}, \xi_{1}}.$$
Here $\kappa_U^\circ$ denotes the degree of $L_U$ on a geometrically connected component of $X_U$.

The contribution of $2\pair{Z_*(g,\phi)_U \xi_{1}, 1} -2\pair{Z_*(g,\phi)_U\xi_{1}, \xi_{1}}$ after Lemma \ref{pseudo} is
\begin{equation}\label{line3}
\kappa_U^\circ\cdot (\pair{\xi_{1}, \xi_1}-\pair{\xi_{1}, 1})\cdot E(0,g,\phi)_U.
\end{equation}

\subsubsection*{Second line}
Now we consider the second line. Denote $c_3'=(1+\log 4)[F:\QQ]$. 
By Proposition \ref{analytic series extra},
$$
\pr'\CJ'(0, g, \phi)
=- (c_0+c_3') E_*(0,g,\phi)- \sum_{v\nmid\infty} C_*(0, g,\phi)(v)
+2\sum_{v\nmid\infty} E'(0,g,\phi)(v).$$
Here we have Eisenstein series
\begin{eqnarray*}
E(s, g, \phi) &=&\sumu \sum _{\gamma \in P(F)\bs \GL_2(F)} \delta(\gamma g)^s r(\gamma g)\phi(0,u),\\
C(s, g, \phi)(v) &=&\sum_{v\nmid\infty} \sumu \sum _{\gamma \in P(F)\bs \GL_2(F)} 
\delta(\gamma g)^s  c_{\phi_v}(\gamma g,0,u)\ r(\gamma g^v)\phi^v(0,u),
\end{eqnarray*}
with 
$$c_{\phi_v}(g,y,u)=r(g)\phi_{1,v}(y,u) W_{0,v}^{\circ}\, '(0,g,u,\phi_{2,v}) 
+ \log \delta(g_v)r(g)\phi_v(y,u);$$
and we have a pseudo-Eisenstein series 
$$ E'(0,g,\phi)(v)
=\sumu\sum_{a\in F^\times} W_a^v(0,g,u,\phi^v)
\left(W_{a,v}'(0,g,u,\phi_v)-\frac12\log|a|_v\cdot W_{a,v}(0,g,u,\phi_v)\right).$$

By Proposition \ref{geometric series extra2},
\begin{eqnarray*}
\pair{Z_*(g,\phi) 1, \xi_{1}}
= -\frac 12 \kappa_U^\circ\ E_*(0,g,\phi)_U\ \pair{1, \xi} 
 + \frac 12  \sum_{v\notin \Sigma} \CF^{(v)}_\phi(g),
\end{eqnarray*}
where the pseudo-Eisenstein series 
$$
\CF^{(v)}_\phi(g)
= \sumu\sum_{a\in F^\times}  W_{a}^v(0,g,u,\phi) \, 
f_{\phi_v, a}(g,u)$$
with
$$
f_{\phi_v, a}(g,u)=(1-N_v^{-2})|d_v|^{\frac32} |au^{-1}|_v\, \kappa_U^\circ
\sum_{y\in \bb_v(au^{-1})/U_v^1} r(g)\phi_v(y,u)  D_0(t_{1,v}^{-1}y t_{2,v}).
$$

The difference gives
\begin{eqnarray*}
&&-\pr'\CJ'(0, g, \phi)_U+ 2\pair{Z_*(g,\phi) 1, \xi_{1}}\\
&=&
(c_0+c_3'-\kappa_U^\circ \pair{1, \xi}) E_*(0,g,\phi)
+ \sum_{v\nmid\infty} C_*(0, g,\phi)(v) \\
&&
-2 \sum_{v\in \Sigma_f} E'(0,g,\phi)(v)
-2\sum_{v\notin \Sigma} (E'(0,g,\phi)(v)-\frac12\CF^{(v)}_\phi(g)).
\end{eqnarray*}
This is a finite sum of Eisenstein series and pseudo-Eisenstein series, by the following considerations using the explicit local results.
\begin{itemize}
\item[(1)] For any $v\in \Sigma_f$, the explicit result of 
Lemma \ref{local explicit}(2) implies that 
$$W_{a,v}'(0,1,u,\phi_v)-\frac12\log|a|_v\cdot W_{a,v}(0,1,u,\phi_v),$$ 
as a function of 
$(a,u)\in  F_v^\times\times  F_v^\times$, satisfies the condition of Lemma \ref{whittaker image2}.
Therefore, $E'(0,g,\phi)(v)$ is a non-singular pseudo-Eisenstein series in this case. Denote by 
$$E(0,g, \phi_v^+\otimes\phi^v)+E(0,g, \phi_v^-\otimes\phi^v)$$
the associated Eisenstein series. 
Note that Lemma \ref{whittaker image1} and 
Lemma \ref{local explicit}(2)
further give
$$
r(w)\phi_v^+(0,u)+r(w)\phi_v^-(0,u)=0, \quad \forall u\in F_v\cross.
$$

\item[(2)]
For any $v\notin \Sigma$, the pseudo-Eisenstein series
$$E'(0,g,\phi)(v)- \frac12\CF^{(v)}_\phi(g)
=\sumu\sum_{a\in F^\times} W_a^v(0,g,u,\phi^v)
\tilde f_{\phi_v, a}(g,u),$$
where 
$$
\tilde f_{\phi_v, a}(g,u)
=\left(W_{a,v}'(0,g,u,\phi_v)-\frac12\log|a|_v\cdot W_{a,v}(0,g,u,\phi_v)\right)
-\frac12 f_{\phi_v, a}(g,u)
$$
for any $a,u\in F_v^\times, \ g\in \gl(F_v)$.
By the explicit results of Lemma \ref{local explicit}(1) and 
Proposition \ref{geometric series extra3}, 
$\tilde f_{\phi_v, a}(1, u)\neq 0$ only if $u\in\ofv\cross$ and $v(a)\geq -v(d_v)$. 
Moreover, for $u\in \ofv\cross$ and $a\in \ofv$,
\begin{eqnarray*}
\tilde f_{\phi_v, a}(1,u)  
= \left(-\zeta_v'(2)/\zeta_v(2) + \log|d_v|\right) W_{a,v}(0,1,u)
+   |d_v|^{\frac{3}{2}}  
\frac{1-|d_v|}{N_v-1}\log N_v.
\end{eqnarray*}
By Lemma \ref{whittaker image2}, we see that 
$E'(0,g,\phi)(v)- \frac12\CF^{(v)}_\phi(g)$ is a non-singular pseudo-Eisenstein series in this case. 
The associated Eisenstein series is of the form
$$ \left(-\zeta_v'(2)/\zeta_v(2) + \log|d_v|\right) E(0,g, \phi)+E(0,g, \phi_v^+\otimes\phi^v)+E(0,g, \phi_v^-\otimes\phi^v).$$ 
The last two terms are 0 for almost all $v\notin\Sigma$. 
Moreover, Lemma \ref{whittaker image2} also gives for all $v\notin\Sigma$,
$$
\phi_v^+(0,u)+\phi_v^-(0,u)=0, \quad \forall u\in F_v\cross.
$$
\end{itemize}
Therefore, the contribution of $-\pr'\CJ'(0, g, \phi)_U+ 2\pair{Z_*(g,\phi) 1, \xi_{1}}$ after Lemma \ref{pseudo} is
\begin{eqnarray}\label{line2}
&&
(c_0+c_3'-\kappa_U^\circ \pair{1, \xi}+2\sum_{v\notin \Sigma} \left(\zeta_v'(2)/\zeta_v(2) - \log|d_v|\right)) E(0,g,\phi) \nonumber\\
&&+ \sum_{v\nmid\infty} C(0, g,\phi)(v) 
- 2 \sum_{v\nmid \infty} \big(E(0,g, \phi_v^+\otimes\phi^v)+E(0,g, \phi_v^-\otimes\phi^v)\big).
\end{eqnarray}

\subsubsection*{First line}

It remains to consider 
$$\pr' I'(0, g,\phi)_U-2\pair{Z_*(g,\phi)_U 1,  1}.$$ 
By Theorem \ref{analytic series}, the current $\pr' I'(0, g,\phi)_U$ has the same expression as the old $\pr I'(0, g,\phi)_U$ in \cite[Theorem 7.2]{YZ}.
By Theorem \ref{height series}, the current $\pair{Z_*(g,\phi)_U 1,  1}$ 
has the expression of the old $Z(g, (1,1),\phi))_U$ in 
\cite[Theorem 8.6]{YZ} with the extra term $\ds-\frac12 [F:\QQ]  E_*(0,g,\phi)$.
Consequently, the current difference
$\pr' I'(0, g,\phi)_U-2\pair{Z_*(g,\phi)_U 1,  1}$
has the expression of the old $\mathcal D(g,\phi)$ in 
\cite[\S9.1]{YZ} with an extra term
$[F:\QQ]  E_*(0,g,\phi)$.

Note that the choice of $\phi$ in \cite{YZ} is slightly different from what we have here. In \cite[\S7.2]{YZ}, it has an extra set $S_2$ of two non-archimedean places 
$v$ of $F$ split in $E$ with certain degenerate $\phi_v$. 
This assumption is made to kill the terms close to $E(s,g,\phi)$. 
However, the computations for $\mathcal D(g,\phi)$ holds for our 
$\pr' I'(0, g,\phi)_U-2\pair{Z_*(g,\phi)_U 1,  1}$
by pretending $S_2=\emptyset$.  

Hence, the translation of the computational results of \cite[\S9.1]{YZ} to the current setting gives
\begin{eqnarray*}
&& \pr' I'(0, g,\phi)_U-2\pair{Z_*(g,\phi)_U 1,  1}\\
&=& -2 \sum_{v\nmid\infty \ \nonsplit} \barint_{C_U}
(\CK^{(v)}_{\phi}(g,(t,t))-\CM^{(v)}_{\phi}(g,(t,t))\log N_v) dt\\
&&+\sum_{v\in \Sigma_f}
(2\log N_v) \barint_{C_U} j_{\bar v}(Z_*(g,\phi)_Ut,t) dt\\
&&+ \sum_{v\nmid\infty} \sumu \sum_{y\in E^\times} 
d_{\phi_v}(g,y,u)\ r(g)\phi^v(y,u)\\
&& +(\frac{2i_0(1,1)}{[O_E\cross:O_F\cross]} -c_1) \sumu \sum_{y\in E^\times} r(g)\phi(y,u)\\
&& +[F:\QQ]  E_*(0,g,\phi)\\
&& -[F:\QQ] (\gamma+\log(4\pi)-1) E_*(0,g,\phi).
\end{eqnarray*}
Here the last line comes from and is equal to  
$$
-2 \sum_{v|\infty} \barint_{C_U}
(\overline\CK^{(v)}_{\phi}(g,(t,t))-\CM^{(v)}_{\phi}(g,(t,t))) dt,
$$ 
which follows from Proposition \ref{archimedean comparison} and was missed in \cite[\S9.1]{YZ}.

The last expression is a sum of finitely many non-singular pseudo-theta series and non-singular pseudo-Eisenstein series. We refer to \S\ref{sec key lemma} for quick review of the notion of pseudo-theta series introduced in \cite[\S6.2]{YZ}. 
To explain why the fourth and the fifth line are pseudo-theta series, take the fourth line as an example. We have the series 
$$
\sumu \sum_{y\in E^\times} 
d_{\phi_v}(g,y,u)\ r(g)\phi^v(y,u). 
$$
Let $B$ be a totally definite quaternion algebra over $F$, with an embedding $E\to B$ and an isomorphism $B\otimes_F \AA^{S_0}\to \BB^{S_0}$ compatible with $E_\AA\to \BB$. 
Here $S_0$ is the set of places $v$ at which $B_v$ is not isomorphic to $\BB_v$. 
Then $\phi^{S_0\cup \{v\}}$ is also viewed as a Schwartz function on 
$B_{\AA^{S_0\cup \{v\}}}\times (\AA^{S_0\cup \{v\}})^\times$. 
It follows that we can view the series as a pseudo-theta series sitting on quadratic spaces 
$0\subset E\subset B$. 

Note that neither the fourth line or the fifth line contributes to the result after Lemma \ref{pseudo} because they are \textit{degenerate} pseudo theta series. 
So we only recall the other lines. 
Recall that 
\begin{eqnarray*}
&&\CK^{(v)}_\phi(g,(t_1,t_2))-\CM^{(v)}_\phi(g,(t_1,t_2))(\log N_v)\\
&=&\sum_{u\in \mu_U^2\bs F\cross}\sum _{y\in B(v)-E}  r(g,(t_1,t_2))\phi^v (y, u)
\bar k_{r(t_1,t_2)\phi _v}(g,y, u)
\end{eqnarray*}
where 
$$
\bar k_{\phi_v}(g,y,u)= k_{\phi_v}(g,y,u)-m_{r(g)\phi_v}(y,u)\log N_v.
$$
In particular.
$$
\bar k_{\phi_v}(y,u)= k_{\phi_v}(1,y,u)-m_{\phi_v}(y,u)\log N_v
$$
extends to a Schwartz function in $\overline\CS(B(v)_v\times F_v\cross)$.
Similarly, 
\begin{eqnarray*}
j_{\bar v}(Z_*(g,\phi)_Ut_1,t_2)
&=&\sum_{u\in \mu_U^2\bs F\cross}\sum _{y\in B(v)^\times}  r(g,(t_1,t_2))\phi^v (y, u)
 r(t_1,t_2)l_{r(g)\phi _v}(y, u)
\end{eqnarray*}
where $l_{\phi_v}(y,u)$ extends to a Schwartz function in $\overline\CS(B(v)_v\times F_v\cross)$.

The contribution of $\pr' I'(0, g,\phi)_U-2\pair{Z_*(g,\phi)_U 1,  1}$ after Lemma \ref{pseudo} is
\begin{eqnarray}
&&  -2 \sum_{v\nmid\infty \ \nonsplit} \barint_{C_U}
\theta(g,(t,t), \bar k_{\phi_v}\otimes\phi^v) dt \nonumber\\
&& +\sum_{v\in \Sigma_f}
(2\log N_v) \barint_{C_U}
\theta(g,(t,t), l_{\phi_v}\otimes\phi^v) dt \nonumber\\
&& 
-[F:\QQ] (\gamma+\log(4\pi)-2) E_*(0,g,\phi). \label{line1}
\end{eqnarray}

\subsubsection*{The sum}

As a conclusion, the difference $\mathcal D(g,\phi)$ is the sum of finitely many 
non-singular pseudo-Eisenstein series and finitely many non-singular pseudo-theta series. We are finally ready to apply Lemma \ref{pseudo}.
For the conclusion, the sum of (\ref{line3}), (\ref{line2}) and (\ref{line1}) gives
\begin{eqnarray}
\mathcal D(g,\phi)
&=&  -2 \sum_{v\nmid\infty \ \nonsplit} \barint_{C_U}
\theta(g,(t,t), \bar k_{\phi_v}\otimes\phi^v) dt
+\sum_{v\in \Sigma_f} (2\log N_v)
\barint_{C_U} \theta(g,(t,t), l_{\phi_v}\otimes\phi^v) dt \nonumber\\
&&+ \sum_{v\nmid\infty} C(0, g,\phi)(v) 
- 2 \sum_{v\nmid\infty} \big(E(0,g, \phi_v^+\otimes\phi^v)+E(0,g, \phi_v^-\otimes\phi^v)\big) \nonumber\\
&&+c_4\cdot E(0,g,\phi)_U, \label{final D}
\end{eqnarray}
where
$$c_4=c_0+c_3'-2\kappa_U^\circ \pair{1, \xi}+\kappa_U^\circ \pair{\xi_1, \xi_1}
+2\sum_{v\notin \Sigma} \left(\zeta_v'(2)/\zeta_v(2) - \log|d_v|\right)
-[F:\QQ] (\gamma+\log(4\pi)-2).$$
Here we have used the identity $ \pair{1, \xi}= \pair{1, \xi_1}$, which holds by considering geometrically connected components. 
Moreover, we have the following result. 
\begin{lem}
$$
\kappa_U^\circ \pair{\xi_1, \xi_1}=-2\, h_{\ol \CL_U}(X_U), \qquad
\kappa_U^\circ \pair{1, \xi}=- h_{\ol\CL_U}(P_U)+[F:\QQ].
$$
\end{lem}
\begin{proof}
Denote $|\pi_0|$ the number of geometrically connected component of $X_U$. 
The first result goes as follows:
$$
\kappa_U^\circ \pair{\xi_1, \xi_1}
= \frac{\kappa_U^\circ}{|\pi_0|}\pair{\xi, \xi}
= \frac{\kappa_U^\circ}{|\pi_0|} \frac{1}{(\kappa_U^\circ)^2} \pair{\ol\CL_U, \ol\CL_U}
= \frac{1}{\deg(L_U)} \pair{\ol\CL_U, \ol\CL_U}
=- \frac{1}{\deg(L_U)} \wh\deg(\hat c_1(\ol\CL_U)^2).
$$
The first equality holds by considering geometrically connected components of $X_U$, and the other equalities holds by definition. 
Note that the negative sign of the last equality is due to different normalizations of the intersection numbers, which is originally from the negative sign in the arithmetic Hodge index theorem. 

The second result goes as follows:
$$
\kappa_U^\circ \pair{1, \xi}
=  \pair{1, \ol\CL_U}
=  \pair{\bar P_U, \ol\CL_U}+[F:\QQ]
= - h_{\ol\CL_U}(P_U)+[F:\QQ].
$$
Here the second equality follows from Lemma \ref{int of Green}.
\end{proof}

By the lemma, $c_4$ contains the height $h_{\ol \CL_U}(X_U)$ which we need to compute. For $h_{\ol\CL_U}(P_U)$, by Theorem \ref{height CM}, the main result of Part II of \cite{YZ},  
$$ h_{\ol\CL_U}(P_U)
=-\frac{L_f'(0,\eta)}{L_f(0,\eta)} +\frac 12 \log  \frac{d_\BB}{d_{E/F}}.
$$
It cancels the major part of 
$$c_0=2\frac{L'(0,\eta)}{L(0,\eta)} +\log|d_E/d_F|
=2\frac{L_f'(0,\eta)}{L_f(0,\eta)} +\log|d_E/d_F| - [F:\BQ](\gamma+\log 4\pi ).$$
More precisely, 
$$
c_0+2\,h_{\ol\CL_U}(P_U)=\log|d_\BB d_F| - [F:\BQ](\gamma+\log 4\pi ).
$$
By $c_3'=(1+\log 4)[F:\QQ]$, we further have 
$$
c_0+c_3'+2\,h_{\ol\CL_U}(P_U)=\log|d_\BB d_F| - [F:\BQ](\gamma+\log \pi-1).
$$
Hence, we can simplify $c_4$ to get 
$$
c_4=-2\, h_{\ol \CL_U}(X_U)
+2\sum_{v\notin \Sigma} \left(\zeta_v'(2)/\zeta_v(2) - \log|d_v|\right)
+\log|d_\BB d_F| 
-[F:\QQ] (2\gamma+2\log(2\pi)-1).
$$

\subsubsection*{The constant terms}

Note that $\mathcal D(g,\phi)$ is a cusp form, so its constant term must be 0. 
Then the constant terms of the right-hand side of (\ref{final D}) should  be 0. 
This will give the result we need.

In the following, we first treat the case $|\Sigma|>1$ and then mention the difference for the easier case $|\Sigma|=1$.
While it is straightforward to write down the constant terms of the theta series, it takes a little extra effort to treat those for the Eisenstein series. We claim that the constant terms of the Eisenstein series 
$$
   E(0,g,\phi)_U,\qquad C(0, g,\phi)(v), \qquad
  E(0,g, \phi_v^+\otimes\phi^v)+E(0,g, \phi_v^-\otimes\phi^v), 
$$
 are respectively equal to 
\begin{eqnarray*}
&& \sumu r(g)\phi(0,u), \\
&& \sumu c_{\phi_v}(g,0,u)\ r(g^v)\phi^v(0,u),\\
&& \sumu 
(r(g_v)\phi_v^+(0,u)+r(g_v)\phi_v^-(0,u))\ r(g^v)\phi^v(0,u). 
\end{eqnarray*}
In other words, the contribution from the intertwining part at $s=0$ is 0.

The claim is a consequence of the assumption $|\Sigma|>1$.
In fact, the result for $E(0,g,\phi)_U$ is immediately a consequence of 
\cite[Proposition 2.9(3)]{YZZ}. 
The other two Eisenstein series are similar, and we take $C(0, g,\phi)(v)$ for example. 
Recall from \S\ref{sec 3.1} that 
$$C(s, g, \phi)(v)= \sumu C(s, g, u,\phi)(v)$$
with 
$$C(s, g, u, \phi)(v)=  \sum _{\gamma \in P^1(F)\bs \SL_2(F)} 
\delta(\gamma g)^s  \Psi(\gamma g, u),$$
where
$$\Psi(g,u)
=c_{\phi_v}(g_v,0,u)\ r(g^v)\phi^v(0,u)
$$
is a principal series in the sense that
$$
\Psi(m(a)n(b)g,u)=|a|_\AA^{2}\Psi(g,u), \quad a\in \AA^\times, \ b\in \AA.
$$
The constant term  
$$
C_0(s, g, u, \phi)(v)=\delta(g)^s  \Psi(g, u) + W_0(s,g,u,\Psi)
$$
with the intertwining part
$$
W_0(s,g,u,\Psi)=\int_\AA \delta(wn(b) g)^s  \Psi(wn(b) g, u) db.
$$
Note that $\Psi=\otimes_w \Psi_w$ is naturally a product of local terms, and we can define $W_{0,w}(s,g,u,\Psi_w)$ similarly. 
Following \cite[\S2.5.2]{YZZ}, 
for the sake of analytic continuation at $s=0$, we write 
$$
W_0(s,g,u,\Psi)=\tilde\zeta_F(s+1)\prod_w W_{0,w}^\circ(s,g,u,\Psi_w),
$$
where the normalized term
$$
W_{0,w}^\circ(s,g,u,\Psi_w)=\frac{1}{\zeta_w(s+1)}W_{0,w}(s,g,u,\Psi_w)
$$
is holomorphic at $s=0$ for any $w$. 
By \cite[Proposition 2.9(1)(a)]{YZZ},  
$$
W_{0,w}^\circ(0,g,u,\Psi_w)=0, \quad w\in \Sigma\setminus\{v\}.$$ 
In other words, $W_{0,w}^\circ(s,g,u,\Psi_w)$ has a zero at $s=0$ for any $w\in \Sigma\setminus\{v\}$. 
On the other hand, there is simple pole of $\zeta_F(s+1)$ at $s=0$.
By the above product expression,  $W_0(s,g,u,\Psi)$ has a zero at $s=0$ of order at least 
$|\Sigma\setminus\{v\}|-1\geq 1$ by $|\Sigma|\geq 3$. 
It follows that $W_0(0,g,u,\Psi)=0$.
This proves the claim. 

Taking the constant terms of (\ref{final D}),
we end up with 
\begin{eqnarray*}
0
&=&  -2 \sum_{v\nmid\infty \ \nonsplit} \sumu
r(g)(\bar k_{\phi_v}\otimes\phi^v)(0,u) \\
&&+\sum_{v\in \Sigma_f} (2\log N_v) \sumu
r(g)(l_{\phi_v}\otimes\phi^v)(0,u)\\
&&+ \sum_{v\nmid\infty} \sumu c_{\phi_v}(g,0,u)\ r(g^v)\phi^v(0,u)\\
&&- 2 \sum_{v\nmid\infty} \sumu 
(r(g_v)\phi_v^+(0,u)+r(g_v)\phi_v^-(0,u))\ r(g^v)\phi^v(0,u)\\
&&+c_4 \sumu r(g)\phi(0,u).
\end{eqnarray*}
The goal is to get a formula of $c_4$ from the expression. 
Then it suffices to take a specific $g\in \gla$ such that 
$$\sumu r(g)\phi(0,u)\neq 0.$$
Note that $g=1$ does not work since $\phi_v(0,u)=0$ for any $v\in \Sigma_f$. 
Define $g=(g_v)_v\in\gla$ by 
$$
g_v=\begin{cases}
w=\matrixx{}{1}{-1}{}, & v\in \Sigma_f, \\
1, & v\notin \Sigma_f. 
\end{cases}
$$ 
Now we simplify the above equality for this $g$. 

By the above discussion, we already have 
$$r(g_v)\phi_v^+(0,u)+r(g_v)\phi_v^-(0,u)=0$$
for any $v\nmid\infty$. 
It follows the fourth line of the right-hand side of the equality is 0. 
The equation becomes 
\begin{eqnarray}
0
&=&  -2 \sum_{v\nmid\infty \ \nonsplit} \sumu
r(g)(\bar k_{\phi_v}\otimes\phi^v)(0,u) \nonumber\\
&&+\sum_{v\in \Sigma_f} (2\log N_v) \sumu
r(g)(l_{\phi_v}\otimes\phi^v)(0,u) \nonumber\\
&&+ \sum_{v\nmid\infty} \sumu c_{\phi_v}(g,0,u)\ r(g^v)\phi^v(0,u) \nonumber\\
&&+c_4 \sumu r(g)\phi(0,u).  \label{constant term}
\end{eqnarray}
Here we recall that
$$
\bar k_{\phi_v}(y,u)= k_{\phi_v}(1,y,u)-m_{\phi_v}(y,u)\log N_v
$$
extends to a Schwartz function in $\overline\CS(B(v)_v\times F_v\cross)$.

Note that each of the first three lines of the right-hand side of (\ref{constant term}) has a sum over places certain places $v$ of $F$.
In the following, for each non-archimedean place $v$, we consider the contribution of this fixed $v$ from our these three lines. 
\begin{itemize}
\item[(1)] 
If $v$ is split in $E$, then only the third line has contribution from $v$. In this case, by \cite[Lemma 7.6]{YZ},
$$
c_{\phi_v}(1,0,u)
=\log |d_v|\ \phi_v(0,u).
$$
 
\item[(2)] If $v$ is nonsplit in $E$ but split in $\BB$, then both the first line  and the the third line has contribution from $v$. In this case, by \cite[Lemma 7.4, Lemma 7.6, Lemma 8.7]{YZ},
$$
-2 \bar k_{\phi_v}(0,u)+c_{\phi_v}(1,0,u)=\log |d_v|\ \phi_v(0,u).
$$
See also \cite[Proposition 9.2]{YZ}. 

\item[(3)] If $v$ is nonsplit in $\BB$, 
by Lemma \ref{derivative of intertwining}, Lemma \ref{average k}, and 
Proposition \ref{averaged j}, 
$$
-2 r(w)\bar k_{\phi_v}(0,u)+2r(w)l_{\phi_v}(0,u)\log N_v+c_{\phi_v}(w,0,u)=(-\log |d_v|+\alpha_v\log N_v)
r(w)\phi_v(0,u),$$
where 
$$
\alpha_v=
1-\frac{N_v-1}{2(N_v+1)}.
$$
Note that the expressions in Lemma \ref{derivative of intertwining} and Lemma \ref{average k} depend on the parity of $v(d_v)$, but their combined expression for $\alpha_v$
happens to be uniform for all $v(d_v)$. 
This can be explained by the fact that Lemma \ref{derivative of intertwining} treats 
$W_{0,v}'(0,g,u, \phi_{2,v})$, while Lemma \ref{average k}
treats $W_{a,v}'(0,g,u, \phi_{2,v})$. 
\end{itemize}

Taking all these into consideration, the equation becomes
\begin{eqnarray*}
0
&=&  \left(\sum_{v\notin\Sigma} \log |d_v|+
\sum_{v\in\Sigma_f}  (-\log |d_v|+\alpha_v\log N_v) +c_4 \right)
\sumu r(g)\phi(0,u).
\end{eqnarray*}
Note that 
$$
\sumu r(g)\phi(0,u)>0
$$
by our choice of $g$.
We get an equation  
\begin{eqnarray}
\sum_{v\notin\Sigma} \log |d_v|+
\sum_{v\in\Sigma_f}  (-\log |d_v|+\alpha_v\log N_v) +c_4=0.\label{c4}
\end{eqnarray}
This is obtained for the case $|\Sigma|>1$.

If $|\Sigma|=1$, we claim that (\ref{c4}) also holds.
In this case, the constant terms for $E(0,g,\phi)_U$ and other similar series might contain nonzero intertwining parts by  
\cite[Proposition 2.9(3)]{YZZ}. 
We may figure our the effect of this by extra argument.
Alternatively, (\ref{final D}) simply implies
$$\mathcal D(g,\phi)
=c_4\cdot E(0,g,\phi)_U,
$$
since the other terms are zero by the computational results.
Comparing the constant terms, we easily have $c_4=0$, since $\mathcal D(g,\phi)$ is cuspidal.
This agrees with (\ref{c4}).

\subsubsection*{Logarithmic derivative}
Recall that 
$$
c_4=-2\, h_{\ol \CL_U}(X_U)
+2\sum_{v\notin \Sigma} \left(\zeta_v'(2)/\zeta_v(2) - \log|d_v|\right)
+\log|d_\BB d_F| -[F:\QQ] (2\gamma+2\log(2\pi)-1)
$$
Then (\ref{c4}) becomes 
\begin{eqnarray*}
-2\, h_{\ol \CL_U}(X_U)
+2\sum_{v\notin \Sigma} \zeta_v'(2)/\zeta_v(2)
+\log|d_\BB d_F^2| -[F:\QQ] (2\gamma+2\log(2\pi)-1)
+\sum_{v\in\Sigma_f}  \alpha_v\log N_v
=0.
\end{eqnarray*}

Note that  
$$
\frac{\zeta_{F}'(2)}{\zeta_{F}(2)}
=\sum_{v\nmid\infty} \frac{\zeta_v'(2)}{\zeta_v(2)}, \qquad
\frac{\zeta_v'(2)}{\zeta_v(2)}
=-\frac{N_v^{-2}}{1-N_v^{-2}}\log N_v.$$
The first equality holds because the Euler product of $\zeta_{F}(s)$ is absolutely convergent for $\Re(s)>1$. 
Hence, we finally end up with
\begin{eqnarray*}
-2\, h_{\ol \CL_U}(X_U)
+2\frac{\zeta_{F}'(2)}{\zeta_{F}(2)}+
\sum_{v\in\Sigma_f}  \left(\alpha_v+\frac{2N_v^{-2}}{1-N_v^{-2}}+1\right)\log N_v
+\log|d_F^2| -[F:\QQ] (2\gamma+2\log(2\pi)-1)=0.
\end{eqnarray*}
Here the local term
$$
\alpha_v+\frac{2N_v^{-2}}{1-N_v^{-2}}+1
=1-\frac{N_v-1}{2(N_v+1)}+\frac{2N_v^{-2}}{1-N_v^{-2}}+1
=\frac32+\frac{1}{N_v-1}
=\frac{3N_v-1}{2(N_v-1)}.
$$
Therefore,
\begin{eqnarray*}
h_{\ol \CL_U}(X_U)
=\frac{\zeta_{F}'(2)}{\zeta_{F}(2)}+
\sum_{v\in\Sigma_f}  \frac{3N_v-1}{4(N_v-1)}\log N_v
+\log|d_F|-(\gamma+\log(2\pi)-\frac12)[F:\QQ] .
\end{eqnarray*}

\subsubsection*{Functional equation}

We can convert the logarithmic derivative at $2$ to that at $-1$ by the functional equation. In fact, the completed Dedekind zeta function 
$$
\tilde\zeta_{F}(s)=\tilde\zeta_{F,\infty}(s)\zeta_{F}(s)$$
with the gamma factor
$$
\tilde\zeta_{F,\infty}(s)=(\pi^{-s/2}\Gamma(s/2))^{[F:\BQ]}
$$
has functional equation 
$$
\tilde\zeta_{F}(1-s)=|d_K|^{s-\frac12} \tilde\zeta_{F}(s).
$$
Note that
\begin{align*}
\frac{\tilde\zeta_{F,\infty}'(2)}{\tilde\zeta_{F,\infty}(2)}
=& - \frac{1}{2}(\gamma+\log \pi )[F:\QQ],\\
\frac{\tilde\zeta_{F,\infty}'(-1)}{\tilde\zeta_{F,\infty}(-1)}
=& - \frac{1}{2}(\gamma+\log (4\pi) )[F:\QQ]+[F:\QQ].
\end{align*}
It follows that
\begin{eqnarray*}
h_{\ol \CL_U}(X_U)
&=&\frac{\tilde\zeta_{F}'(2)}{\tilde\zeta_{F}(2)}
+\log|d_F| -\frac12[F:\QQ] (\gamma+\log(4\pi)-1)
+\sum_{v\in\Sigma_f}  \frac{3N_v-1}{4(N_v-1)}\log N_v\\
&=&-\frac{\tilde\zeta_{F}'(-1)}{\tilde\zeta_{F}(-1)} -\frac12[F:\QQ] (\gamma+\log(4\pi)-1)
+\sum_{v\in\Sigma_f}  \frac{3N_v-1}{4(N_v-1)}\log N_v\\
&=&-\frac{\zeta_{F}'(-1)}{\zeta_{F}(-1)} -\frac12[F:\QQ] 
+\sum_{v\in\Sigma_f}  \frac{3N_v-1}{4(N_v-1)}\log N_v.
\end{eqnarray*}
This prove Theorem \ref{main}.

\

\noindent \small{Beijing International Center for Mathematical Research, Peking University, Beijing 100871, China}

\noindent \small{\it Email: yxy@bicmr.pku.edu.cn}

\end{document}